\documentclass[11pt,preprint,authoryear]{elsarticle}

\usepackage{pgfplots}
\usepackage{lipsum}
\usepackage{amssymb}
\usepackage{amsthm}
\usepackage{lineno}
\usepackage{pgfplots}
\usepackage{comment}
\usepackage{geometry}
\usepackage{caption}
\usepackage{subcaption}

\pgfkeys{/pgfplots/MyAxisStyle/.style={xmin=0,xmax=5, ymin=0,ymax=3,height=4cm,width=0.9\linewidth,scale only axis}}
\pgfkeys{/pgfplots/MyLineStyle/.style={samples=50, smooth, ultra thick}}

\usepackage{natbib}
\setcitestyle{authoryear}
\biboptions{round,semicolon}

\usepackage{pgfplotstable}
\usepackage{pgfplots}
\usepackage{amsmath}
\usepackage{stackrel}
\usepackage{mathtools}
\usepackage{glossaries}
\usepackage{graphicx}
\usepackage{algcompatible}
\usepackage{float}
\usepackage{url}
\usepackage{tikz}
\usepackage{pgfplots}
\usetikzlibrary{pgfplots.statistics,intersections,positioning,calc}
\pgfplotsset{width=7.5cm}
\usetikzlibrary{patterns}
\usetikzlibrary{plotmarks}
 \usepackage{setspace}
\pgfplotsset{compat=1.16}

\newcommand{\R}{\mathbb{R}}
\newcommand{\X}{\mathcal{X}} 
\newcommand{\Z}{\mathcal{Z}} 
\newcommand{\ratio}{\sigma} 

\newcommand{\D}{\mathcal{D}} 
\newcommand{\A}{\mathcal{A}} 

\usepackage{array,calc}
\usepackage[linesnumbered,lined,boxed,commentsnumbered]{algorithm2e}
\newtheorem{definition}{Definition}
\newtheorem{remark}{Remark}
\newtheorem{example}{Example}
\newtheorem{proposition}{Proposition}
\newtheorem{lemma}{Lemma}
\usepackage{algcompatible}
\usepackage{amsthm}
\usepackage{amsmath}
\usepackage{lipsum}
\usepackage{pgfplots}
\usepackage{filecontents}
\usepackage{xcolor}

\usepackage{setspace}
\newacronym{LP}{LP}{Linear Programme}
\newacronym{BLP}{BLP}{Bi-objective Linear Programme}
\newacronym{BBSA}{BBSA}{Bi-objective Benders Simplex Algorithm}
\newacronym{FCTP}{FCTP}{Fixed Charge Transportation Problem}
\newacronym{BFCTP}{Bi-FCTP}{Bi-objective Fixed Charge Transportation Problem}
\newacronym{BM}{BM}{Bi-objective Master Problem}
\newacronym{BS}{BS}{Bi-objective Subproblem}
\newacronym{PP}{PP}{Portfolio Problem}
\newacronym{BPP}{Bi-PP}{Bi-objective Portfolio Problem}
\newacronym{MKP}{MKP}{Multiple Knapsack Problem}
\newacronym{BMKP}{Bi-MKP}{Bi-objective Multiple Knapsack Problem}
\newacronym{DBA}{DBA}{Dichotomic Benders Algorithm}

\usepackage{csvsimple}
\begin{document}

\begin{frontmatter}

\title{Benders Decomposition for Bi-objective Linear Programs}

\author[First]{Andrea Raith\corref{cor1}}
\cortext[cor1]{Corresponding author a.raith@auckland.ac.nz}
\author[Second]{Richard Lusby}
\author[First]{Ali Akbar Sohrabi Yousefkhan}

\affiliation[First]{
organization={Department of Engineering Science, The University of Auckland},
city={Auckland},
country={New Zealand}
}
\affiliation[Second]{organization={Department of Technology, Management and Economics, Operations Research Section, Danish Technical University},
city={Lyngby},
country={Denmark}
}

\begin{abstract}
In this paper, we develop a new decomposition technique for solving bi-objective linear programming problems. The proposed methodology combines the bi-objective simplex algorithm with Benders decomposition and can be used to obtain a complete set of efficient extreme solutions, and the corresponding set of extreme non-dominated points, for a bi-objective linear program. Using a Benders-like reformulation, the decomposition approach decouples the problem into a bi-objective master problem and a bi-objective subproblem, each of which is solved using the bi-objective parametric simplex algorithm. The master problem provides candidate efficient solutions that the subproblem assesses for feasibility and optimality. As in standard Benders decomposition, optimality and feasibility cuts are generated by the subproblem and guide the master problem solve. This paper discusses bi-objective Benders decomposition from a theoretical perspective, proves the correctness of the proposed reformulation and addresses the need for so-called weighted optimality cuts. Furthermore, we present an algorithm to solve the reformulation
and discuss its performance for three types of bi-objective optimisation problems.
\end{abstract}

\end{frontmatter}

{\bf Declaration of Interest:} This research has been partially supported by the Marsden Fund, grant number 16-UOA-086. The funder had no involvement in the research itself, in writing this article or in the decision to submit it to this journal.

\section{Introduction}
\label{sec:introduction}
A wide variety of problems in different fields such as engineering and industry involve several conflicting objectives that should be taken into account simultaneously in optimisation problems. It is generally impossible to find an optimal solution of such a multi-objective problem that optimises all the objectives concurrently.
Many real-world decision making problems not only have multiple objectives,  they are also modelled by complex mathematical programs that are often large-scale.  Different methods can be employed to solve these problems. In particular, the interest in decomposition methods to solve large-scale linear optimisation problems has been growing remarkably due to their many applications.  Benders decomposition is a widely used decomposition technique
that was first introduced by \citet{benders1962} to solve problems with a dual block-angular structure. This approach has been used on a broad range of applications such as production planning \citep{adulyasak2015}, facility location problems \citep{boland2016}, airline scheduling \citep{cordeau2001}, capacity expansion problems \citep{bloom1983}, investment planning \citep{oliveira2014}, and healthcare planning \citep{luong2015}. In this paper, we combine two areas of research to describe a novel solution approach for \glspl{BLP} that integrates a Benders decomposition approach and the bi-objective parametric simplex algorithm.

There is a scarcity of literature on the application of problem decomposition approaches in the context of bi-objective or multi-objective optimisation problems. Most of the published literature in this context applies a scalarisation to the multi-objective problem, and the resulting single-objective problem is then solved using conventional problem decomposition techniques.
A lot of research that applies decomposition techniques does in fact solve problems with objectives made up of a weighted-sum of multiple objectives, which could be considered scalarisations of multi-objective problems. However, the multiple objectives are usually not a focus of the articles
\citep[e.g.][]{amberg2019a,cacchiani2017a,lusby2017a}.

Scalarisation-based methods that use column generation for bi-objective or multi-objective problems are covered first.
Several authors solve a sequence of one or more problems having a weighted-sum objective with a column generation approach: \citet{SU13} (radiotherapy treatment planning), \citet{franco2016} (inventory routing), and \citet{GJN18} (vehicle routing). An $\varepsilon$-constraint scalarisation is applied in combination with column generation by the following authors: \citet{ET07} (integer programs), \citet{TERZ11} (airline crew scheduling), \citet{SAJ13} (vehicle routing), \citet{AJS18} (combinatorial problems) and \citet{glize2020} (multi-vehicle covering tour).

Other scalarisation-based methods use Benders decomposition to solve bi-objective or multi-objective problems.
Weighted-sum scalarisation is applied, and the resulting problems are solved with Benders decomposition in \citet{mardan} (supply chain network design)
and \citet{PTG20} (stochastic facility location). 
In the context of bi-objective multi-stage stochastic programming, \citet{dowson2022} describe a bi-objective Benders decomposition approach that iteratively performs iterations of a weighted-sum, single-objective Benders approach. The weighting that is used in both the master and the subproblem is updated iteratively akin to the parametric simplex method, thus identifying efficient solutions of the current master problem. It should be noted that several sweeps of the full efficient set (i.e.~across all possible weights) are required until there is no change in the set of non-dominated points. It is also assumed that the subproblem is always feasible.
Benders decomposition is used to solve $\varepsilon$-constraint problems in: \citet{KA93} (electrical power distribution), \citet{attari2018} (stochastic optimisation model for blood collection and testing), and \citet{momeni2022} (coordinating trucks and drones to monitor forests for bushfires).

Some researchers have integrated the concept of column generation into the steps of the bi-objective parametric simplex method. \citet{RMES12} describe how the entering variable selection in the bi-objective parametric simplex method can be reformulated as a column generation subproblem for \glspl{BLP}.
\citet{MRE15} extend this research to perform Dantzig-Wolfe decomposition by considering multiple subproblems, one for each commodity of a bi-objective multi-commodity network flow problem.

In this paper, we provide an introduction to bi-objective linear optimisation problems and the Benders decomposition technique in Sections~\ref{sec:BLP} and~\ref{subsec:single-obj-benders}.
We present an approach that performs Benders decomposition for \glspl{BLP} within the iterations of the bi-objective parametric simplex algorithm.
Given a \gls{BLP}, we first show that we can always find a bi-objective Benders reformulation
consisting of only a subset of variables, two new variables to capture the first and second objective values of the remaining problem and cuts to ensure feasibility of the master variables as well as optimal objective values. This reformulation has a complete set of efficient solutions that is equivalent to that of the original \gls{BLP}, see Section~\ref{sec:full-Benders-reformulation-BLP}. Benders master and subproblems are formulated and
an algorithm is developed that integrates the concept of Benders decomposition within the steps of the bi-objective parametric simplex method in Section~\ref{sec:algorithms}. Section~\ref{sec:numerical} summarises the results of an extensive computational study. Section~\ref{sec:conclusion} concludes the paper and identifies promising directions for future work.

This paper makes the following contributions to the literature. Firstly, we show how a full Benders reformulation can be derived in the bi-objective case, i.e.~we introduce a bi-objective Benders reformulation of a non-decomposed \gls{BLP} consisting of both feasibility cuts and what we term weighted optimality cuts. Secondly, we devise a bi-objective Benders decomposition approach and prove the correctness of the reformulation by showing that it can be used to obtain a complete set of efficient extreme solutions to the original problem. Finally, we present an algorithm, which we term the \gls{BBSA}, for solving the bi-objective Benders reformulation. This algorithm efficiently combines the bi-objective parametric simplex algorithm with Benders decomposition and enables solving the problem in a single sweep (with some backtracking).
The \gls{BBSA} is applied to solve three different types of \gls{BLP} and extensive numerical results are discussed, including a comparison with a dichotomic approach that solves weighted single-objective Benders problems.

\section{Benders Decomposition for Bi-objective Linear Optimisation Problems}
\label{sec:theory}

In this section, we consider Benders decomposition in a bi-objective setting. To provide the necessary background, we begin in Section~\ref{sec:BLP} by introducing a generic \gls{BLP} and some relevant notation and terminology.
In Section~\ref{subsec:single-obj-benders}, we review how Benders decomposition can be used to solve a single-objective \gls{LP}. Section~\ref{sec:full-Benders-reformulation-BLP} then extends the single-objective methodology to the bi-objective case. In particular, we introduce the bi-objective Benders reformulation of a \gls{BLP} supported by an illustrative example.

\subsection{Bi-objective Linear Programs}
\label{sec:BLP}

\label{subsec:formulation_defs}
A \gls{BLP} in standard form is given by
\begin{equation}
\label{eq:BLPstandard}
\begin{array}{rrcl}
   \min & \left ( \begin{array}{c} (\hat{c}^1)^{\top} x  \\ (\hat{c}^2)^{\top} x \end{array} \right ) & = & \left ( \begin{array}{c} z_1(x) \\ z_2(x) \end{array} \right )\\
   \text{s.t.} & \hat{A}x  &=& \hat{b}, \\
               &  x &\geqq& 0,
\end{array}
\end{equation}
where $x\in \R^n$ and $\hat{c}^1,\hat{c}^2,\hat{A},\hat{b}$ have appropriate dimensions.
In the following we define some terms related to solutions of \glspl{BLP}.
\begin{definition}
For a \gls{BLP} we call $\R^n$ the \emph{decision space} and $\R^2$ the \emph{objective space}.
The set of feasible solutions in decision space $x\in \R^n$ that satisfy the constraints of \gls{BLP} is called the \emph{feasible set in decision space} $\X \subseteq \R^n$, whereas its image is called the \emph{feasible set in objective space} $\Z \subseteq \R^2$. For~\eqref{eq:BLPstandard} the feasible set is $\X = \{ x \in \R^n | \hat{A}x = \hat{b}, x \geqq 0 \}$ and its image in objective space is
$\Z = \left \{ z(x) \in \R^2 | z(x) = (z_1(x), z_2(x))^{\top} \text{ and } x \in \X \right \}$.
\end{definition}

\begin{definition}
When comparing objective vectors $z,\bar{z} \in \R^2$ the following notation is used. We define symbols $<, \leq$, and $\leqq$ as follows: $z < \bar{z}$ if $z_i < \bar{z}_i$ for $i=1,2$; $z \leqq \bar{z}$ if $z_i \leqq \bar{z}_i$ for $i=1,2$; $z \leq \bar{z}$ if $z_i \leqq \bar{z}_i$ for $i=1,2$ and $z \neq \bar{z}$. When $z \leq \bar{z}$, $z$ \emph{dominates} $\bar{z}$, and when $z < \bar{z}$, $z$ \emph{strictly dominates} $\bar{z}$.

\end{definition}
\begin{definition}
When solving a bi-objective optimisation problem we seek to identify the set of \emph{efficient} solutions $\X_e$ that consists of all solutions $x \in \X$ whose objective vector $z(x)$ is not dominated by that of another solution $\bar{x} \in \X$.
The corresponding set $\Z_n$ in objective space is called the set of \emph{non-dominated} points $\Z_n = \{ z(x) \in \R^2 | x \in \X_e \}$. Similarly, a \emph{weakly efficient} solution is defined as a solution $x \in \X$ whose objective vector $z(x)$ is not \emph{strictly dominated} by that of another solution $\bar{x} \in \X$. 
We aim to identify a \emph{complete set of efficient (extreme) solutions} that contains at least one efficient (extreme) solution per non-dominated point $z \in \Z_n$.
\end{definition}

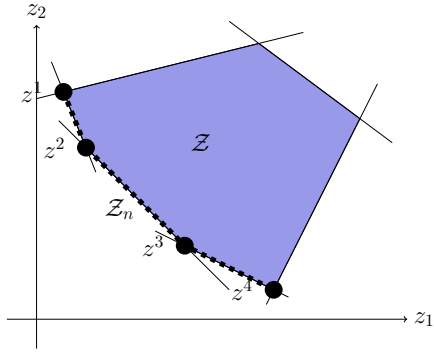
\begin{figure}
\begin{center}
\scalebox{0.65}{
  \begin{tikzpicture}[x=30mm, y=30mm]
    \draw[->,name path=xaxis] (-0.2,0) -- (2.5,0) node[right] {\Large $z_1$};
    \draw[->,name path=yaxis] (0,-0.2) -- (0,2.0) node[above] {\Large $z_2$};

    \draw[name path=line1,domain=0.1:0.4] plot (\x,{2 - 2.5 * \x}) node[above right] {};
    \draw[name path=line2,domain=0.15:1.3] plot (\x,{1.5 - 1 * \x}) node[above right] {};
    \draw[name path=line3,domain=0.8:1.7] plot (\x,{1 - 0.5 * \x}) node[above right] {};
    \draw[name path=line4,domain=1.55:2.3] plot (\x,{-3+ 2 * \x}) node[above right]{};
    \draw[name path=line5,domain=0:1.8] plot (\x,{1.5 + 0.25 * \x}) node[above right] {};
    \draw[name path=line6,domain=1.3:2.4] plot (\x,{3- 0.75 * \x}) node[above right] {};

    \node[name intersections={of=line1 and line2}] (a) at (intersection-1) {};
    \node[name intersections={of=line2 and line3}] (b) at (intersection-1) {};
    \node[name intersections={of=line3 and line4}] (c) at (intersection-1) {};
    \node[name intersections={of=line4 and line5}] (d) at (intersection-1) {};
    \node[name intersections={of=line4 and line6}] (e) at (intersection-1) {};
    \node[name intersections={of=line5 and line6}] (f) at (intersection-1) {};
    \node[name intersections={of=line5 and line1}] (g) at (intersection-1) {};

    \node (h) at (1.6,1.2) {};
    \node (i) at (0.75,0.75) {};

    \filldraw[ultra thin, fill=blue!80!black,fill opacity=0.4] (a.center)  -- (b.center) -- (c.center)  -- (d.center) -- (e.center) -- (f.center) -- (g.center) -- cycle;

     \path let \p0 = (h) in node [left=1cm of h] {\Large $\mathcal{Z}$};
     \path let \p0 = (i) in node [left = 0cm of i] {\Large $\mathcal{Z}_n$};

    \draw[line width=3, black, dashed] (g.center) node[circle,fill=black, label=left:\textcolor{black}{\Large $z^1$}]{} -- (a.center) node[circle,fill=black, label=left:\textcolor{black}{\Large $z^2$}]{} -- (b.center) node[circle,fill=black, label=left:\textcolor{black}{\Large $z^3$}]{} -- (c.center) node[circle,fill=black, label=left:\textcolor{black}{\Large $z^4$}] {};

  \end{tikzpicture}
  }
  \caption{Feasible set $\Z$ in objective space, set of non-dominated points $\Z_n$ (black dotted) and non-dominated extreme  points $z^1,z^2,z^3,z^4$.} \label{fig:objective_space}
\end{center}
\label{fig::NDs}
\end{figure}

The feasible set in decision space of an \gls{LP} (and a \gls{BLP}) is a polyhedron, and so is the feasible set $\Z$ of \gls{BLP} in objective space, see also Figure \ref{fig:objective_space}. The set of non-dominated points $\Z_n$ is shown as a dotted line in the figure, and it lies on the lower-left boundary of $\Z$.

Different scalarisations for bi-objective and multi-objective optimisation problems can be used to convert the original problem into a related single-objective problem. A widely used scalarisation transforms the original problem of the form \eqref{eq:BLPstandard} by weighting and summing the objectives. Given a weight $\lambda \in [0,1]$ the weighted-sum scalarisation $\lambda$-\gls{BLP} derived from \eqref{eq:BLPstandard} is
\begin{equation}
\label{eq:BLPstandard_lambda}
\begin{array}{rrcrcl}
   \min & \left ( \lambda(\hat{c}^1) + (1-\lambda) (\hat{c}^2) \right )^{\top} x &&\\
   \text{s.t.} & \hat{A}x  &=& \hat{b}, \\
               &  x &\geqq& 0.
\end{array}
\end{equation}
In the following we will refer to the weighted-sum scalarisation of any bi-objective optimisation problem, such as \gls{BLP} or \eqref{eq:BLPstandard}, as $\lambda$-\gls{BLP} or $\lambda$-\eqref{eq:BLPstandard} for $\lambda \in [0,1]$.

It is known that a feasible solution $x \in \X$ of \gls{BLP} \eqref{eq:BLPstandard} is efficient if and only if there exists a weight $\lambda \in (0,1)$ such that $x$ is an optimal solution of $\lambda$-\eqref{eq:BLPstandard}, \cite{Isermann74}. As illustrated in Figure~\ref{fig:objective_space}, all non-dominated points lie on facets of the polyhedron in objective space, where any non-dominated point that is not an extreme point corresponds to an optimal solution of $\lambda$-\eqref{eq:BLPstandard} for exactly one value of $\lambda$. Some non-dominated points are optimal for a range of $\lambda$-values and lie on the intersection of adjacent facets, such as points $z^1,z^2,z^3,z^4$ in Figure \ref{fig:objective_space}. These non-dominated points are termed non-dominated \emph{extreme} points and they correspond to efficient \emph{extreme} solutions. Note that, for a problem with bounded set of efficient solutions (Remark \ref{remark:bounded}), all non-dominated points of a \gls{BLP} can be obtained as convex combinations of neighbouring non-dominated extreme points, and therefore a complete set of efficient solutions (and corresponding non-dominated points) can be obtained by taking all convex combinations of neighbouring efficient extreme solutions (and corresponding non-dominated extreme points). The non-dominated extreme points are $z^1, z^2, z^3, z^4$ in Figure \ref{fig:objective_space}.

To solve a \gls{BLP} by weighted-sum scalarisation one must iteratively solve
single-objective weighted-sum problems to obtain a complete set of efficient extreme solutions, and methods such as the dichotomic approach \cite[e.g.][]{Aneja,CCS79} provide a streamlined approach for this. The bi-objective parametric simplex algorithm is an alternative approach that does not iteratively solve
single-objective weighted-sum problems. This method determines a complete set of efficient extreme solutions $\X_e = \{x^1,x^2,\ldots x^{q} \}$ with increasing first objective value, i.e.~$z_1(x^1) < z_1(x^2) < \ldots < z_1(x^{q})$. Each efficient extreme solution $x^i$ is an optimal solution of $\lambda$-BLP for all $\lambda \in [ a_{i+1},a_{i}]$ with $a_{i} \geqq a_{i+1}$ and $a_1 = 1, a_{q+1} = 0$ for $i=1,\ldots,q$.
We provide an overview of the bi-objective parametric simplex algorithm and further details can be found in~\cite{Ehrgott}. Through a sequence of basis updates, the algorithm iteratively moves from one non-dominated point to the next, which corresponds to proceeding from left to right in Figure~\ref{fig:objective_space}.
The algorithm starts by finding an efficient extreme solution $x^*$ that is optimal with respect to the first objective function. Associated with this solution is a set of basic variable indices and a set of non-basic variable indices.
As per the standard simplex algorithm, the bi-objective version iteratively performs basis exchanges by identifying entering and leaving variable pairs. In a bi-objective setting, the entering variable is identified to be a non-basic variable with maximum improvement of the second objective function over the deterioration of the first objective function. This ratio can be computed using the reduced costs, $\overline{c}^{k}$, for objective functions $k=1,2$.
If an entering variable exists, then the leaving variable can be determined (if it exists) as per the simplex algorithm for single-objective \gls{LP}.
This iterative procedure of identifying an entering variable
terminates when no entering variable can be found (i.e.~the basis is optimal with respect to the second objective function).

\begin{remark} \label{remark:bounded}
We assume $\X \neq \emptyset$. In general the set of efficient solutions of a \gls{BLP} is made up of efficient extreme solutions, and / or unbounded efficient edges, or rays, representing directions in which the set of non-dominated points in objective space is unbounded \citep{RUV17}. 
We also assume that  the set of efficient solutions $\X_e$ of all considered bi-objective optimisation problems is bounded (no unbounded efficient edges exist). This means that the weighted-sum scalarisation of the \gls{BLP} is also bounded,
i.e.~that an optimal solution of $\lambda$-\eqref{eq:BLPstandard} exists for $\lambda \in [0,1]$.  It should be noted that the parametric simplex algorithm (and our proposed algorithm) can be adapted to the unbounded case for \gls{BLP} following \citet{RUV17}.
\end{remark}

\subsection{Single-Objective Benders Decomposition}
\label{subsec:single-obj-benders}

Problem size is critical when solving optimisation problems. In many practical applications, the number of decision variables or constraints is prohibitively large. Benders decomposition is a technique that can be used to deal with large-scale \glspl{LP}. As this paper extends this approach to a bi-objective setting, we provide a comprehensive overview of the approach. Further details can be found in e.g.,~\cite{benders1962,Taskin10}.
Consider a single-objective \gls{LP} of the form
\begin{equation}
\label{eq:LP}
\begin{array}{rrcrcl}
   \min & c^{\top} x &+& f^{\top} y && \\
   \text{s.t.} & Ax  &+&  By & \geqq & b, \\
               &     & &  Dy & \geqq & d,\\
               &&&  x,y &\geqq& 0,
\end{array}
\end{equation}
where again $x\in \R^n$, $y \in \R^q$ and vectors $c,f,b,d$ and matrices $A,B,D$ have appropriate dimensions. One can observe that \eqref{eq:LP} reduces to a problem in terms of the $x$ variables once the  the $y$ variables have been fixed to certain values, and can therefore be restated as \eqref{eq:LP-restated}  where $\theta(\bar{y})$ is the optimal objective value of \eqref{eq:SubP}, which must be considered for all $\bar{y}\in \mathcal{Y}$, where $\mathcal{Y}=\{y \,|\, Dy\geqq d, y\geqq 0\}$.

\begin{minipage}{0.45\textwidth}
\begin{equation}
\label{eq:LP-restated}
\begin{array}{rrcrcl}
   \min & \theta(y) &+& f^{\top} y && \\
   \text{s.t.} & & &   Dy & \geqq& d,\\
               & & & y &\geqq& 0,
\end{array}
\end{equation}
\end{minipage}
\begin{minipage}{0.09\textwidth}
\end{minipage}
\begin{minipage}{0.45\textwidth}
\begin{equation}
\label{eq:SubP}
\begin{array}{rrcrcl}
   \min & c^{\top} x && \\
   \text{s.t.} &  Ax & \geqq & b &-& B\bar{y},\\
               &  x &\geqq& 0, &&
\end{array}
\end{equation}
\end{minipage}

\bigskip

If~\eqref{eq:SubP} is bounded, then the value of $\theta(\bar{y})$ can also be obtained by optimising its dual formulation. Assuming that $\pi$ denotes the vector of dual variables associated with the constraints of~\eqref{eq:SubP}, then the corresponding dual formulation can be stated as:
\begin{equation}
\label{eq:SubD}
\begin{array}{rrcl}
   \max &(b - B\bar{y})^{\top} \pi  \\
   \text{s.t.} & A^{\top} \pi &\leqq& c,\\
               & \pi & \geqq& 0.
\end{array}
\end{equation}
An advantage of solving the dual formulation is that its feasible region is independent of the values of the $y$ variables. The feasible region $\D$ of \eqref{eq:SubD} can be stated in terms of its set of extreme points $\D_p$ and extreme rays $\D_r$, i.e.~$\D = \{\pi \in \R^m | A^{\top} \pi \leqq c, \pi \geqq 0 \} = \text{conv}(\D_p) + \text{cone}(\D_r)$. Therefore \eqref{eq:SubD} can be restated to give:
\begin{equation}
\label{eq:SubD-restated}
\begin{array}{rrcccll}
   \min & \theta && && \\
   \text{s.t.} &&& (b - B\bar{y})^{\top} \pi_r &\leqq& 0 & \; \forall \pi_r \in \D_r, \\
               &-\theta&+& (b - B\bar{y})^{\top} \pi_p &\leqq& 0 & \; \forall\pi_p \in \D_p,
\end{array}
\end{equation}
where $\theta\in\mathbb{R}$ is the only decision variable. The first set of constraints, termed {\em feasibility cuts}, ensures that any $\bar{y}\in\mathcal{Y}$ that results in an unbounded direction for~\eqref{eq:SubD} is infeasible. The second set of constraints, termed {\em optimality cuts}, ensures that the $\theta$ variable correctly captures the optimal objective value of the $x$ variables.
Using~\eqref{eq:SubD-restated}, the following Benders reformulation of the original problem~\eqref{eq:LP} can be obtained.
\begin{equation}
\label{eq:LP-Benders}
\begin{array}{rrcccll}
   \min & \theta &+& f^{\top} y && \\
   \text{s.t.} &&&    Dy & \geqq & d, &\\
               &&& (b - By)^{\top} \pi_r &\leqq& 0 & \; \forall \pi_r \in \D_r, \\
               &-\theta&+& (b - By)^{\top} \pi_p &\leqq& 0 & \; \forall \pi_p \in \D_p, \\
               &\theta&,&  y &\geqq& 0. &
\end{array}
\end{equation}
Practical applications of Benders decomposition iteratively identify and add violated feasibility and optimality cuts rather than enumerating the sets $\D_r$ and $\D_p$ a priori, which may be prohibitively large. The {\it master problem} is a relaxation of~\eqref{eq:LP-Benders} that includes only subsets of feasibility and optimality cuts. An iterative approach finds candidate solutions $(\overline{y}, \overline{\theta})$ that are then subsequently checked for infeasibility/optimality by the {\it subproblem}~\eqref{eq:SubD}. If~\eqref{eq:SubD} is unbounded, a feasibility cut is added to the master problem, while, if~\eqref{eq:SubD} is feasible, a check is made to see whether or not $\overline{y}$ is optimal. An optimality cut can be generated and added to the master problem to remove the provably suboptimal solution $\overline{y}$ if $\theta(\overline{y})>\overline{\theta}$.
Lower and upper bounds on the optimal objective value of the master problem are obtained at each iteration, and the problem is solved once the bounds 
are considered sufficiently close in value.

\subsection{Benders Reformulation of \gls{BLP}}
\label{sec:full-Benders-reformulation-BLP}
We now focus on developing a bi-objective variant of Benders decomposition and consider \glspl{BLP} of the form
\begin{equation}
\label{eq:BLP}
\begin{array}{rrcrcl}
   \min & \left ( \begin{array}{c} (c^1)^{\top} x \\ (c^2)^{\top} x \end{array} \right . & \left . \begin{array}{c} + \\ + \end{array} \right . & \left . \begin{array}{c} (f^1)^{\top} y \\ (f^2)^{\top} y \end{array} \right ) & & \\
   \text{s.t.} & Ax  &+&  By &\geqq& b, \\
               &     & &  Dy &\geqq& d,\\
               &&&  x,y &\geqq& 0,
\end{array}
\end{equation}
where $x\in \R^n$, $y \in \R^q$ and vectors $c^1,c^2,f^1,f^2,b,d$ and matrices $A,B,D$ have appropriate dimensions. Note that \gls{BLP}~\eqref{eq:BLPstandard} was originally formulated in its standard form, whereas we use the more convenient inequality form \eqref{eq:BLP} here.
To apply Benders decomposition to this problem in which we leave the $y$ variables in the master problem and place the $x$ variables in the subproblem, we proceed as in~Section~\ref{subsec:single-obj-benders} and reformulate~\eqref{eq:BLP} as a \gls{BM} in \eqref{eq:BLP-restated} where $\theta_1(\bar{y})$ and $\theta_2(\bar{y})$ bound the objective values of the \gls{BS} in \eqref{eq:BSub}, which must be considered for all $\bar{y}\in\mathcal{Y}$: \\
\begin{minipage}{0.45\textwidth}
\begin{equation}
\label{eq:BLP-restated}
\begin{array}{rrcrcl}
    \min & \left ( \begin{array}{c} \theta_1(y) + (f^1)^{\top} y  \\ \theta_2(y) + (f^2)^{\top} y \end{array} \right ) & & \\

   \text{s.t.}  &  Dy &\geqq& d,\\
               &  y &\geqq& 0
\end{array}
\end{equation}
\end{minipage}
\begin{minipage}{0.09\textwidth}
\end{minipage}
\begin{minipage}{0.45\textwidth}
\begin{equation}
\label{eq:BSub}
\begin{array}{rrcrcl}
  \min & \left ( \begin{array}{c} (c^1)^{\top} x \\ (c^2)^{\top} x \end{array} \right ) \\
  \text{s.t.} &  Ax &\geqq& b &-& B \bar{y},\\
               &  x &\geqq& 0 &&
\end{array}
\end{equation}
\end{minipage}

It should be noted that \eqref{eq:BSub} will generally have infinitely many efficient solutions in the decision space and non-dominated points in the objective space. Variables $\theta_1$ and $\theta_2$ are used to capture the contribution of variables $x$ to the two individual objective functions in \eqref{eq:BLP-restated}, enabling the correct representation of non-dominated points in objective space.

How to handle a bi-objective subproblem \eqref{eq:BSub} in the context of Benders decomposition is an interesting question to ask. We first consider the impact of using two separate subproblems, one for each objective, and investigate the conditions under which such an approach is sufficient for obtaining a complete set of efficient extreme solutions to \gls{BLP} \eqref{eq:BLP}. Both subproblems are of the form~\eqref{eq:SubP} with the exception that the first (respectively, second) has the objective function $(c^1)^{\top} x$ (respectively, $(c^2)^{\top}x$). Continuing as in the single-objective case, we identify the set of extreme points for the dual formulation of each subproblem. These are denoted $\mathcal{D}_p^1$ and $\mathcal{D}_p^2$, respectively. Note that the set of extreme rays to each of the two dual formulations is the same as the set of rays $\mathcal{D}_r$ is independent of the subproblem objective in~\eqref{eq:SubP}. This leads to the following Benders reformulation:
\begin{equation}
\label{eq:BLP-Bmaster-Benders-SO-subproblems-prelim}
\begin{array}{rrcrcrcll}
   \min & \left ( \begin{array}{c} \theta_1 \\ \phantom{.} \end{array} \right . &\left . \begin{array}{c} + \\ \phantom{.} \end{array} \right . & \left . \begin{array}{c} \\ \theta_2 \end{array} \right . & \left . \begin{array}{c}  \\ + \end{array} \right . & \left . \begin{array}{c} (f^1)^{\top} y \\ (f^2)^{\top} y \end{array} \right )&&& \\
   \text{s.t.} & & &     & &  Dy &\geqq& d,&\\
               &&&&& (b - By)^{\top} \pi_r &\leqq& 0 & \; \forall \pi_r \in \D_r, \\
               &-\theta_1&+& && (b - By)^{\top} \pi_p &\leqq& 0 & \; \forall \pi_p \in \D_p^1, \\
               &&&-\theta_2&+&  (b - By)^{\top} \pi_p &\leqq& 0 & \; \forall \pi_p \in \D_p^2, \\
               &&& &&  y &\geqq& 0.&
\end{array}
\end{equation}
If one then proceeded
as in the single-objective case~\eqref{eq:LP-Benders}, dynamically separating violated optimality and feasibility cuts using two subproblems instead of one,  one
might expect that~\eqref{eq:BLP-Bmaster-Benders-SO-subproblems-prelim} correctly describes a complete set of efficient extreme solutions of the original problem~\eqref{eq:BLP}.
Our first observation is  that \eqref{eq:BLP-Bmaster-Benders-SO-subproblems-prelim}, which only relies on single-objective subproblems, is not able to fully describe all efficient solutions, as the following example illustrates.

\begin{example} \label{ex:BLP-needs-weighted-cuts}
Consider the following \gls{BLP}:
\begin{equation}
\begin{array}{rrcrcrcrcrll}
\label{eq:ex1-NP1}
  \min & \left(\begin{array}{r} 4x_1 \\ -2x_1 \end{array}\right. & \left . \begin{array}{c} - \\ - \end{array} \right .& \left . \begin{array}{c} x_2 \\ 2x_2 \end{array} \right .&\left . \begin{array}{c} + \\ + \end{array} \right .&\left . \begin{array}{c} 2y_1 \\ 4y_1 \end{array} \right .&\left . \begin{array}{c} \phantom{+} \\ - \end{array} \right .&\left . \begin{array}{c} \phantom{.} \\ 6y_2 \end{array} \right .&\left . \begin{array}{c} - \\ - \end{array} \right .&\left . \begin{array}{c} 4 y_3 \\ 3y_3 \end{array} \right )&&\\
   \text{s.t.} &&&&& y_1 &-& 6y_2& -& 4y_3 &\geqq& -2, \\
     & -2 x_1 &-& 6x_2 &-& 4y_1 &-& 3 y_2 &-& 6y_3 &\geqq& -5,  \\
     & -5 x_1 && && &-& 3y_2 &-& 5y_3 &\geqq& -2, \\
     & x_1, && x_2, &&y_1, && y_2,&& y_3 &\geqq& 0.
\end{array}
\end{equation}
The non-dominated points of~\eqref{eq:ex1-NP1} are shown in Figure~\ref{fig:NP1-reformulation}. The set of non-dominated extreme  points is given by
$\Z_n=\{ z^1=\left(- \textstyle\frac{61}{30}, -\frac{62}{30}\right), z^2=\left(-\textstyle\frac{16}{9}, -\frac{23}{9}\right), z^3=\left(-\textstyle\frac{2}{3}, -\frac{10}{3}\right), z^4=\left(\textstyle\frac{1}{5},-\frac{18}{5}\right)\}$,
and the corresponding set of efficient extreme solutions is
$ \X_e= \{ \left( 0, \frac{13}{30}, 0,0,\textstyle\frac{4}{10} \right ),\linebreak[1] \left (0, \frac{4}{9}, 0, \textstyle\frac{1}{9}, \frac{1}{3} \right ),\linebreak[1] \left (0, \frac{2}{3}, 0, \textstyle\frac{1}{3}, 0 \right ),\linebreak[1] \left (\frac{1}{5}, \frac{3}{5}, 0, \textstyle\frac{1}{3},0 \right )  \}$,
where the variables in each efficient extreme solution appear in the order $(x_1, x_2,y_1, y_2, y_3 )$.
\begin{figure}[tbp]
    \begin{minipage}[adjusting]{0.48\textwidth}
      \begin{center}
      \begin{tikzpicture}[x=30mm, y=30mm]
        \draw[very thin,color=gray] (0.0,0.0) grid (1.5,1.0);
        \draw[->,thick,name path=xaxis] (-0.2,0) -- (1.5,0) node[below] {$\pi^1_1$};
        \draw[->,thick,name path=yaxis] (0,-0.2) -- (0,1.0) node[left] {$\pi^1_2$};
        \draw[name path=line3] (1/6,1.0) -- (1/6, -0.2);
        \node[name intersections={of=line3 and xaxis}] (c) at (intersection-1) {};
        \node (d) at (1.5,0) {};
        \node (e) at (1.5,1.0) {};
        \node (f) at (1/6, 1.0) {};
        \filldraw [red] (0.17,0) circle (3pt);
        \filldraw[ultra thin,color=blue!10,fill=blue!80!black,fill opacity=0.1] (c.center)  -- (d.center) -- (e.center)  -- (f.center) -- cycle;
        \path let \p0 = (c) in node [above right=0cm and 0.025cm of c.center] {($\frac{1}{6}$,\;0)};
        \draw[color=red,thick] (0.2,0) -- (1.5,0) node[above] {};
        \draw[color=red,thick] (1/6,1.0) -- (1/6, -0.2);
      \end{tikzpicture}
    \end{center}
    \end{minipage}
    \begin{minipage}[adjusting]{0.48\textwidth}
      \begin{center}
      \begin{tikzpicture}[x=30mm, y=30mm]
        \draw[very thin,color=gray] (0.0,0.0) grid (1.5,1.0);
        \draw[->,thick,name path=xaxis] (-0.2,0) -- (1.5,0) node[below] {$\pi^2_1$};
        \draw[->,thick,name path=yaxis] (0,-0.2) -- (0,1.0) node[left] {$\pi^2_2$};
        \draw[name path=line2,domain=-0.1:1.5] plot (\x,{0.4-0.4 * \x}) node[above right] {};
        \draw [dashed, color=red, name path=line3] (2/6,1.0) -- (2/6, -0.2);
        \node[name intersections={of=line2 and xaxis}] (b) at (intersection-1) {};
        \node[name intersections={of=line2 and line3}] (a) at (intersection-1) {};
        \node (c) at (1.5,0) {};
        \node (d) at (1.5,1.0) {};
        \node (e) at (2/6, 1.0) {};
        \filldraw[ultra thin, color=green!10,fill=green!80!black,fill opacity=0.4] (a.center)  -- (b.center) -- (c.center)  -- (d.center) -- (e.center) -- cycle;
        \path let \p0 = (a) in node [above right=-0.15cm of a]  {($\frac{1}{3}$,\;$\frac{4}{15}$)};
        \path let \p0 = (b) in node [above right=-0.10cm of b] {(1,\;0)};
        \filldraw [red] (1,0) circle (3pt);
        \filldraw [red] (0.333,0.27) circle (3pt);
        \draw[->,name path=xaxis] (-0.2,0) -- (1.5,0) node[above] {};
        \draw[color=red, thick] (-0.2,0) -- (1.5,0) node[above] {};
        \draw[color=red, thick, domain=-0.1:1.5] plot (\x,{0.4- 0.4 * \x}) node[above right] {};
        \draw[color=red, thick] (2/6,1.0) -- (2/6, -0.2);
      \end{tikzpicture}
      \end{center}
\end{minipage}
    \caption{Feasible region of dual formulation for the first subproblem (left) and the second (right).}
    \label{fig:NP1-SO-subproblems}
\end{figure}
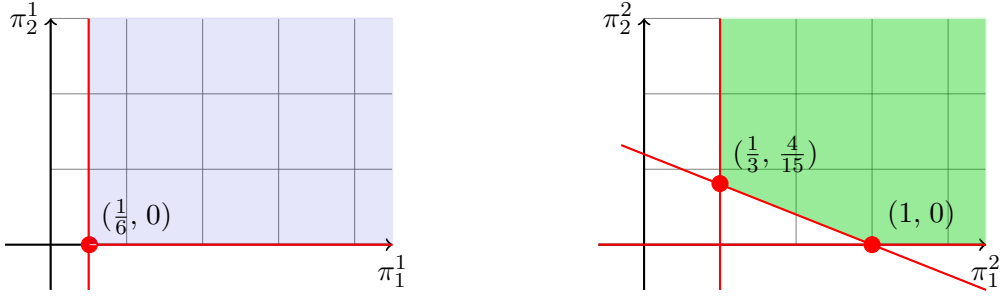

We apply a Benders decomposition in which the $y$ variables are placed in the master problem, and the $x$ variables are placed in the two single-objective subproblems. The feasible regions of the resulting dual subproblem formulations are illustrated in Figure~\ref{fig:NP1-SO-subproblems}; the blue (respectively, green) area in the left (respectively, right) plot  corresponds to the first (respectively, second) objective.
The sets of extreme points for each feasible region are $\mathcal{D}^1_p:=\left\{ \left(\textstyle\frac{1}{6},0\right) \right\}$ and $\mathcal{D}^2_p:=\left\{(1,0), \left(\textstyle\frac{1}{3},\frac{4}{15}\right) \right \}$.
There are also two extreme rays $\D_r = \left \{ (1,0), (0,1) \right \}$ for each space. The Benders reformulation can therefore be stated in the form of~\eqref{eq:BLP-Bmaster-Benders-SO-subproblems-prelim} as
\begin{equation}
\begin{array}{rlll}
\label{eq:ex1-NP1-reformulation}
  \min & \left ( \begin{array}{l} \theta_1 + 2y_1 - 4 y_3 \\ \theta_2 + 4y_1 - 6y_2 - 3y_3 \end{array} \right ) & \\
   \text{s.t.} & y_1 - 6y_2 - 4y_3 &\geqq -2, \\
     & (4y_1+3y_2+6y_3-5)\cdot (\pi^r_1) + (3y_2+5y_3-2)\cdot (\pi^r_2) &\leqq 0, & \forall (\pi^r_1,\pi^r_2)\in \mathcal{D}_r\\
     &(4y_1+3y_2+6y_3-5)\cdot \left(\pi^1_1\right) + (3y_2+5y_3-2)\cdot (\pi^1_2) &\leqq \theta_1, &\forall (\pi^1_1,\pi^1_2)\in \mathcal{D}_p^1  \\
     &(4y_1+3y_2+6y_3-5)\cdot \left(\pi^2_1\right) + (3y_2+5y_3-2)\cdot (\pi^2_2) &\leqq \theta_2, &\forall (\pi^2_1,\pi^2_2)\in \mathcal{D}_p^2  \\
    & \theta_1, \theta_2, y_1, y_2, y_3 &\geqq 0.
\end{array}
\end{equation}
Reformulation~\eqref{eq:ex1-NP1-reformulation} only has three non-dominated extreme points $\{z^1,z^2,z^5\}$. It correctly identifies points $z^1$ and $z^2$ as shown in Figure~\ref{fig:NP1-reformulation} but incorrectly identifies point  $z^5=\left(-\textstyle\frac{2}{3}, -\frac{18}{5}\right)$. This point is defined by the first objective value for point $z^3$ and the second objective value for point $z^4$. Points $z^3$ and $z^4$ have identical master variable values $(y_1,y_2,y_3) = \left (0, \textstyle\frac{1}{3}, 0 \right )$ but different subproblem variable values $(x_1,x_2)$.

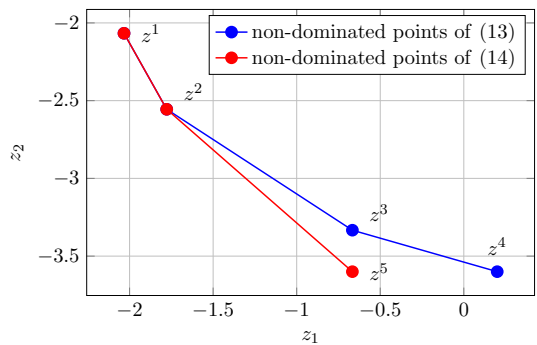
\begin{figure}
\begin{center}
\scalebox{0.7}{
  \begin{tikzpicture}
    \begin{axis}[
        xlabel=$z_1$,
        ylabel=$z_2$,
        width=10cm,height=7cm,
        grid
      ]

        \node[label={[label distance=0.05cm]0:$z^1$}] (a) at   (-2.0333, -2.06667) {};
        \node[label={[label distance=0.05cm]10:$z^2$}] (a) at   (-1.7778, -2.5556) {};
        \node[label={[label distance=0.05cm]10:$z^3$}] (a) at  (-0.6667, -3.3333) {};
        \node[label={[label distance=0.05cm]90:$z^4$}] (a) at    (0.2, -3.6) {};
        \node[label={[label distance=0.05cm]0:$z^5$}] (a) at     (-0.6667, -3.6) {};

      \addplot[thick, mark=*, color=blue,mark size=3pt] coordinates {
        (-2.0333, -2.06667)
        (-1.7778, -2.5556)
        (-0.6667, -3.3333)
        (0.2, -3.6)
      };
      \addplot[thick, mark=*, color=red,mark size=3pt] coordinates {
        (-2.0333, -2.06667)
        (-1.7778, -2.5556)
        (-0.6667, -3.6)
        };
    \legend{non-dominated points of \eqref{eq:ex1-NP1},non-dominated points of \eqref{eq:ex1-NP1-reformulation}}
    \end{axis}
  \end{tikzpicture}
  }
\end{center}
    \caption{Non-dominated points of the original problem~\eqref{eq:ex1-NP1} and a Benders reformulation \eqref{eq:ex1-NP1-reformulation}.}
    \label{fig:NP1-reformulation}
\end{figure}

This example demonstrates that it is insufficient to only add optimality cuts corresponding to extreme dual solutions to the single-objective subproblems.
To understand why, we analyse the efficient extreme solution associated with point $z^4$. In this solution, the subproblem variables take the values $(x_1,x_2) = \left(\frac{1}{5},\frac{3}{5} \right)$.
The corresponding dual solutions are $\left(\pi^1_1,\pi^1_2\right)=\left(\frac{1}{6},-\frac{13}{15}\right)$ for the first subproblem and $\left(\pi^2_1,\pi^2_2\right)=\left(\frac{1}{3},\frac{4}{15}\right)$ for the second, respectively. We note that the first dual solution is infeasible. The explanation for this is that the solution $\left(x_1,x_2\right)=\left(\frac{1}{5}, \frac{3}{5}\right)$ is not optimal for the first objective. The second dual solution is feasible (and hence optimal for the second objective). To ensure $\theta_1$
accurately captures the relevant subproblem objective function value, the following
cut could be defined:
\begin{equation}
\label{eq:new-cut}
\begin{array}{rl}
    (4y_1+3y_2+6y_3-5)\cdot \left(\frac{1}{6}\right) + (3y_2+5y_3-2)\cdot \left(-\frac{13}{15}\right) & \leqq \theta_1 
\end{array}
\end{equation}
This is not present in~\eqref{eq:BLP-Bmaster-Benders-SO-subproblems-prelim}
since it is derived from a dual solution that is not an extreme point of the dual feasible region for the subproblem associated with the first objective. Evaluating this constraint at the master solution $(y_1,y_2,y_3) = \left (0, \textstyle\frac{1}{3}, 0 \right )$ enforces $\theta_1\geq \frac{1}{5}$, the correct bound for the first objective function value in $z_4$. However, the bound implied in~\eqref{eq:ex1-NP1-reformulation} at the same master solution is $\theta_1\geq -\frac{2}{3}$, (due to the single optimality cut for $\theta_1$). This is the correct bound for the first objective function value in $z_3$. Being more restrictive, the addition of~\eqref{eq:new-cut} to~\eqref{eq:ex1-NP1-reformulation} would render it impossible to generate the third non-dominated extreme  point $\left(z_1,z_2\right)= \left(-\frac{2}{3}, -\frac{10}{3}\right)$. From this we conclude that single-objective optimality cuts alone are insufficient.

A weighted cut, which weights the values of $\theta_1$ and $\theta_2$ in an optimality cut, would allow $\theta_1$ and $\theta_2$ to take different values for the same master solution $(y_1,y_2,y_3) = \left (0, \textstyle\frac{1}{3}, 0 \right )$. Such a cut enables the Benders reformulation to correctly identify the two efficient extreme solutions with $(x_1,x_2) = \left(\frac{1}{5},\frac{3}{5} \right)$ and $(x_1,x_2) = \left (0,\frac{2}{3} \right)$, respectively.
We demonstrate this by adding a weighted cut with $\lambda=\frac{4}{17}$. This weight is specifically chosen since it defines the weight for which the non-dominated extreme points $z^3$ and $z^4$ are optimal for $\lambda$-\eqref{eq:ex1-NP1}. The weighted cut is of the form:
$$(4y_1+3y_2+6y_3-5)\cdot \left(\frac{5}{17} \right) + (3y_2+5y_3-2)\cdot (0)  \leqq \lambda \cdot \theta_1+(1-\lambda) \cdot \theta_2,$$
where $\left(\pi_1,\pi_2\right)=\left(\frac{5}{17},0\right)$ is the optimal dual solution of the $\lambda=\frac{4}{17}$ weighted subproblem that is induced by the master problem solution $(y_1,y_2,y_3) = \left (0, \textstyle\frac{1}{3}, 0 \right )$. Adding just this cut to reformulation~\eqref{eq:ex1-NP1-reformulation} yields the correct set of efficient extreme solutions (and hence non-dominated extreme  points), as shown in Figure~\ref{fig:NP1-reformulation}.
\end{example}

In summary, Example~\ref{ex:BLP-needs-weighted-cuts} illustrates, based on problem \eqref{eq:ex1-NP1}, that the Benders reformulation \eqref{eq:BLP-Bmaster-Benders-SO-subproblems-prelim}
 of~\eqref{eq:BLP} requires weighted optimality cuts.
Assume that $(\bar{x},\bar{y})$ is an efficient solution of \eqref{eq:BLP}. Then there exists $\lambda \in (0,1)$ \citep{Isermann74} such that $(\bar{x},\bar{y})$ is an optimal solution of the weighted-sum problem \eqref{eq:BLP} with objective function $(c^{\lambda})^{\top} x + (f^{\lambda})^{\top} y$ with $c^{\lambda} = \lambda c^1 + (1-\lambda) c^2$ and $f^{\lambda} = \lambda f^1 + (1-\lambda) f^2$. Applying this weighting factor $\lambda$ to \gls{BM} \eqref{eq:BLP-restated} with
$c^{\lambda},f^{\lambda}$ as before, leads to a weighted (primal) subproblem \eqref{eq:BSub-weighted}, which we will also refer to as $\lambda$-\eqref{eq:BSub}:
\begin{equation}
\label{eq:BSub-weighted}
\begin{array}{rrcrcl}
   \min & (c^{\lambda})^{\top} x \\
   \text{s.t.} &  Ax &\geqq& b &-& B \bar{y},\\
               &  x &\geqq& 0. &&
\end{array}
\end{equation}
The feasible region of the dual of $\lambda$-\eqref{eq:BSub} can be described as $\D^{\lambda} = \{\pi \in \R^m | A^{\top} \pi \leqq  c^{\lambda}, \pi \geqq 0 \} = \text{conv}(\D_p^{\lambda}) + \text{cone}(\D_r)$.
The (single-objective) Benders reformulation (see also \eqref{eq:SubD-restated}) of the $\lambda$-weighted problem $\lambda$-\eqref{eq:BLP} therefore includes Benders feasibility cuts of the form $(b-B\bar{y})^{\top} \pi_r \leqq 0$  as well as optimality cuts of the form $-\theta + (b-B\bar{y})^{\top} \pi_p^{\lambda} \leqq 0$ for extreme rays $\pi_r \in \D_r$ and extreme points $\pi_p^{\lambda} \in \D_p^{\lambda}$ of the dual of the $\lambda$-weighted subproblem. 

\begin{definition} \label{def:weighted_optimality_cut}
    Given $\lambda \in [0,1]$ and a corresponding $\pi_p^{\lambda} \in \D_p^{\lambda}$.
    Then a \emph{weighted optimality cut} takes the following form:
    \begin{equation} \label{eq:weighted_optimality_cut}
    (b - By)^{\top} \pi^{\lambda}_p \leqq \lambda \theta_1 + (1-\lambda) \theta_2.
    \end{equation}

\end{definition}

We therefore obtain the following Benders reformulation of \gls{BLP}~\eqref{eq:BLP}, which includes weighted cuts of the form~\eqref{eq:weighted_optimality_cut} with $\lambda\in[0,1]$. As such, the single-objective cuts corresponding to $\lambda=0$ and $\lambda=1$ are included in the set of weighted cuts.
\begin{equation}
\label{eq:BLP-complete-reformulation}
\begin{array}{rrcrcrcll}
   \min & \left ( \begin{array}{c} \theta_1 \\ \phantom{.} \end{array} \right .&\left . \begin{array}{c} + \\ \phantom{.} \end{array} \right .&\left . \begin{array}{c}  \\ \theta_2 \end{array} \right .&\left . \begin{array}{c} \\ + \end{array} \right .&\left . \begin{array}{c} (f^1)^{\top} y \\ (f^2)^{\top} y \end{array} \right )&&& \\
   \text{s.t.} & & &     & &  Dy &\geqq& d,&\\
               &&&&& (b - By)^{\top} \pi_r &\leqq& 0 & \; \forall \pi_r \in \D_r, \\
               & - \lambda \theta_1 &-& (1-\lambda) \theta_2 &+&
               (b - By)^{\top} \pi^{\lambda}_p &\leqq&  0 & \; \forall \pi^{\lambda}_p \in \D_p^{\lambda}, \lambda\in[0,1], \\
               &&& &&  y &\geqq& 0.&
\end{array}
\end{equation}
Formulation~\eqref{eq:BLP-complete-reformulation} represents the complete Benders \gls{BLP} reformulation. While there is an infinite number of weighted cuts, we will show that only a finite set of them is required in~\eqref{eq:BLP-complete-reformulation}, and that it is always possible to generate this finite Benders \gls{BLP} reformulation. We use~\eqref{eq:BLP-complete-reformulation}-$\{ \bar{\D}_r,\Lambda,\bar{\D}_p\}$ to refer to a variant of \eqref{eq:BLP-complete-reformulation} that replaces $\D_r$ by $\bar{\D}_r$ and $\D_p^{\lambda}$ for $\lambda \in [0,1]$ by $ \bar{\D}_p(\lambda)$ for $\lambda \in \Lambda$, with $\bar{\D}_r \subseteq \D_r$ and $\bar{\D}_p(\lambda) \subseteq \D_p^{\lambda}$. In the following, we prove that, while the efficient solutions in decision space of \eqref{eq:BLP} and \eqref{eq:BLP-complete-reformulation} may be different, their respective sets of nondominated points in objective space are equal. Furthermore, the reformulation only requires a finite subset of cuts $\bar{\D}_r$ and $\bar{\D}_p (\lambda)$ for $\lambda \in \Lambda$.

\begin{proposition}
\label{prop:1}
Consider a \gls{BLP} of the form~\eqref{eq:BLP}. Given a complete set of efficient extreme solutions $\X_e=\left\{(x^1,y^1),(x^2,y^2),\ldots,(x^q,y^{q})\right\}$, corresponding $\lambda$ optimality ranges $[a_{i+1},a_{i}]$ for each $(x^i,y^i)\in\X_e$, where $a_1=1$ and $a_{q+1}=0$,  and the corresponding set of non-dominated extreme  points $\Z_n=\left\{(z^1_1,z^1_2),(z^2_1,z^2_2),\dots,(z^{q}_1,z^{q}_2)\right\}$, then $\X^r_e=\left\{(\theta_1^1,\theta_2^1,y^{1}),(\theta_1^2,\theta_2^2,y^{2}),\ldots,(\theta_1^{q},\theta_2^{q},y^{q})\right\}$, where $\theta_1^i=(c^1)^\top x^i$ and $\theta_2^i=(c^2)^\top x^i$ for $i=1,2\dots,q$, comprises a complete set of efficient extreme solutions to a Benders reformulation~\eqref{eq:BLP-complete-reformulation}-$\{ \bar{\D}_r,\Lambda,\bar{\D}_p\}$ and has the same set of non-dominated extreme points $\Z_n$, where $(z_1^i,z_2^i) = \left ((f^1)^\top y^{i}+\theta_1^i, (f^2)^\top y^{i}+\theta_2^i \right)$. In particular, we can replace  $\lambda\in[0,1]$ in~\eqref{eq:BLP-complete-reformulation} with $\lambda\in \Lambda = \left\{a_1,a_2,\dots,a_{q}\right\}$ in \eqref{eq:BLP-complete-reformulation}-$\{ \bar{\D}_r,\Lambda,\bar{\D}_p\}$.
\end{proposition}

\proof
We begin by observing that the constraint set of~\eqref{eq:BLP-complete-reformulation} is the union of the constraint sets of an infinite number of Benders reformulations of the single-objective $\lambda_k$-\eqref{eq:BLP} problem.
The same weighted-sum scalarisation can be applied to the two objective functions of \eqref{eq:BLP} and \eqref{eq:BLP-complete-reformulation} for some $\lambda_k\in[0,1]$. From the theory on single-objective Benders decomposition we know that the optimal objective value of $\lambda_k$-\eqref{eq:BLP} is equal to the optimal objective value of $\lambda_k$-\eqref{eq:BLP-complete-reformulation}.
Consequently, for any $\lambda_k\in[0,1]$ the following relationship between an optimal solution $(x^\prime, y^\prime)$ of $\lambda_k$-\eqref{eq:BLP} and an optimal solution $(\theta_1^{\prime\prime},\theta_2^{\prime\prime}, y^{\prime\prime})$ of $\lambda_k$-\eqref{eq:BLP-complete-reformulation} exists:
\begin{equation}\label{eq:decomp-reformulation-connection} (f^{\lambda_k})^\top y^\prime + (c^{\lambda_k})^\top x^\prime = f^\top y^{\prime\prime} +\lambda_k\theta_1^{\prime\prime}+(1-\lambda_k)\theta_2^{\prime\prime},\end{equation}
where $f^{\lambda_k} = (\lambda_k f^1+(1-\lambda_k) f^2)$ and $c^{\lambda_k} = (\lambda_k c^1+(1-\lambda_k) c^2)$.
Due to the weighting of $\theta_1''$ and $\theta_2''$, there is no guarantee that $\theta_1''=(c^1)^\top x^{\prime\prime}$ and $\theta_2''=(c^2)^\top x^{\prime\prime}$, where $x^{\prime\prime}$ is an optimal solution to $\min \; \{(c^{\lambda_k})^\top x | Ax \geqq b - B y^{\prime\prime}, x \geqq 0\}$, i.e., the subproblem \eqref{eq:BSub-weighted} that is induced by $y^{\prime\prime}$.

Each efficient extreme solution $(x^i,y^i)\in\X_e$ for $i=1,2,\dots,q$ is an optimal solution to $\lambda_k$-\eqref{eq:BLP} for all $\lambda_k \in [a_{i+1},a_{i}]$. Using~\eqref{eq:decomp-reformulation-connection}, we construct the solution $(\theta_1^i,\theta_2^i,y^{i})$, where $\theta_1^i=(c^1)^\top x^i$ and $\theta_2^i=(c^2)^\top x^i$ and with $z_1^i=(f^1)^\top y^{i}+\theta_1^i$ and $z_2^i=(f^2)^\top y^{i}+\theta_2^i$. This must be an optimal solution to $\lambda_k$-\eqref{eq:BLP-complete-reformulation} for all $\lambda \in [a_{i+1},a_{i}]$. If this were not the case, then for some $i=1,2,\dots,q$, there must exist a $\lambda_j\in \left(a_{i+1},a_{i}\right)$ and an optimal solution $(\theta_1^j,\theta_2^j,y^{j}$) to $\lambda_j$-\eqref{eq:BLP-complete-reformulation} for which 
$$(f^{\lambda_j})^\top y^{j}+\lambda_j\theta_1^{j}+(1-\lambda_j)\theta_2^{j} < (f^{\lambda_j})^\top y^i + (c^{\lambda_j})^\top x^i. $$
However, this contradicts the observation in~\eqref{eq:decomp-reformulation-connection}. If the solution $(\theta_1^j,\theta_2^j,y^{j}$) were an optimal solution to $\lambda_j$-$\eqref{eq:BLP-complete-reformulation}$, then there must exist a feasible solution $(x^j,y^j)$ to $\lambda_j$-$\eqref{eq:BLP}$ with the same objective value. This contradicts the assumption that $(x^i,y^i)$ is an optimal solution to $\lambda_k$-$\eqref{eq:BLP}$ for all $\lambda_k\in [a_{i+1},a_{i}]$. We can conclude that $(\theta_1^j,\theta_2^j,y^{j})$ does not exist.
Each solution $(\theta_1^i,\theta_2^i,y^{i})\in \X_e^r$ for $i=1,2,\dots,q$ is therefore an optimal solution to $\lambda_k$-\eqref{eq:BLP-complete-reformulation} for all $\lambda_k\in [a_{i+1},a_{i}]$ and has individual objective function values $z_1^i$ and $z_2^i$. We also note that any $y^i$, where $i=1,2,\dots, q$, can be used to provide an upper bound to $\lambda_k$-$\eqref{eq:BLP}$, and hence $\lambda_k$-$\eqref{eq:BLP-complete-reformulation}$, with $i\neq k$, since the feasibility of $y^i$ is independent of the objective function.
Feasibility cuts, defined by $\mathcal{D}_r$, do not depend on the $\lambda_k$ value.

For a given $\lambda_k\in[0,1]$ the constraints of~\eqref{eq:BLP-complete-reformulation} ensure that any feasible solution cannot have a better $\lambda_k$-weighted objective value than the optimal objective value of the weighted single objective problem $\lambda_k$-$\eqref{eq:BLP}$.
Thus, any solution that is efficient in \eqref{eq:BLP} will have a corresponding efficient solution in~\eqref{eq:BLP-complete-reformulation}. This is also true for efficient extreme solutions. Thus, $\X_e^r$ is a set of efficient extreme solutions to~\eqref{eq:BLP-complete-reformulation}.
The set of non-dominated extreme  points remains $\Z_n$ with this reformulation. Since, each solution $(\theta_1^i,\theta_2^i,y^{i})$ is optimal for the range $\lambda\in[a_{i+1},a_{i}]$ we can remove all $\lambda\in (a_{i+1},a_{i})$ for $i=1,2,\dots,q$ without changing the set $\X_e^r$. Therefore we can choose $\Lambda = \{ a_1, a_2, \ldots, a_q \}$, $\bar{\D}_r = \D_r$, and  $\bar{\D}_p(\lambda) =  \D_p^{\lambda}$ for $\lambda \in \Lambda$ to obtain \eqref{eq:BLP-complete-reformulation}-$\{ \bar{\D}_r,\Lambda,\bar{\D}_p\}$.
\endproof

We conclude this section by identifying conditions under which single-objective optimality cuts, defined by $\mathcal{D}_p^1$ and $\mathcal{D}_p^0$, are sufficient to define the Benders reformulation of~\eqref{eq:BLP}.

\begin{remark} \label{remark1}
Consider a Benders reformulation of~\eqref{eq:BLP} of the form~\eqref{eq:BLP-complete-reformulation}, with master variables $y$ and subproblem variables $x$. If the optimal solution to the weighted subproblem~\eqref{eq:BSub-weighted} is the same for $\lambda=1$ and $\lambda=0$ for any master \eqref{eq:BLP-restated} solution $\bar{y}$, i.e, if the two subproblem objectives in \eqref{eq:BSub} are not contradictory, then weighted cuts, obtained using any $\lambda_k\in(0,1)$ are not necessary in the Benders reformulation.

Please refer to \ref{app:proof_remark1} for a proof of Remark \ref{remark1}.
\end{remark}

\begin{remark} \label{remark2}
Consider a Benders reformulation of~\eqref{eq:BLP} of the form~\eqref{eq:BLP-complete-reformulation}, with master variables $y$ and subproblem variables $x$. Assume that the subproblem variables are only present in one objective. If, w.l.o.g., this is the first objective, we have $c^2 = 0$ in \eqref{eq:BLP}. Therefore, the subproblem \eqref{eq:BSub} reduces to a single-objective optimisation problem, variable $\theta_2$ can be omitted from \eqref{eq:BLP-complete-reformulation} and, again, only single-objective cuts for $\lambda=1$ are required.
\end{remark}

\section{Algorithm for Bi-objective Benders Decomposition}
\label{sec:algorithms}

Proposition \ref{prop:1} confirms that for a \gls{BLP} \eqref{eq:BLP}, a bi-objective Benders reformulation \eqref{eq:BLP-complete-reformulation}-$\{ \bar{\D}_r,\Lambda,\bar{\D}_p\}$ with a finite number of weighted cuts exists. To obtain such a reformulation one could iteratively solve
single-objective weighted-sum problems with Benders decomposition, for instance by applying a dichotomic approach (see also Section \ref{subsec:formulation_defs}). This scalarisation approach would yield a single-objective Benders reformulation \eqref{eq:LP-Benders} of each weighted problem $\lambda$-\eqref{eq:BLP} for some $\lambda \in [0,1]$. The cuts from all scalarised reformulations with $\lambda \in \Lambda$ are combined to obtain \eqref{eq:BLP-complete-reformulation}-$\{ \bar{\D}_r,\Lambda,\bar{\D}_p\}$ where the optimality cuts in \eqref{eq:LP-Benders} are converted to weighted cuts by replacing $\theta = \lambda \theta_1 + (1-\lambda) \theta_2$, see \eqref{eq:weighted_optimality_cut}.
Instead, we describe an algorithm in the following that integrates the generation of Benders feasibility and weighted optimality cuts into the parametric bi-objective simplex algorithm to solve \gls{BM} \eqref{eq:BLP-restated} and \gls{BS} \eqref{eq:BSub}, thereby not solving single-objective scalarised problems. The two approaches are compared in Section \ref{sec:Comparison}.

In this section, the proposed \gls{BBSA} is introduced.
The \gls{BBSA} works by iteratively identifying non-dominated points  of a Benders master problem \gls{BM} by applying bi-objective parametric simplex iterations, where the initial BM is a problem of the form~\eqref{eq:BLP-restated}, which is \eqref{eq:BLP-complete-reformulation}-$\{\emptyset,\emptyset,\emptyset \}$ obtained from~\eqref{eq:BLP-complete-reformulation} by removing all feasibility and optimality cuts. At any iteration of the algorithm, this Benders reformulation is made up of a subset of feasibility cuts $\bar{\D}_r$ and optimality cuts $\bar{\D}_p(\lambda)$ for $\lambda \in \Lambda$ that have been added up to the current iteration of the \gls{BBSA}, and is the problem we refer to as \gls{BM} \eqref{eq:BLP-complete-reformulation}-$\{ \bar{\D}_r,\Lambda,\bar{\D}_p\}$ in the context of the \gls{BBSA}. \gls{BM} \eqref{eq:BLP-complete-reformulation}-$\{ \bar{\D}_r,\Lambda,\bar{\D}_p\}$  is iteratively updated by adding all identified feasibility and optimality cuts.

We term an efficient solution $(\bar{\theta}_1,\bar{\theta}_2,\bar{y})$ with objective vector $\bar{z}$ of \gls{BM} \eqref{eq:BLP-complete-reformulation}-$\{ \bar{\D}_r,\Lambda,\bar{\D}_p\}$ \emph{explored} 
if a \gls{BS} \eqref{eq:BSub} has been solved for $\bar{y}$. It is termed \emph{unexplored} otherwise. An explored solution that remains a feasible solution of \gls{BM} \eqref{eq:BLP-complete-reformulation}-$\{ \bar{\D}_r,\Lambda,\bar{\D}_p\}$ is \emph{confirmed} as efficient solution.
In each iteration, the \gls{BBSA} identifies at most one unexplored solution $(\bar{\theta}_1,\bar{\theta}_2,\bar{y})$ of \gls{BM} \eqref{eq:BLP-complete-reformulation}-$\{ \bar{\D}_r,\Lambda,\bar{\D}_p\}$ and explores it, i.e.~applies bi-objective simplex to solve \gls{BS} \eqref{eq:BSub} with $\bar{y}$. If \gls{BS} \eqref{eq:BSub} is infeasible, $\bar{y}$ is removed from \gls{BM} \eqref{eq:BLP-complete-reformulation}-$\{ \bar{\D}_r,\Lambda,\bar{\D}_p\}$ by adding a feasibility cut of the following form:
\begin{equation}
\label{eq:feasibility_cut_form}
    (b- B \bar{y})^{\top} \pi_{r} \leqq 0,
\end{equation}
where $\pi_{r} \in \mathcal{D}_{r}$ and $\bar{\D}_r$ is updated to include $\pi_r$ (i.e.~$\bar{\D}_r = \bar{\D}_r \cup \{ \pi_r \} $). Otherwise, the values of $\theta_{1}$ and $\theta_{2}$ in \gls{BM} \eqref{eq:BLP-complete-reformulation}-$\{ \bar{\D}_r,\Lambda,\bar{\D}_p\}$ are modified, if necessary, by adding one or more weighted optimality cuts for $\pi_p^{\lambda^i} \in \D_p(\lambda^i)$, as explained later in this section, to obtain an updated \gls{BM} \eqref{eq:BLP-complete-reformulation}-$\{ \bar{\D}_r,\Lambda,\bar{\D}_p\}$ (i.e.~$\Lambda$ and $\bar{\D}_p$ are updated). The \gls{BBSA} terminates when there are no remaining unexplored non-dominated points of \gls{BM}, at which point the final \eqref{eq:BLP-complete-reformulation}-$\{ \bar{\D}_r,\Lambda,\bar{\D}_p\}$ is obtained.

The \gls{BBSA} investigates at most one unexplored non-dominated point of~\gls{BM} \eqref{eq:BLP-complete-reformulation}-$\{ \bar{\D}_r,\Lambda,\bar{\D}_p\}$ in each iteration.
It is therefore not necessary to fully solve \gls{BM} \eqref{eq:BLP-complete-reformulation}-$\{ \bar{\D}_r,\Lambda,\bar{\D}_p\}$ with the bi-objective simplex algorithm 
in each iteration of the \gls{BBSA}. Let us suppose that, at a given iteration, amongst all non-dominated points of \gls{BM} \eqref{eq:BLP-complete-reformulation}-$\{ \bar{\D}_r,\Lambda,\bar{\D}_p\}$, $z^{last}$ is the most recently explored (and therefore confirmed) point and that $(\theta^{last}_{1},\theta_{2}^{last},y^{last})$ is the corresponding efficient solution.
In the subsequent iteration, the bi-objective simplex algorithm starts from the explored non-dominated point $z^{last}$ and moves to a neighbouring unexplored non-dominated point.

\begin{algorithm}[htp]
\caption{\gls{BBSA} to solve bi-objective linear optimisation problems}
\label{alg:BBSA}

\newcommand\mycommfont[1]{\small\textcolor{darkgray}{#1}}
\SetCommentSty{mycommfont}

\SetAlgoLined
\setcounter{AlgoLine}{0}
\linespread{1.2}\selectfont
\SetKwRepeat{Do}{do}{while}
\SetKwRepeat{Repeat}{repeat}{until}
\KwIn{$A,\,B,\,D,\,C=(c^{1},c^{2}),\,F=(f^{1},f^{2}),\,b,\,d$}
\KwOut{$\mathcal{Z}_{nx}$, $\mathcal{X}_{ex}$}

$\bar{\mathcal{Z}}_{n} \leftarrow \{ \}$; $\bar{\mathcal{X}}_{e} \leftarrow \{ \}$; $\Lambda = \{1 \}$; $\bar{\lambda} \leftarrow 1$   \label{alg:initZX}

$[(\overline{\theta}_{1}, \overline{\theta}_{2} , \overline{y}),\bar{z},\bar{\D}_r,\bar{\D}_p(1)] \leftarrow$ \texttt{SolveSingleObjectiveBenders}($\bar{\lambda}$-\eqref{eq:BLP-complete-reformulation}-$\{ \emptyset,\emptyset, \emptyset \}$)\label{alg:single_solve}

\Do{$\ell > 0$}
{
$[ \D_r^{+}, \D_p^+, \Lambda^+ ] \leftarrow$ \texttt{Explore}($\bar{y}$) \label{alg:Solve BS} 

\eIf{$\D_r^+ \neq \emptyset $}{\label{alg:startcheck cuts} 
$\bar{\D}_r = \bar{\D}_r \cup \D_r^+ $ \tcp{Add feasibility cut $(b - B\bar{y})^{\top} \pi_r \leqq 0$ \eqref{eq:feasibility_cut_form} for $\D_r^+ = \{ \pi_r \}$}  \label{alg:feasibility cut}
}{ 

$\Lambda \leftarrow \Lambda \cup \Lambda^+ \cup \{ 0 \}$

\For{$a_i \in \Lambda^+$}{ \label{alg:StartWeightedOptCuts}
    $\bar{\D}_p (a_{i}) = \bar{\D}_p (a_{i}) \cup \{\pi_p^{a_{i}} \}$  \tcp{Add weighted optimality cut $(b-By)\pi^{\lambda}_{p} \leqq \lambda \theta_{1} + (1-\lambda)\theta_{2}$ \eqref{eq:weighted_optimality_cut} for $\lambda = a_i$ where $\pi^{\lambda}_{p} = \lambda \pi^{i}_{1} + (1- \lambda) \pi^{i}_{2}$} 
} \label{alg:EndWeightedOptCuts}

$\bar{\D}_p (0) = \bar{\D}_p (0) \cup \{\pi_p^{0} \}$ \tcp{Add weighted optimality cut for $\lambda=0$, $\pi_p^{0} = \pi^k_2$} \label{alg:weightedOptCut0} 

}\label{alg:end check cuts}

\tcp{\gls{BM} \eqref{eq:BLP-complete-reformulation}-$\{ \bar{\D}_r,\Lambda,\bar{\D}_p\}$ was updated by adding new cuts}

\If {first iteration}{ \label{alg:start_cheking lambda}
$[(\overline{\theta}_{1}, \overline{\theta}_{2} , \overline{y}),\bar{z}] \leftarrow$ \texttt{FindFirstEfficient}(\eqref{eq:BLP-complete-reformulation}-$\{ \bar{\D}_r,\Lambda,\bar{\D}_p\}$) \label{alg:firstEff}
} 

\If{$(\overline{\theta}_{1}, \overline{\theta}_{2} , \overline{y})$ feasible solution of \gls{BM} \eqref{eq:BLP-complete-reformulation}-$\{ \bar{\D}_r,\Lambda,\bar{\D}_p\}$}{ \label{alg:startStore} 

$\bar{\Z}_{n} \leftarrow \bar{\Z}_{n} \cup \{ \overline{z} \}$; $\bar{\X}_{e} \leftarrow \bar{\X}_{e} \cup \{(\overline{\theta}_{1}, \overline{\theta}_{2}, \overline{y} )\}$ \label{alg:storeFirstEff}

$(\theta_1^{last} ,\theta_2^{last},y^{last}) \leftarrow (\bar{\theta}_1,\bar{\theta}_2,\bar{y})$; $z^{last} \leftarrow \bar{z}$ \label{alg:initLast}
\tcp{Update last explored solution} 
} \label{alg:endStore} 

$\ell \leftarrow 1$; 
$(\theta_1^{\ell},\theta_2^{\ell},y^{\ell}) \leftarrow (\theta_1^{last},\theta_2^{last},y^{last}) $; $z^{\ell} \leftarrow z^{last}$ \label{alg:beforeStartPivots}

\Repeat{$(\theta_1^{\ell},\theta_2^{\ell},y^{\ell}),z^{\ell}$ unexplored and $z_2^{\ell} < z_2^{last}$}
{ \label{alg:startPivots}
$[(\theta_1^{\ell + 1},\theta_2^{\ell + 1},y^{\ell + 1}),stop] \leftarrow$ 
\texttt{OneParaSimplexPivot}($(\theta_1^{\ell},\theta_2^{\ell},y^{\ell})$, \eqref{eq:BLP-complete-reformulation}-$\{ \bar{\D}_r,\Lambda,\bar{\D}_p\}$) \label{alg:OnePivot}

\eIf{$!stop$} {
$z^{\ell + 1}_i \leftarrow \theta_i^{\ell+1} + (f^i)^{\top} y^{\ell+1}$ for $i = 1,2$ 

$\ell = \ell + 1$
} {
$\ell \leftarrow 0$; break
}
}\label{alg:endPivots}

$(\bar{\theta}_1,\bar{\theta}_2,\bar{y}) \leftarrow (\theta_1^{\ell} ,\theta_2^{\ell},y^{\ell})$; $\bar{z} \leftarrow z^{\ell}$ \tcp{Update next solution to explore} \label{alg:updateNext}
}

$[\Z_{nx},\X_{ex}] \leftarrow$ 
\texttt{FilterExtreme}($\bar{\Z}_{n},\bar{\X}_{e}$) \label{alg:discard not extreme}

\Return{$\Z_{nx},\X_{ex}$}
\end{algorithm}

Algorithm~\ref{alg:BBSA} provides an overview of the \gls{BBSA} algorithm. The input of Algorithm~\ref{alg:BBSA} is problem data, and its output is the set of non-dominated extreme points $\mathcal{Z}_{nx}$ for problem \eqref{eq:BLP} and \eqref{eq:BLP-complete-reformulation}-$\{ \bar{\D}_r,\Lambda,\bar{\D}_p\}$, and corresponding efficient solutions $\mathcal{X}_{ex}$ of \eqref{eq:BLP-complete-reformulation}-$\{ \bar{\D}_r,\Lambda,\bar{\D}_p\}$. Sets $\bar{\Z}_n$ and $\bar{\X}_e$ track non-dominated points and efficient solutions found while the algorithm runs (initialised in line \ref{alg:initZX}). Weight set $\Lambda$ and weight $\bar{\lambda} = 1$ are also initialised in line \ref{alg:initZX}.

 The \gls{BBSA} starts by finding the first non-dominated extreme  point, which is optimal for the first objective function, by solving \gls{BM} i.e.~$\bar{\lambda}$-\eqref{eq:BLP-complete-reformulation}-$\{ \emptyset,\emptyset, \emptyset \}$ for $\bar{\lambda} = 1$ and without cuts, for the first objective function using the standard Benders decomposition technique, see line \ref{alg:single_solve}, function \texttt{SolveSingleObjectiveBenders}. Solution $(\bar{\theta}_1, \bar{\theta}_2,\bar{y})$ and objective cost vector $\bar{z}$ are returned as well as the resulting sets of extreme rays and extreme points denoted $\bar{\D}_r$, $\bar{\D}_p(1)$ that form the initial cuts for problem \eqref{eq:BLP-complete-reformulation}-$\{ \bar{\D}_r,\Lambda, \bar{\D}_p \}$.

The \gls{BBSA} iterates between exploring a solution in line \ref{alg:Solve BS}, i.e.~solving the subproblem \gls{BS} \eqref{eq:BSub} with the parametric simplex method (if it is feasible) and adding cuts (lines \ref{alg:startcheck cuts}-\ref{alg:end check cuts}), storing solutions that are confirmed as efficient (lines \ref{alg:startStore}-\ref{alg:endStore}) and partially solving the updated \gls{BM} (lines \ref{alg:startPivots}-\ref{alg:endPivots}). In the following we describe these steps in detail.

In line \ref{alg:Solve BS}, a solution $\bar{y}$ is explored by applying function \texttt{Explore} (Algorithm \ref{alg:Explore}). Inputs to Algorithm \ref{alg:Explore} are $\bar{y}$ and the problem parameters required to define \gls{BS} \eqref{eq:BSub}. It returns either a feasibility cut $\D_r^+$ or weighted optimality cuts $\D_p^+, \Lambda^+$ to be added to \gls{BM}. At first the \gls{BLP} \gls{BS} \eqref{eq:BSub} is solved in line \ref{alg:SolveBLP} of Algorithm \ref{alg:Explore} where function \texttt{SolveBLP} applies the bi-objective parametric simplex \citep[e.g.][]{Ehrgott}, returning the set of $k$ efficient extreme points $\X_{BS} = \{ x^1, x^2, \ldots, x^k \} $ of \gls{BS} as well as the corresponding non-dominated extreme points $\Z_{BS} = \{ z^1_{BS}, z^2_{BS}, \ldots, z^k_{BS} \}$.
Set $\A = \{ a_1, a_2, \ldots, a_k, a_{k+1} \}$ with $a_{k+1}=0$ contains the weight limits such that for $ \lambda^i \in [a_{i+1}, a_i]$ each $x_{BS}^i \in \X_{BS}$ is an optimal solution of $\lambda^i$-\eqref{eq:BSub}. 
Set $\Pi = \{ (\pi_1^1, \pi_2^1), (\pi_1^2, \pi_2^2), \ldots, (\pi_1^k, \pi_2^k) \} $ is also returned and contains, for each extreme point $x^i \in \X_{BS}$, vectors of dual variables $(\pi_1^i, \pi_2^i)$ associated with the constraints of \gls{BS} \eqref{eq:BSub}, the values of which are determined by the respective objective function.
Vector $\pi_r$ is an unbounded extreme ray (or a null vector if none exists). If the problem is infeasible ($\X_{BS} = \emptyset$), only the extreme ray is returned in $\D_r^+$ in line \ref{alg:extremeray} and the others sets remain empty (lines \ref{alg:Dp0} and \ref{alg:L0}). If, on the other hand, the problem is feasible, $\D_r^+$ remains empty (line \ref{alg:Dr0}), and the weighted duals are returned as $\D_p^+$ (line \ref{alg:Dp}) and corresponding weights are $\Lambda^+$ (line \ref{alg:L}). We note that it is assumed the ordering of elements in $\X_{BS}$, $\Z_{BS}$, $\Pi$, $\A$ is maintained throughout the algorithm.
 
\begin{algorithm}[tb]
\caption{\texttt{Explore}$(\bar{y})$}
\label{alg:Explore}

\newcommand\mycommfont[1]{\small\textcolor{darkgray}{#1}}
\SetCommentSty{mycommfont}

\SetAlgoLined
\setcounter{AlgoLine}{0}
\linespread{1.2}\selectfont
\SetKwRepeat{Do}{do}{while}
\KwIn{$\bar{y}$ and (implicitly) $A,\,B,\,C=(c^{1},c^{2}),\,b$}
\KwOut{$\D_r^+$, $\D_p^+$, $\Lambda^+$}

$[\X_{BS}, \Z_{BS}, \A, \Pi, \pi_r] \leftarrow$ \texttt{SolveBLP}(\eqref{eq:BSub}, $\bar{y}$) \label{alg:SolveBLP}

    \eIf {$\X_{BS} ==  \emptyset $} {
    \tcp{Problem infeasible}
    $\D_r^+ \leftarrow \{ \pi_r \}$  \label{alg:extremeray}

    $\D_p^+ = \emptyset$ \label{alg:Dp0}

    $\Lambda^+ = \emptyset$ \label{alg:L0}
    
    } 
    {  
        \tcp{Problem feasible}

    $\D_r^+ = \emptyset$ \label{alg:Dr0}

    $\D_p^+ = \{ \pi_p^{a_i} = a_i \pi_1^i + (1- a_i) \pi_2^i \; | \; (\pi_1^i, \pi_2^i) \in \Pi , a_i \in \A \} $ \label{alg:Dp}
    
    $\Lambda^+ = \A$ \label{alg:L}
       
    }
    
    \Return{$[\D_r^+, \D_p^+, \Lambda^+]$}

\end{algorithm}

In Algorithm \ref{alg:BBSA} the cut sets are updated next. If the \gls{BS} \eqref{eq:BSub} was found to be infeasible in \texttt{Explore}($\bar{y}$), $\bar{\D}_r$ is updated by adding a feasibility cut for the identified extreme ray (line \ref{alg:feasibility cut}).
Otherwise the \gls{BS} \eqref{eq:BSub} was found to be feasible and, for each weight $a_i \in \Lambda^+$, a weighted optimality cut \eqref{eq:weighted_optimality_cut} is added (lines \ref{alg:StartWeightedOptCuts}-\ref{alg:EndWeightedOptCuts}),
where $\lambda = a_{i}$, and the weighted dual is $\pi^{\lambda}_{p} = \lambda \pi^{i}_{1} + (1-\lambda) \pi^{i}_{2}$. Note that an additional optimality cut for $\theta_{2}$ is generated in line \ref{alg:weightedOptCut0} for $\lambda=0$ in~\eqref{eq:weighted_optimality_cut}. {Here, $\pi^{\lambda}_{p} = \pi_{2}^{k}$}
since $(\pi_{1}^{k},\pi_{2}^{k} )$ is the vector of duals corresponding to the last non-dominated point in \gls{BS} \eqref{eq:BSub} for $\lambda \in [a_{k+1},a_{k}]$ where $a_{k+1}=0$. After completing lines \ref{alg:startcheck cuts}-\ref{alg:end check cuts} of Algorithm \ref{alg:BBSA} \gls{BM} \eqref{eq:BLP-complete-reformulation}-$\{ \bar{\D}_r,\Lambda, \bar{\D}_p \}$ has been updated with the new cuts.

It should be noted that the initial solution $(\bar{\theta}_1, \bar{\theta}_1, \bar{y})$ was obtained as a minimiser of the first objective in line \ref{alg:single_solve}, and is therefore feasible, hence weighted optimality cuts are added in the first iteration. The initial solution $(\bar{\theta}_1,\bar{\theta}_2,\bar{y})$ is not necessarily an efficient solution of  \eqref{eq:BLP-complete-reformulation}-$\{ \bar{\D}_r,\Lambda,\bar{\D}_p\}$ as the second objective value may not be captured correctly. Therefore the first efficient solution is found in the first iteration after the initial set of optimality cuts are added. This happens in line \ref{alg:firstEff}  where function \texttt{FindFirstEfficient} performs steps of the bi-objective parametric simplex method \citep[e.g.][]{Ehrgott} until the first efficient solution is found.

If solution $(\overline{\theta}_{1}, \overline{\theta}_{2} , \overline{y})$ remains a feasible solution of \gls{BM} \eqref{eq:BLP-complete-reformulation}-$\{ \bar{\D}_r,\Lambda,\bar{\D}_p\}$ (which will always be the case in the first iteration), it is stored in $\bar{\Z}_n, \bar{\X}_e$ (line \ref{alg:storeFirstEff}). The last explored solution $(\theta_1^{last}, \theta_2^{last},y^{last})$ with objective vector $z^{last}$ is updated in line \ref{alg:initLast}.

From $(\theta_1^{last}, \theta_2^{last},y^{last})$, pivot steps of the bi-objective simplex algorithm  are utilised to move to the first unexplored efficient solution $(\theta_1^{\ell}, \theta_2^{\ell},y^{\ell})$ with cost vector $z^{\ell}$ in lines \ref{alg:beforeStartPivots}-\ref{alg:endPivots}. This is done by setting index $\ell \leftarrow 1$, initialising $(\theta_1^{\ell}, \theta_2^{\ell},y^{\ell}) \leftarrow (\theta_1^{last}, \theta_2^{last},y^{last})$ in line \ref{alg:beforeStartPivots} and then conducting a series of single parametric bi-objective simplex pivots (function \texttt{OneParaSimplexPivot} in line \ref{alg:OnePivot}) \citep[e.g.][]{Ehrgott} each obtaining $(\theta_1^{\ell + 1}, \theta_2^{\ell + 1},y^{\ell+1})$ from $(\theta_1^{\ell}, \theta_2^{\ell},y^{\ell})$ for problem \gls{BM} \eqref{eq:BLP-complete-reformulation}-$\{ \bar{\D}_r,\Lambda, \bar{\D}_p \}$. This process continues until an unexplored efficient solution is found with decrease in the second objective function (line \ref{alg:endPivots}), or there are no more pivots (in which case function \texttt{OneParaSimplexPivot} returns $stop = true$ in line \ref{alg:OnePivot}) and therefore the minimiser of the second objective has been reached.
In line \ref{alg:updateNext} the next solution to explore is updated to $(\theta_1^{\ell}, \theta_2^{\ell},y^{\ell})$.

Finally, the non-extreme and weakly non-dominated points and their corresponding efficient solutions are discarded to keep only extreme solutions (function \texttt{FilterExtreme} in line \ref{alg:discard not extreme}) to obtain $\X_{ex}$ and $\Z_{nx}$.  It should be noted that the complete set of efficient extreme solutions $\X_{ex}$ obtained by the algorithm may not be identical to that of the non-decomposed problem $(\X_e)$ but they are equivalent, i.e.~the image of $\X_{ex}$ and $\X_e$ is the same set of non-dominated extreme points. While Algorithm \ref{alg:BBSA} does not explicitly track the ranges of weights associated with each efficient solution in $\X_{ex}$, these weights are easily obtained in the process of the bi-objective parametric simplex algorithm and can be tracked throughout the algorithm. This detail is omitted here to keep the presentation as concise as possible.

\medskip For a worked example of Algorithm \ref{alg:BBSA} for the problem from Example \ref{ex:BLP-needs-weighted-cuts}, refer to \ref{subsec:example}.

\subsection{Correctness of Algorithm \ref{alg:BBSA}}\label{sec:proof_algorithm}

This section addresses correctness and finiteness of the \gls{BBSA} (Algorithm \ref{alg:BBSA}). Lemma \ref{lemma:BBSA_finds_efficient} confirms that all solutions $\X_{ex}$ found by the \gls{BBSA} algorithm do indeed correspond to efficient solutions of \gls{BLP} \eqref{eq:BLP}. Proposition \ref{prop:allfound} confirms that a complete set of efficient extreme solutions is found, i.e.~that the set of non-dominated extreme points $\Z_{nx}$ obtained by the \gls{BBSA} matches that of \gls{BLP} \eqref{eq:BLP}. Finally, Proposition \ref{prop:finiteness} confirms that the algorithm terminates after a finite number of steps.

\begin{lemma} \label{lemma:BBSA_finds_efficient}

 A solution $(\bar{\theta}_1,\bar{\theta}_2,\bar{y}) \in \X_{ex}$ with image $\bar{z} \in \Z_{nx}$ identified by Algorithm \ref{alg:BBSA} always corresponds to an efficient solution of \eqref{eq:BLP}, i.e.~there exists an efficient solution of \eqref{eq:BLP} with non-dominated image in objective space that is identical to $\bar{z}$.
 
\end{lemma}
\proof
The solution with image $\bar{z}$, and corresponding solution $(\bar{\theta}_1,\bar{\theta}_2,\bar{y})$ is a solution of \gls{BM} \eqref{eq:BLP-complete-reformulation}-$\{ \bar{\D}_r,\Lambda,\bar{\D}_p\}$ identified in some iteration of the \gls{BBSA},
 and it is optimal for $\lambda \in [\alpha,\beta]$. 
 This means that lines \ref{alg:startcheck cuts}-\ref{alg:end check cuts} of Algorithm \ref{alg:BBSA} did not generate cuts rendering $(\bar{\theta}_1,\bar{\theta}_2,\bar{y})$
 infeasible. The solution $(\bar{\theta}_1,\bar{\theta}_2,\bar{y})$ therefore defines an optimal solution to the corresponding problem $\bar{\lambda}$-\eqref{eq:BLP}, $\bar{\lambda} \in [\alpha,\beta]$, as the corresponding weighted subproblem $\bar{\lambda}$-\eqref{eq:BSub} cannot yield an optimality cut for $\bar{\lambda}$ to remove this solution,
  as shown in the following.

If \gls{BS} \eqref{eq:BSub} is feasible, it  is fully solved in line \ref{alg:Solve BS} of Algorithm \ref{alg:BBSA} (with details in Algorithm \ref{alg:Explore}) with all weighted cuts added in lines \ref{alg:StartWeightedOptCuts}-\ref{alg:weightedOptCut0}. For each efficient extreme solution of \gls{BS}, there is an interval of weights $[a_{k+1},a_k]$ for which the efficient solution is an optimal solution for $\lambda$-BS, with $\lambda \in [a_{k+1},a_k]$. In lines \ref{alg:StartWeightedOptCuts}-\ref{alg:weightedOptCut0} a weighted cut is generated for each obtained efficient extreme solution of \gls{BS} and corresponding $\lambda \in \{ a_1, a_2, \ldots, a_{k+1} \}$, with $a_1=1$ and $a_{k+1}=0$. In particular, there cannot exist a weighted cut for $\bar{\lambda} \in [\alpha,\beta]$ that would render solution $(\bar{\theta}_1,\bar{\theta}_2,\bar{y})$ infeasible.
To see this choose any $\bar{\lambda} \in [\alpha,\beta]$, which will be a convex combination $\bar{\lambda} = \gamma a_{k+1} + (1-\gamma) a_k, \gamma \in [0,1]$ of the interval limits of one of the weight intervals $[a_{k+1},a_k]$ identified when solving \gls{BS}. In particular the following two weighted cuts were added for $a_{k+1}$ and $a_k$ (amongst others):
\begin{align}
    (b - By)^\top \pi_p^{a_{k+1}} & \leqq a_{k+1} \theta_1 + (1-a_{k+1}) \theta_2 & \label{eq:weightedcut1} \\
    (b - By)^\top \pi_p^{a_{k}} & \leqq a_{k} \theta_1 + (1-a_{k}) \theta_2 & \label{eq:weightedcut2}
\end{align}
If solution $(\bar{\theta}_1,\bar{\theta}_2,\bar{y})$ is feasible for weighted optimality cuts \eqref{eq:weightedcut1} and \eqref{eq:weightedcut2}, it is also feasible for their convex combination as the following is a convex combination of \eqref{eq:weightedcut1} and \eqref{eq:weightedcut2}:
\begin{equation}
    (b - By)^\top \underbrace{\left ( \gamma \pi_p^{a_{k+1}} + (1-\gamma) \pi_p^{a_{k}} \right )}_{\pi_p^{\bar{\lambda}}} \leqq \bar{\lambda} \theta_1 + (1-\bar{\lambda}) \theta_2, \label{eq:weighted_cut_convex_combination}
\end{equation}
where $\pi_p^{\bar{\lambda}}$ is an optimal dual solution of the weighted subproblem $\bar{\lambda}$-\eqref{eq:BSub} (see also definition of a weighted optimality cut in \eqref{eq:weighted_optimality_cut}). Please refer also to the proof of Remark \ref{remark1} in \ref{app:proof_remark1}, which uses a similar argument.

We can conclude that the solution $(\bar{\theta}_1,\bar{\theta}_2,\bar{y})$ is an optimal solution to the corresponding problem $\bar{\lambda}$-\eqref{eq:BLP}, $\bar{\lambda} \in [\alpha,\beta]$ by applying single-objective Benders to problem $\bar{\lambda}$-\eqref{eq:BLP}. The Benders solve can be initialised with the set of cuts $\{ \bar{\D}_r,\Lambda,\bar{\D}_p\}$ generated by Algorithm \ref{alg:BBSA}, where all weighted optimality cuts \eqref{eq:weighted_optimality_cut} are transformed into single-objective optimality cuts of the form $(b-By)^{\top} \pi^{\lambda}_p \leqq \theta$.
One subproblem solve of this single-objective Benders problem will immediately confirm optimality (and hence efficiency) of solution $(\bar{\theta},\bar{y})$ with $\bar{\theta} = \bar{\lambda} \bar{\theta}_1 + (1- \bar{\lambda}) \bar{\theta}_2$, as \gls{BS} \eqref{eq:BSub} has already verified that there does not exist a violated cut of the form  \eqref{eq:weighted_cut_convex_combination} that would render this solution infeasible.
\endproof

\begin{proposition} \label{prop:allfound}
    Algorithm~\ref{alg:BBSA} terminates with a bi-objective master problem \gls{BM} \eqref{eq:BLP-complete-reformulation}-$\{ \bar{\D}_r,\Lambda,\bar{\D}_p\}$ whose set of non-dominated extreme points $\Z_{nx}$ is equal to that of the original \gls{BLP}~\eqref{eq:BLP}.
\end{proposition}
\proof
Observe that at any iteration of the \gls{BBSA}, the master problem that is considered, \gls{BM}~\eqref{eq:BLP-complete-reformulation}-$\{ \bar{\D}_r,\Lambda,\bar{\D}_p\}$, constitutes a relaxation of \eqref{eq:BLP-complete-reformulation} with $\bar{\D}_r\subseteq \D_r$ and $\cup_{\lambda\in\Lambda}{\bar{\D}}_p(\lambda) \subseteq \cup_{\lambda\in[0,1]} \D_p^\lambda$.
Being a relaxation, any feasible $(\theta_1,\theta_2,y)$ of \eqref{eq:BLP-complete-reformulation} is also a feasible solution of \gls{BM} \eqref{eq:BLP-complete-reformulation}-$\{ \bar{\D}_r,\Lambda,\bar{\D}_p\}$; however, the reverse is not necessarily true. We note that, in particular, the set of efficient solutions $\X_e$ to \eqref{eq:BLP-complete-reformulation} remains feasible for \eqref{eq:BLP-complete-reformulation}-$\{ \bar{\D}_r,\Lambda,\bar{\D}_p\}$ at every iteration. Given a solution $(\theta_1,\theta_2,y) \in \X_e$ there does not exist a $\pi_p^\lambda \in \cup_{\lambda\in[0,1]}\D_p^\lambda$
 such that
$$ (b-By)^\top \pi_p^\lambda > \lambda \theta_1+(1-\lambda)\theta_2,$$
i.e., no violated, weighted optimality cut exists. The \gls{BBSA} cannot identify any dual extreme point of $\lambda$-\eqref{eq:BSub}, for a particular choice of $\lambda\in[0,1]$, that is not contained in $\cup_{\lambda\in[0,1]} \D_p^{\lambda}$, since the latter is an enumerated set. If it were possible to identify such a cut, then we would conclude that \eqref{eq:BLP-complete-reformulation} does not contain an enumerated set of extreme points to  $\lambda$-\eqref{eq:BSub} for $\lambda\in[0,1]$ or that $(\theta_1,\theta_2,y)$ is not an efficient extreme solution. As the \gls{BBSA} is guaranteed to move from a current explored efficient solution to a neighbouring efficient solution of \eqref{eq:BLP-complete-reformulation}-$\{ \bar{\D}_r,\Lambda,\bar{\D}_p\}$ following the pivoting steps of the bi-objective parametric simplex method, we conclude that the \gls{BBSA} cannot ``miss" an efficient extreme solution as it proceeds and Lemma~\ref{lemma:BBSA_finds_efficient} confirms that every solution found is efficient. To see this we briefly outline two cases. In Case 1, a particular non-dominated extreme point $\hat{z}$ (and corresponding efficient solution) is skipped, and its neighbouring points are correctly identified. This could be due to a cut that removes $\hat{z}$. However, this cut cannot exist (as also outlined in the second paragraph of the proof of Lemma \ref{lemma:BBSA_finds_efficient}). As we also know the parametric simplex method will correctly visit all non-dominated extreme points, $\hat{z}$ cannot be missed. In Case 2, non-dominated point $\hat{z}$ is missing because it is dominated by the objective vector $z^*$ of another solution that is incorrectly identified as efficient by the algorithm -- this contradicts Lemma \ref{lemma:BBSA_finds_efficient}, and would indicate that cuts exist that would remove $z^*$, hence $z^*$ could not have been confirmed as non-dominated. The cases are also illustrated in Figure \ref{fig:Prop2_cases} in \ref{sec:Appendix_Prop2}.

The \gls{BBSA} therefore finds a complete set of efficient extreme solutions with images corresponding to non-dominated extreme points of \gls{BLP} \eqref{eq:BLP}.
\endproof

\begin{proposition} \label{prop:finiteness}
    The \gls{BBSA} Algorithm \ref{alg:BBSA} applied to \gls{BLP} \eqref{eq:BLP} terminates in a finite number of steps.
\end{proposition}
\proof
    We start by making the following observations:
    \begin{enumerate}
        \item At each iteration of the \gls{BBSA} the master problem \gls{BM} defined by \eqref{eq:BLP-complete-reformulation}-$\{ \bar{\D}_r,\Lambda,\bar{\D}_p\}$ has a finite set of efficient extreme solutions, and so does each subproblem \gls{BS}, if it is feasible.
        \item Any feasible subproblem \gls{BS} has a finite number of efficient extreme solutions from which a finite number of weighted optimality cuts is generated.
        \item Both \gls{BM} defined by \eqref{eq:BLP-complete-reformulation}-$\{ \bar{\D}_r,\Lambda,\bar{\D}_p\}$ and \gls{BS} \eqref{eq:BSub} have a finite set of efficient extreme solutions that can be obtained by the bi-objective parametric simplex method in a finite number of steps (or basis pivots) assuming appropriate strategies to prevent cycling are employed.
    \end{enumerate}
    Initially, the \gls{BBSA} solves a single-objective Benders problem $\lambda$-\eqref{eq:BLP-complete-reformulation}-$\{\emptyset,\emptyset,\emptyset \}$ in a finite number of steps by adding a finite number of cuts. 
    At each iteration the \gls{BBSA} solve of \gls{BM} \eqref{eq:BLP-complete-reformulation}-$\{ \bar{\D}_r,\Lambda,\bar{\D}_p\}$ in lines \ref{alg:startPivots}-\ref{alg:endPivots} starts from the most recently confirmed efficient solution $(\bar{\theta}_1,\bar{\theta}_2,\bar{y})$ with image $\bar{z}$, which is an optimal solution to $\bar{\lambda}$-\gls{BM}. The bi-objective simplex method performs one or more basis pivots thereby iterating, in a finite number of steps, to an unexplored efficient solution of \gls{BM} \eqref{eq:BLP-complete-reformulation}-$\{ \bar{\D}_r,\Lambda,\bar{\D}_p\}$, $(\bar{\theta}_1^{\ell},\bar{\theta}_2^{\ell},\bar{y}^{\ell})$ with image $z^{\ell}$, associated weight $\lambda' \leqq \bar{\lambda}$, and strictly better second objective value (if it exists; otherwise the algorithm terminates).
    The corresponding efficient solution $(\bar{\theta}_1^{\ell},\bar{\theta}_2^{\ell},\bar{y}^{\ell})$ is then tested for feasibility and efficiency at the start of the next iteration of the BBSA, and two cases may occur after cuts are generated in lines \ref{alg:startcheck cuts}-\ref{alg:end check cuts}.

    Case 1: If solution $(\bar{\theta}_1^{\ell},\bar{\theta}_2^{\ell},\bar{y}^{\ell})$ remains feasible, and the algorithm successfully confirms this solution as efficient, then a new efficient solution with lower second objective value compared to $(\bar{\theta}_1,\bar{\theta}_2,\bar{y})$ is obtained.

    Case 2: If solution $(\bar{\theta}_1^{\ell},\bar{\theta}_2^{\ell},\bar{y}^{\ell})$ is neither feasible nor efficient, a violated feasibility cut or a finite number of weighted optimality cuts will be found to remove the solution from further consideration. Updating the master problem \gls{BM} \eqref{eq:BLP-complete-reformulation}-$\{ \bar{\D}_r,\Lambda,\bar{\D}_p\}$ with these new cuts will yield a new candidate unexplored solution with weight $\lambda'' \leqq \lambda'$, while the last explored efficient solution $(\bar{\theta}_1,\bar{\theta}_2,\bar{y})$ remains feasible.

    In Case 2, the algorithm returns to this  last explored efficient solution $(\bar{\theta}_1,\bar{\theta}_2,\bar{y})$ of \gls{BM} \eqref{eq:BLP-complete-reformulation}-$\{ \bar{\D}_r,\Lambda,\bar{\D}_p\}$ and performs a finite number of bi-objective parametric simplex iterations again to identify another candidate efficient solution in a finite number of steps (if it exists), as noted above.
   Case 2 will only occur a finite number of times until a neighbouring efficient extreme solution of $(\bar{\theta}_1,\bar{\theta}_2,\bar{y})$ is reached. This is because, firstly, there exists only a finite number of feasibility cuts. Secondly, only a finite number of weighted optimality cuts are added for a particular weight $\lambda'$:  The weight $\lambda'$ corresponds to optimal solutions  of the weighted sum problem for points on the facet between images $\bar{z}$ and $\bar{z}^{\ell}$ of $(\bar{\theta}_1,\bar{\theta}_2,\bar{y})$ and $(\bar{\theta}_1^{\ell},\bar{\theta}_2^{\ell},\bar{y}^{\ell})$ in objective space. This weight $\lambda'$ may remain unchanged between multiple subsequent iterations of Case 2. A finite set of weighted optimality cuts are added after solving \gls{BS} that, in particular, include cuts for weight $\lambda'$. Therefore, single-objective Benders theory, see \cite{benders1962}, confirms only a finite number of cuts will be added here by solving a finite number of subproblems to confirm $\lambda'$ is the facet's true weight in the original problem, hereby confirming solution $(\bar{\theta}_1,\bar{\theta}_2,\bar{y})$ is optimal for $\lambda \in [\lambda', \bar{\lambda}] $, and that $(\bar{\theta}_1^{\ell},\bar{\theta}_2^{\ell},\bar{y}^{\ell})$ is efficient. Alternatively, a cut alters the facet and its associated weight to $\lambda'' < \lambda'$ and the same argument can be made for this next $\lambda''$.
    
    In summary, to identify each (of the finitely many) efficient basic feasible solution in $\bar{\X}_e$ the \gls{BBSA} conducts a finite number of steps in Case 2 to eventually arrive at Case 1 thus concluding efficiency of a solution in a finite number of iterations.
\endproof

\section{Computational Experiments}
\label{sec:numerical}
Benders decomposition has been successfully applied to a wide range of  single-objective optimisation problems. We test the \gls{BBSA} on three different types of \gls{BLP} with decomposable structure. Sets of problem instances of varying size are tested for bi-objective variants of the \gls{FCTP} \citep[e.g.][]{zetina2019}, the \gls{PP} \citep[e.g.][]{elci2022} and the \gls{MKP} \citep[e.g.][]{angulo2016,maher2021}. Extensive numerical testing is conducted to analyse the performance of the proposed BBSA. We are primarily concerned with solving \glspl{BLP} at this stage and, as such, consider continuous relaxations of our problem instances in the testing. We compare the \gls{BBSA} to a dichotomic approach that is referred to as \gls{DBA} where a series of single-objective weighted-sum problems are solved with Benders decomposition.
The numerical experiments were carried out in MATLAB R2020b and performed on a 64-bit Windows 10 PC system with Intel Xeon W-2145 CPU 3.7 GHz processor and 32 GB RAM. The open-source lpsolve package is used whenever we need to solve weighted \gls{BM} and \gls{BS} optimisation problems.
An implementation of the bi-objective simplex algorithm is used to find a set of non-dominated extreme points in \gls{BM} and \gls{BS}.

Before we discuss results in Section \ref{sec:numerical}, the problem instances are briefly introduced in Sections \ref{sec:FCTP}-\ref{sec:BMKP}.
All instances are available in \url{https://bitbucket.org/araith/bbsa-instances}. The results of the \gls{BBSA} for the different instance types are discussed in Section \ref{sec:Roberti_instances_results}-\ref{sec:BMKP_results}. Finally, the comparison of runtime between the \gls{BBSA} and the \gls{DBA} is presented in Section \ref{sec:Comparison}.

\subsection{\gls{BFCTP}} \label{sec:FCTP}
In the \gls{FCTP} it is assumed that there are $m$ sources (e.g, factories) and $n$ sinks (e.g., customers). Each source can transfer commodities to any of the sinks with a per-unit transportation cost and a fixed charge. A solution to the \gls{FCTP}  involves satisfying the demand of each sink given the available capacity of each source. The first objective of the \gls{BFCTP} considered here minimises the total transportation cost (including variable and fixed costs) and the second one minimises the transportation time, where the fixed-cost component of the time objective is considered to be the preparation time required to set up the connection. The full formulation is shown in \ref{sec:Appendix_FCTP}.
We adapt instances from the literature:
\gls{FCTP} Data sets 1 and 2 are from~\citet{agarwal2012}, and data set 3 is from~\citet{roberti2015}.
Data set 1 contains 12 instances that have between 4 and 17 sources and between 3 and 43 sinks. Data sets 2 and 3 each have 30 instances, with $n=m=15$ and $n=m=20$, respectively.
The unit costs are scaled to approximately maintain a pre-defined ratio $\ratio$ between the total variable and fixed costs. It should be noted that this ratio was denoted $\theta$ in the original paper, but we use $\ratio$ here to avoid confusion with our use of variables $\theta_1,\theta_2$ in Section \ref{sec:full-Benders-reformulation-BLP}. Three different scaling factors are considered, $\ratio=0, 0.2, 0.5$. Instances with $\ratio=0$ are pure \gls{FCTP} in which all unit costs are set to zero. For $\ratio=0.2$ ($\ratio=0.5$), the unit costs are scaled so that their total constitutes $20\%$ ($50\%$) of the total costs.
The second objective function for \gls{BFCTP} is randomly generated with values in the same range as the first objective.

\subsection{\gls{BPP}} \label{sec:BPP}
In the \gls{PP}, there are $n$ potential assets, and the problem is how to distribute a unit investment among the $n$ assets to minimise cost. The original \gls{PP} model \citep{qiu2014} considers an uncertain return, where the overall return should be at least $r$. Here we consider the case where all $m$ possible return scenarios are included simultaneously. The formulation allows violation of some of the overall return constraints, and the total amount of violation is limited to be at most $k$. The full model is shown in \ref{sec:Appendix_BPP}. Since the original problem has a single-objective function, we randomly create a second one to obtain \glspl{BPP}.
In \citet{qiu2014}, $10$ different instances are generated where $m = 200$, and $n=20$ and the number of the violated constraints is $k=15$. We generated $100$ instances, where the total amount by which the return constraints are violated ranges from $k=1$ to $k=10$ with 10 different instances per $k$.

\subsection{\gls{BMKP}} \label{sec:BMKP}
The \gls{MKP} is the third type of benchmark instance and is based on \citet{angulo2016} who formulate a stochastic binary version of the \gls{MKP}. In the problem, three sets of items are packed into two sets of knapsacks with the aim of minimising overall cost, and with separate weight constraints associated with each knapsack.
Here we solve a relaxation of the \gls{MKP} in the deterministic case. \citet{angulo2016} generated 30 instances of identical problem dimensions, which are available in \citet{Ahmed2015}. We follow the same approach as for \gls{BFCTP} and randomly generate coefficients for the second objective function of the \gls{BMKP} instances. However, in the original instances, the subproblem objective function coefficients are all one. Instead, we randomly generated the second objective function coefficients with uniform distribution over a larger range of values. The model is shown in \ref{sec:Appendix_BMKP}.

\subsection{\gls{BFCTP} Results} \label{sec:Roberti_instances_results}
An overview of average statistics within each instance group is given in Table \ref{tab:Roberti}. The data set is shown in column DS, the instance group name is listed in column $n\_m$. Instance groups are considered for instances $15\_15$ and $20\_20$ with ten instances for each value of $\ratio$, and their averages are reported in Table \ref{tab:Roberti}, whereas for the instances in data set 1 only a single instance exists. The table lists the number of non-dominated extreme points $|\mathcal{Z}_n|$, and the number of feasibility, optimality and active cuts (FC, OC and AC). Column ``it" lists the number of iterations, tBBSA is the total runtime of the \gls{BBSA} (in sec), and columns Ben\%, BM\% and BS\% list the corresponding percentage of time taken by the initial single-objective Benders problem solves (Step 3 of Algorithm \ref{alg:BBSA}), \gls{BM} and \gls{BS}, respectively.

\begin{table}[t]
\centering
\caption{Summary of run statistics for \gls{BFCTP} instance groups. Averages are shown for instance groups $15\_15$ and $20\_20$ over sets of 10 instances for each value of $\ratio$.} \label{tab:Roberti}.
\begin{tabular}{clrrrrrrrrrr}
  \hline
\multicolumn{1}{c}{DS} & \multicolumn{1}{c}{$n\_m$} & \multicolumn{1}{c}{$\ratio$} & \multicolumn{1}{c}{$|\mathcal{Z}_n|$} & \multicolumn{1}{c}{FC} &\multicolumn{1}{c}{OC} & \multicolumn{1}{c}{AC} & \multicolumn{1}{c}{it} & \multicolumn{1}{c}{tBBSA} & \multicolumn{1}{c}{Ben\%} & \multicolumn{1}{c}{BM\%} & \multicolumn{1}{c}{BS\%} \\ 
  \hline
  1 & \phantom{4}4\_3 &- &  4 &  10 &  14 &  18 &   7 & 1.58 & 74.05 & 7.59 & 12.03 \\ 
   1 & 10\_10 &- & 28 &  76 &  85 &  78 &  63 & 5.17 & 37.33 & 25.73 & 32.11\\ 
   1 & 10\_12 &- & 22 &  87 &  85 & 131 &  87 & 6.48 & 26.85 & 29.32 & 38.27\\ 
   1 & 12\_12 &- & 33 & 175 & 161 & 155 & 168 & 19.99 & 14.36 & 53.83 & 28.01\\ 
   1 & 13\_13 &- & 33 & 178 & 164 & 161 & 171 & 25.37 & 13.72 & 53.37 & 29.09\\ 
   1 & 16\_16 &- & 54 & 408 & 311 & 260 & 406 & 195.35 & 3.56 & 79.18 & 14.70\\ 
   1 & \phantom{4}8\_32 &- & 60 & 308 & 453 & 110 & 327 & 205.06 & 4.60 & 78.78 & 15.13 \\ 
  1 & 10\_26 &- & 43 & 384 & 380 & 215 & 416 & 243.16 & 2.91 & 83.23 & 11.76\\ 
  1 & 12\_21 &- & 60 & 459 & 414 & 252 & 485 & 324.73 & 2.38 & 84.62 & 10.98\\ 
  1 & \phantom{4}6\_43 &- & 59 & 328 & 497 & 212 & 475 & 332.67 & 1.42 & 85.42 & 10.97 \\ 
  1 & 14\_18 &- & 68 & 531 & 390 & 254 & 534 & 344.30 & 2.43 & 84.03 & 11.23\\ 
  1 & 17\_17 &- & 46 & 907 & 839 & 670 & 987 & 1821.60 & 0.66 & 92.71 & 4.95\\\hline 
  2 & 15\_15 & 0.0 & 49.50 & 174.00 & 2.00 & 105.40 & 69.80 & 11.89 & 43.34 & 37.92 & 15.70\\ 
  2 & 15\_15 & 0.2 & 43.50 & 303.80 & 203.40 & 97.00 & 149.80 & 64.70 & 15.47 & 67.09 & 16.10\\ 
   2 & 15\_15 & 0.5 & 46.00 & 621.80 & 571.90 & 189.30 & 498.00 & 487.11 & 2.81 & 90.18 & 5.56\\ \hline 
  3 & 20\_20 & 0.0 & 75.50 & 495.80 & 2.00 & 247.00 & 82.50 & 73.68 & 48.47 & 39.61 & 9.55\\ 
  3 & 20\_20 & 0.2 & 79.50 & 621.10 & 482.50 & 201.00 & 304.50 & 701.73 & 6.02 & 71.70 & 21.10\\ 
  3 & 20\_20 & 0.5 & 72.50 & 1055.30 & 1411.00 & 188.50 & 942.40 & 5123.32 & 1.05 & 92.11 & 5.88\\ 
    \hline
\end{tabular}
\end{table}

\begin{figure}[tb]
    \centering
    \includegraphics[width=\textwidth]{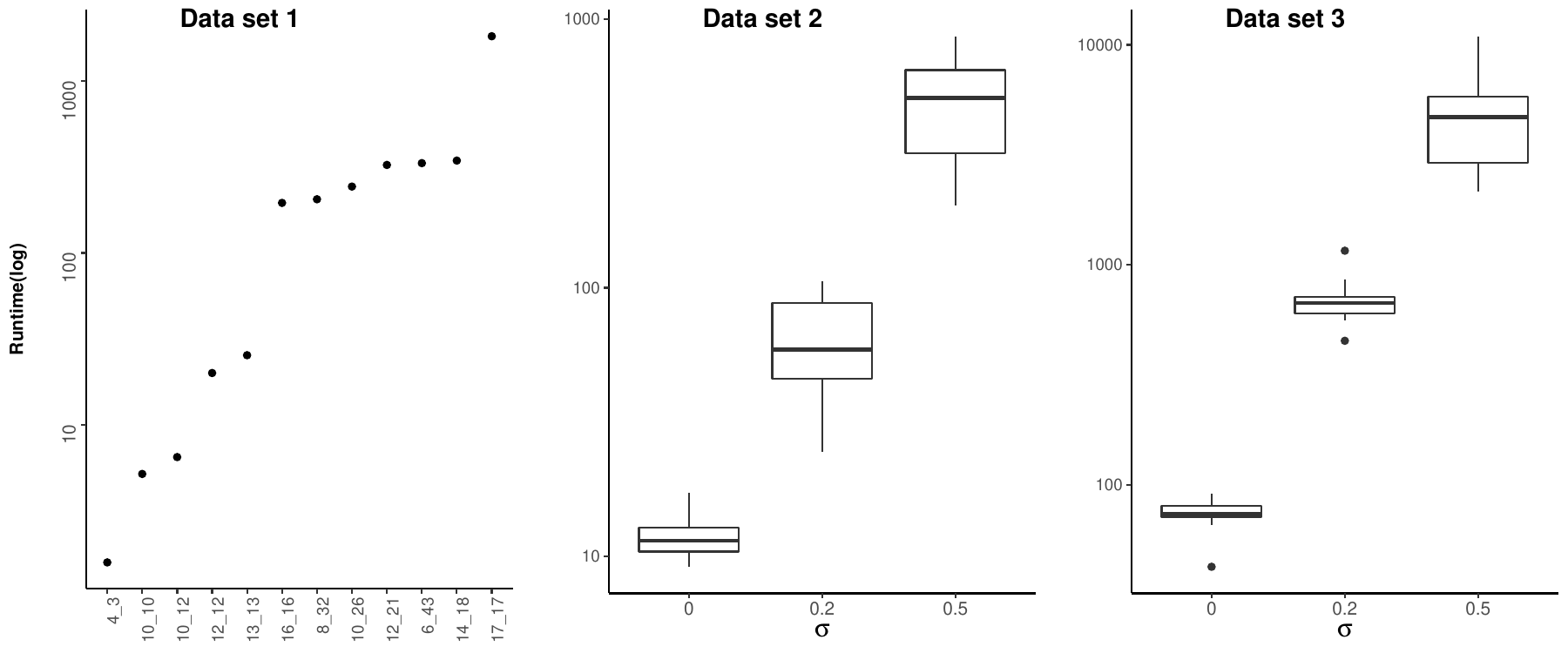}
    \caption{Runtimes per \gls{BFCTP} instance group}
    \label{fig:Roberti_times}
\end{figure}

The number of non-dominated extreme points found varies between 4 and 100 for the instances, where the average number of solutions in data sets 2 and 3 is fairly similar for different values of $\ratio$.
Increasing problem size with larger number of sources $n$ and sinks $m$ leads to higher runtime, and so does an increase in $\ratio$ as unit costs are only introduced into the problem when $\ratio>0$ and increase as $\ratio$ increases (Figure \ref{fig:Roberti_times}).
We observe that there is no strong trend relating the number of non-dominated extreme points with the runtime, except in the instances where $\ratio=0.5$.

\begin{figure}[tb]
    \centering
    \includegraphics[width=1\textwidth]{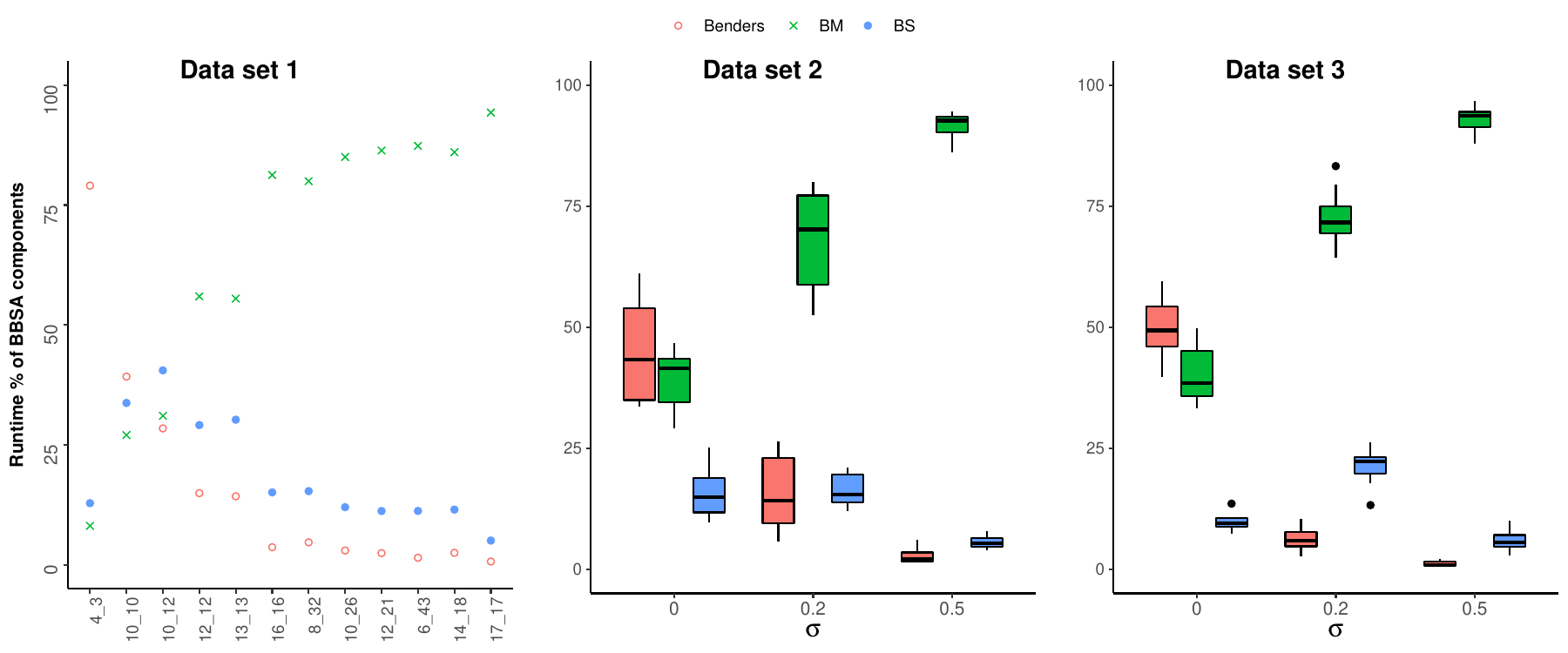}
    \caption{Runtime percentages of initial single-objective Benders, \gls{BM} and \gls{BS} of the BBSA.} \label{fig:Roberti_runtime_percentages}
\end{figure}

Figure \ref{fig:Roberti_runtime_percentages} shows the percentage of the runtime spent in different components of the \gls{BBSA} algorithm for the different \gls{BFCTP} data sets. In data set 1 the results appear to vary; however, it should be noted that some of the instances are quite easy to solve, and thus run quickly (see Table \ref{tab:Roberti}). For the more challenging instances runtime in \gls{BM} exceeds that in \gls{BS}, often quite significantly. For data sets 2 and 3 and $\ratio=0$, runtime of \gls{BM} and the initial single-objective Benders solve are similar, and exceed the percentage of time spent in \gls{BS}.
The percentage of runtime spent in \gls{BM} increases with increasing $\ratio$ as the number of optimality and feasibility cuts tends to increase with increasing master problem complexity.

Only a small portion of optimality cuts are generated in the initial single-objective Benders problem, with the exception of the instances with $\ratio=0$ in data sets 2 and 3. Here, two optimality cuts are generated in the initialisation and no further optimality cuts are generated as objective functions coefficients of \gls{BS} are zero. For data sets 2 and 3, on average 73.88\% and 78.07\% of feasibility cuts are generated in the initialisation, whereas only 45.52\%, on average, are generated for the instances in data set 1. For data set 3 with $\ratio=0$ most feasibility cuts are generated in the initialisation (98.65\% on average, whereas it is only 86.50\% for data~set~2).

In the \gls{BBSA} iterations, for data set 1, feasibility cuts are added in 52.5\% of all iterations, and
2.11 optimality cuts are generated per iteration in which optimality cuts are added. For data sets 2 and 3 feasibility cuts are added in 38.83\% and 30.92\% of all iterations. No optimality cuts are added for instances with $\ratio=0$, whereas on average 2.33 and 2.71 optimality cuts are added per iteration for data sets 2 and 3 when $\ratio > 0 $. The low number of optimality cuts per iteration is due to there being one or only a few non-dominated points in \gls{BS}. A large number of feasibility cuts is required here since the total supply equals the total demand in the \gls{BFCTP} instances, making it difficult to find a feasible solution when decomposing the problem.
For more challenging instances, such as data sets 2 and 3 with $\ratio=0.5$ a large number of cuts is generated but only a fairly small number of active cuts is needed for the Benders reformulation of the bi-objective problems.
There are between 31 and 378 active cuts in data sets 2 and 3 with an average of 2.85 cuts per non-dominated extreme point (5.05 for data set 1).

\begin{figure}
    \centering
    \includegraphics[width=1\textwidth]{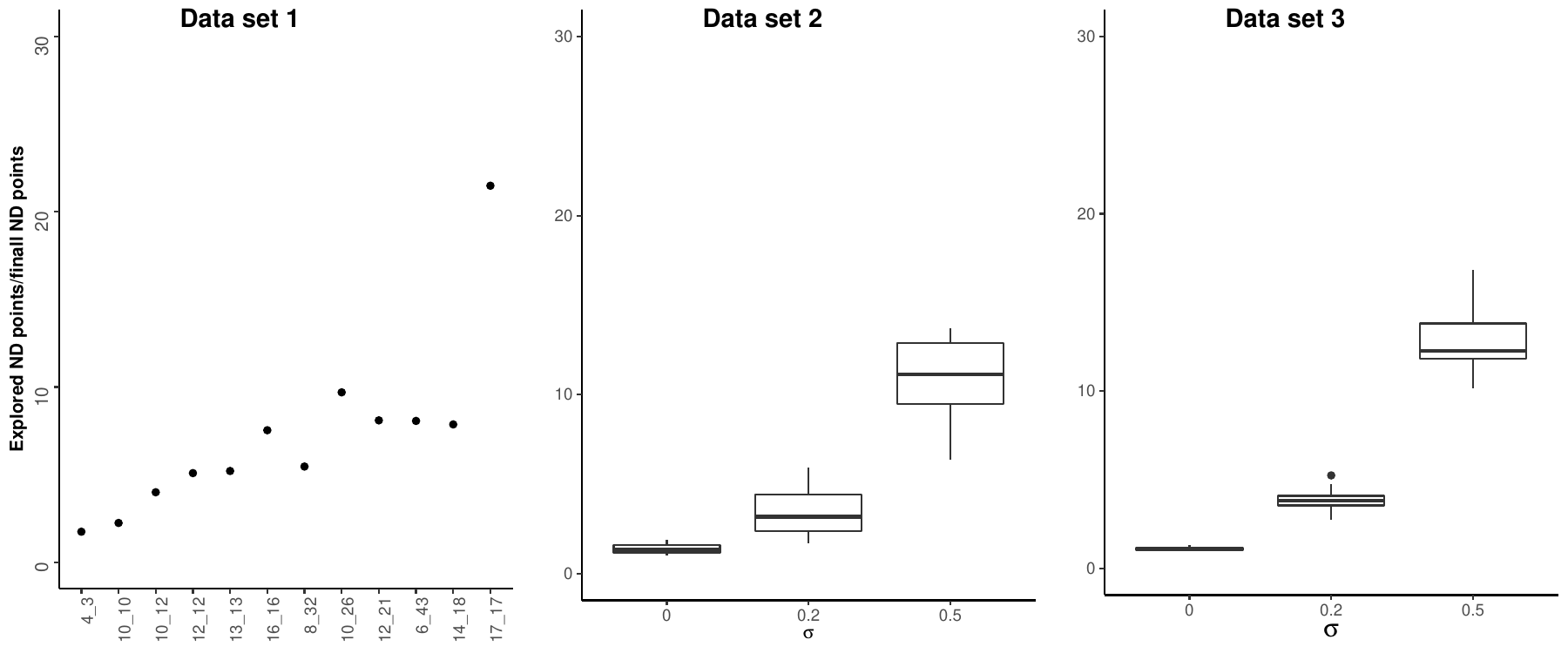}
    \caption{Ratio of explored non-dominated (ND) over final non-dominated extreme points.}
    \label{fig:Roberti_nondom}
\end{figure}

We can also compare the number of non-dominated points that were explored during the \gls{BBSA} to the number of final non-dominated points, as identified at the end of Algorithm \ref{alg:BBSA} (line \ref{alg:discard not extreme}), as shown in the ratio of explored over final non-dominated points in Figure \ref{fig:Roberti_nondom}.
We see a clear increase in this ratio based on problem complexity in \gls{BFCTP} data sets 2 and 3 for higher values of $\ratio$ indicating increasing problem complexity.

\begin{figure}[tb]
    \centering
    \includegraphics[width=1\textwidth]{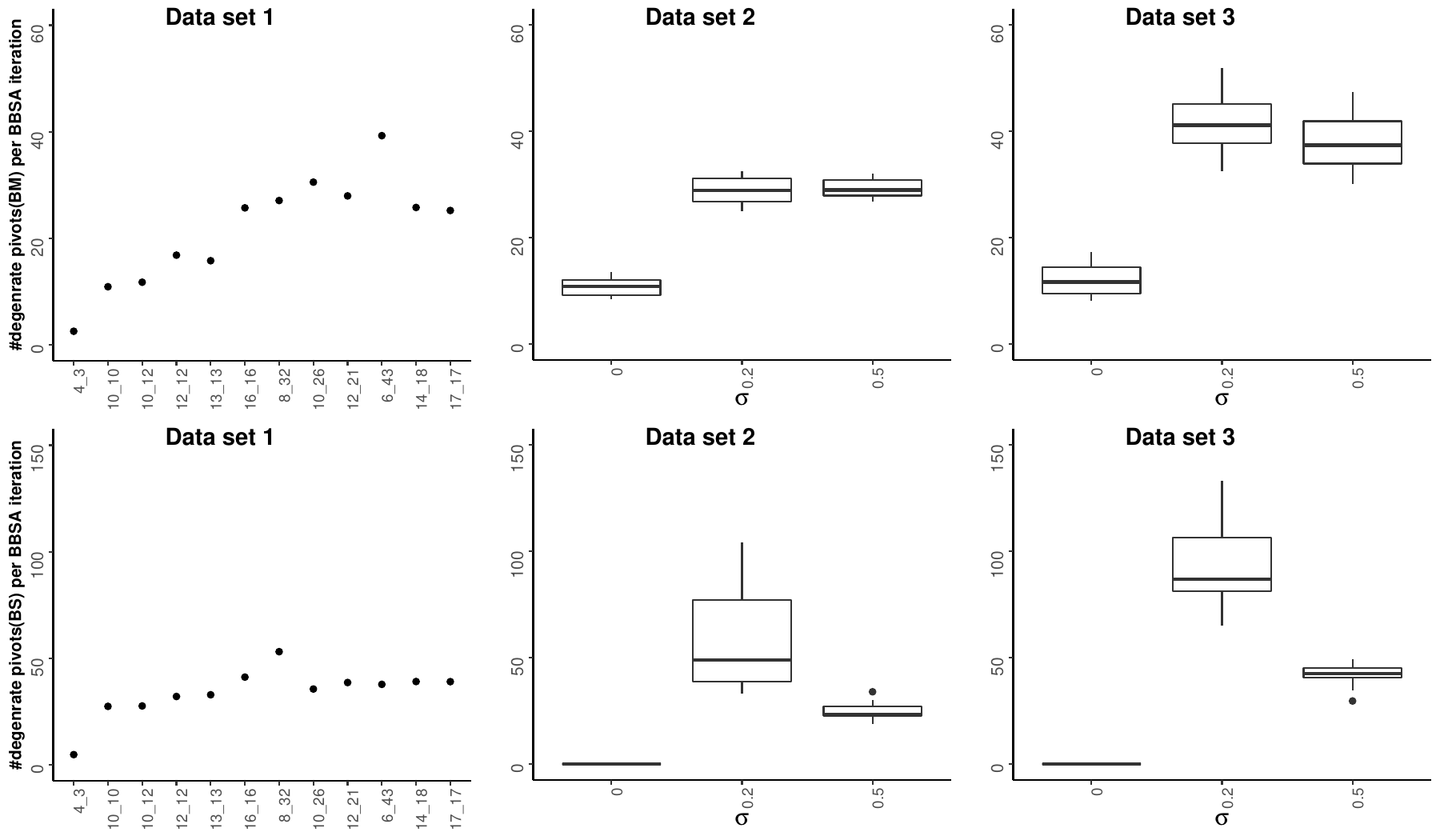}
    \caption{The average number of degenerate pivots per iteration of the \gls{BBSA} in which \gls{BM} is solved is shown in the top row; The average number of degenerate pivots in \gls{BS} in the bottom row.}
    \label{fig:Roberti_pivots}
\end{figure}

The considered problem instances suffer from degeneracy in \gls{BM} and \gls{BS}. Figure \ref{fig:Roberti_pivots} shows the average number of degenerate pivots in \gls{BM} per \gls{BBSA} iteration, and the average number of degenerate pivots per \gls{BS} solved. \gls{BS} generally has few non-dominated extreme points. Due to the degeneracy, many iterations are, however, needed to solve \gls{BS} with the bi-objective simplex algorithm. In a degenerate pivot, a variable enters the basis at value zero, and the objective values do not change. To confirm all non-dominated extreme points, hundreds of degenerate pivots may be required, and this is time-consuming.

\subsection{\gls{BPP} Results}
\begin{figure}[tbp]
    \centering
    \includegraphics[width=0.6\textwidth]{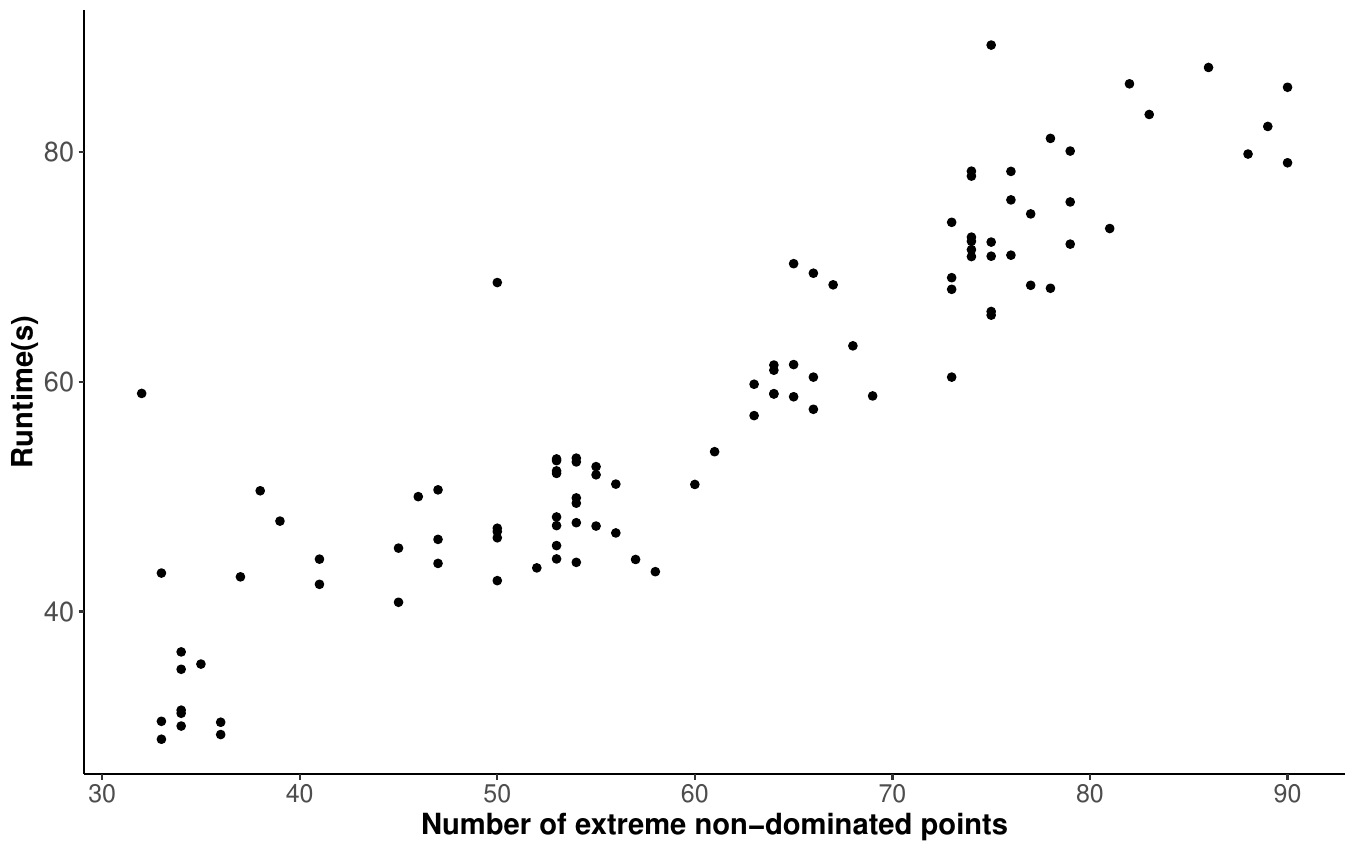}
    \caption{Number of non-dominated extreme points and runtime of \gls{BPP} instances.}
    \label{fig:nondom_BBP}
\end{figure}

There are $100$ instances in this class, $10$ for each value of $k$.
Table~\ref{tab:BPP} reports averaged statistics for each instance group. $\text{\gls{BPP}}\_k$ indicates the $k$-th group of \gls{BPP} instances where the total amount by which the return constraints can be violated is $k$.
The number of non-dominated extreme points varies from $32$ to $90$, with an average of $59.98$, which shows a fair spread in the number of non-dominated extreme  points. As $k$ increases, the number of non-dominated extreme  points also tends to grow. As expected, with increasing number of non-dominated extreme  points, the time taken to solve the problems also rises, as demonstrated in Figure \ref{fig:nondom_BBP}. However, the difference between the average time taken to solve instances with $k=1$ and instances $k=10$ is not large.

\begin{table}[t]
\centering
\caption{Summary of run statistics for \gls{BPP} instance groups, where the average values are shown.} \label{tab:BPP}.
\begin{tabular}{lrrrrrrrrr}
\hline
\multicolumn{1}{c}{$Instances$} & \multicolumn{1}{c}{$|\mathcal{Z}_n|$} & \multicolumn{1}{c}{FC} & \multicolumn{1}{c}{OC} & \multicolumn{1}{c}{AC} & \multicolumn{1}{c}{it} & \multicolumn{1}{c}{tBBSA} & \multicolumn{1}{c}{Ben\%} & \multicolumn{1}{c}{BS\%} & \multicolumn{1}{c}{BM\%} \\ \hline 
\gls{BPP}\_1& 43.40 & 519.50 & 1079.40 & 685.80 & 474.40 & 50.93 & 34.68 & 57.14 & 6.82 \\ 
\gls{BPP}\_2& 57.20 & 174.90 & 1436.20 & 667.20 & 654.40 & 53.16 & 16.59 & 74.51 & 7.21\\ 
\gls{BPP}\_3& 58.20 & 85.70  & 1537.50 & 551.70 & 706.60 & 53.90 & 11.39 & 79.33 & 7.51 \\ 
\gls{BPP}\_4& 62.22 & 37.44  & 1632.00 & 552.89 & 761.67 & 56.16 & 8.58  & 81.96 & 7.65 \\ 
\gls{BPP}\_5& 61.91 & 27.73  & 1652.27 & 392.27 & 770.18 & 56.73 & 8.13  & 82.33 & 7.74 \\ 
\gls{BPP}\_6& 62.60 & 14.00  & 1711.00 & 380.70 & 808.10 & 58.83 & 7.05  & 83.35 & 7.77 \\ 
\gls{BPP}\_7& 65.91 & 7.45   & 1821.36 & 387.18 & 878.91 & 63.54 & 6.18  & 84.14 & 7.84 \\ 
\gls{BPP}\_8& 60.33 & 4.89   & 1641.67 & 368.56 & 772.89 & 59.13 & 7.19  & 82.65 & 8.38\\ 
\gls{BPP}\_9& 63.40 & 3.60   & 1742.80 & 276.90 & 836.00 & 60.47 & 6.13  & 84.27 & 7.77 \\ 
\gls{BPP}\_10& 64.10 & 2.90   & 1806.60 & 312.30 & 877.00 & 63.23 & 6.01  & 84.37 & 7.78 \\ \hline 
\end{tabular}
\end{table}

\begin{figure}[tb]
    \centering
    \includegraphics[width=1\textwidth]{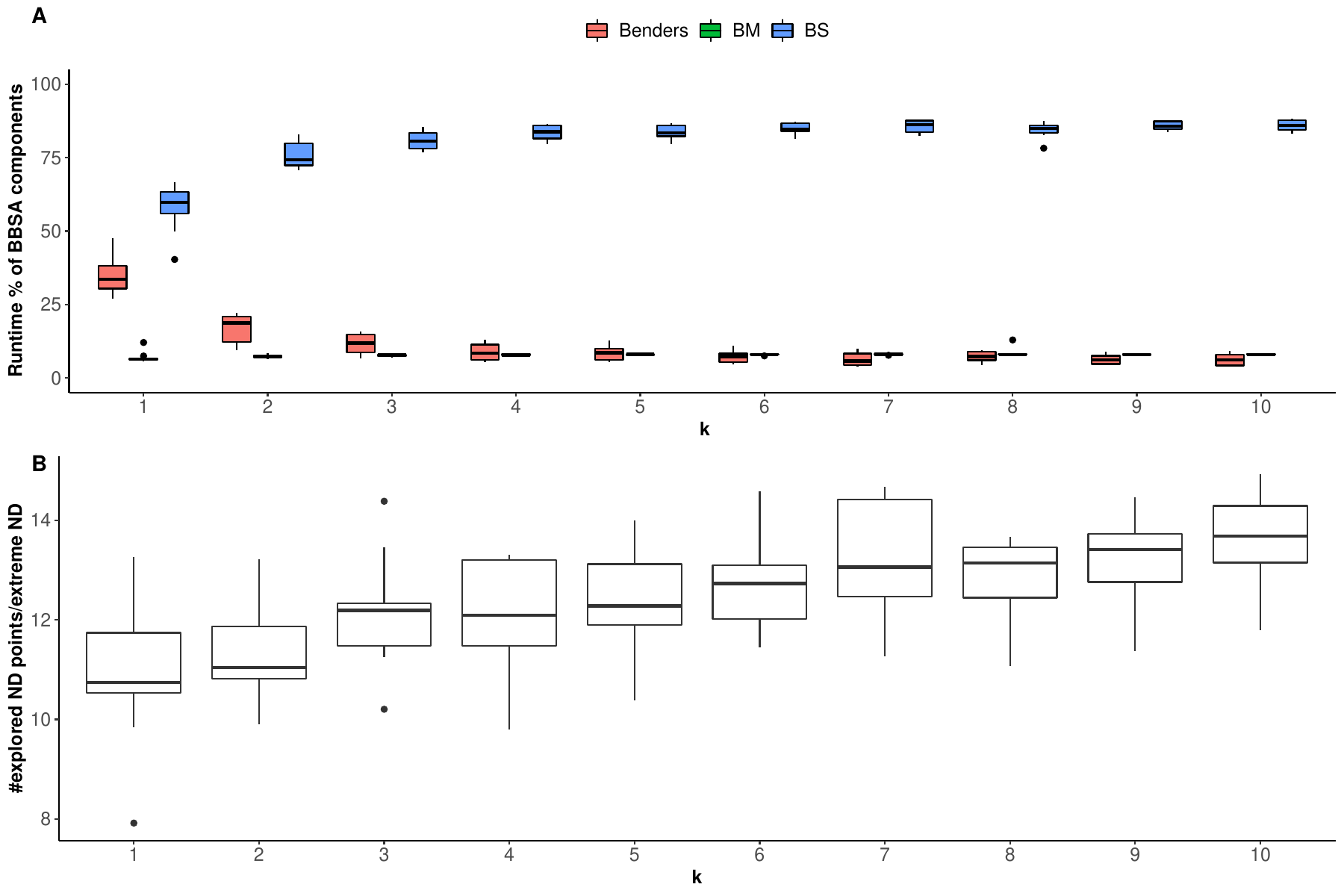}
    \caption{Comparison of the \gls{BBSA} components runtime for \gls{BPP} instances (A). The number of explored non-dominated points per extreme points for \gls{BPP} instances (B).}
    \label{fig:runtime_com_exp_BPP}
\end{figure}

The percentage of time taken in different components of the \gls{BBSA} for \gls{BPP} is depicted in Figure \ref{fig:runtime_com_exp_BPP}.A, which shows that solving \gls{BS} takes significantly more time than the other components of the \gls{BBSA}. The high number of $m=200$ return constraints in the \gls{BPP} model \eqref{eq:bi-objective-portfolio} belonging to \gls{BS} is the reason for comparably long solve time, whereas \gls{BM} contains only the one budget constraint of \eqref{eq:bi-objective-portfolio}. An interesting observation is that solving the initial single-objective Benders problem is more time consuming when $k$ is small, as shown in Figure \ref{fig:runtime_com_exp_BPP}.A.

\begin{figure}[tb]
    \centering
    \includegraphics[width=1\textwidth]{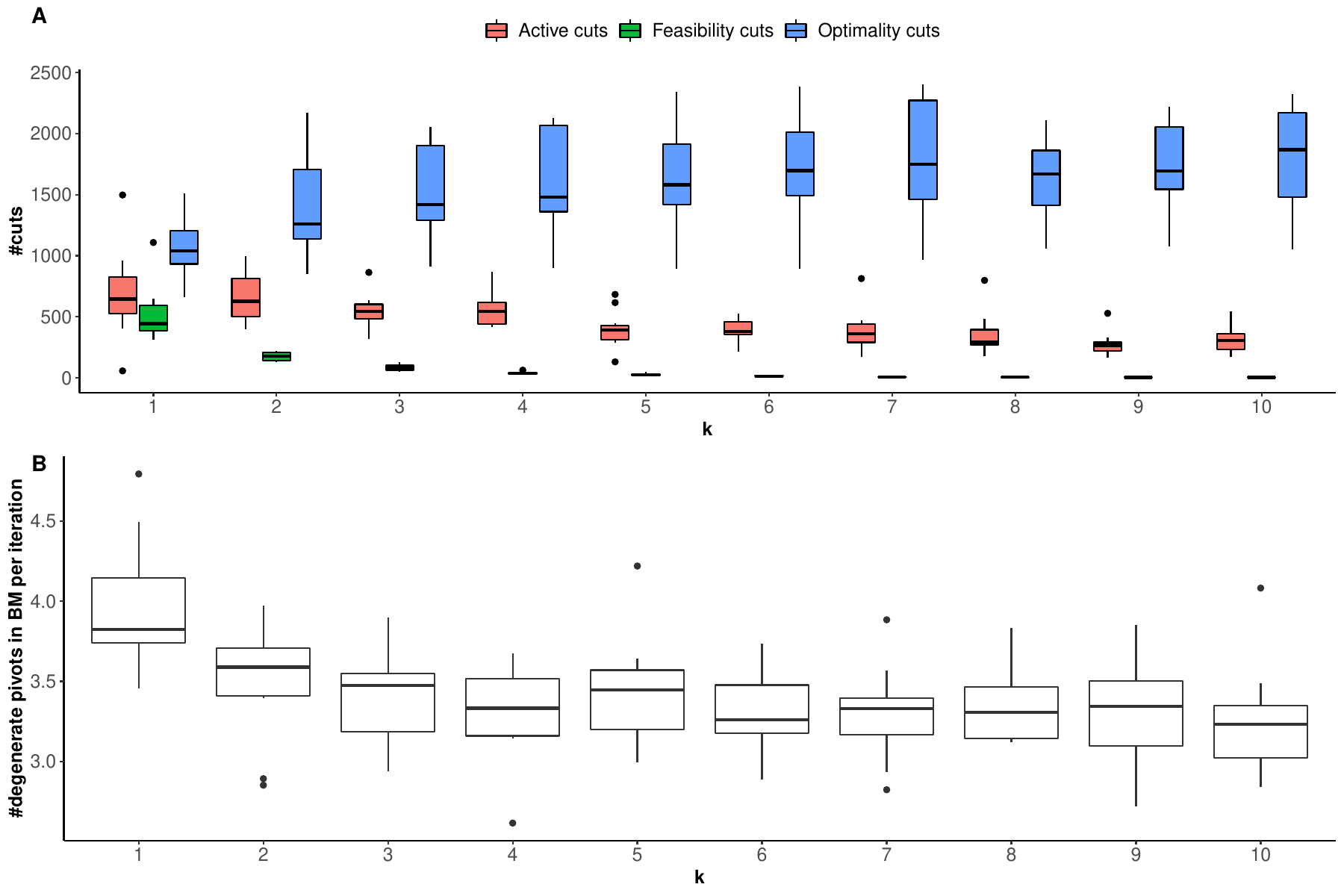}
    \caption{The number of cuts of the \gls{BBSA} for \gls{BPP} instances (A). The number of degenerate pivots of \gls{BM} for \gls{BPP} instances (B).}
    \label{fig:cuts_deg_BPP}
\end{figure}

Figure \ref{fig:runtime_com_exp_BPP}.B shows the average number of explored non-dominated points per extreme point for different groups of instances. We observe that the number of explored non-dominated points per non-dominated extreme point is spread from $3.17$ to $14.93$ with a mean of $12.36$. The average number of non-dominated points that are explored to identify an non-dominated extreme point increases slightly with increasing $k$.

It should also be noted that there is a single optimal solution in \gls{BS} for all \gls{BPP} instances.
When the values of the $y$ variables in \gls{BM} are fixed, the values of the $x$ variables can be uniquely derived and a single non-dominated extreme point is identified in \gls{BS}. Therefore, two optimality cuts are generated by \gls{BS} in each iteration corresponding to $\theta_1$ and $\theta_2$ with $\lambda=0$ and $\lambda=1$.
\gls{BPP} instances are therefore examples of Remark \ref{remark1}, where the optimal solution corresponding to $\lambda=1$ and $\lambda=0$ of \gls{BS} are the same (i.e.~there is only one non-dominated point in \gls{BS}). Therefore, to identify a set of non-dominated extreme  points, it is sufficient to generate single optimality cuts for $\theta_{1}$ and $\theta_{2}$ from \gls{BS} in each iteration.

Figure \ref{fig:cuts_deg_BPP}.A gives an overview of the number of cuts generated for \gls{BPP} instances. An interesting observation is that when $k$ is small, more feasibility cuts are generated as it is more likely to encounter infeasibility when the permitted constraint violation in \gls{BS} is small.
However, when $k$ is larger, the number of feasibility cuts generated from \gls{BS} drops significantly.
There is a large variation of between $1.78$ and $22.34$ active cuts, with an average of $7.99$ per non-dominated extreme point. More initial optimality and feasibility cuts are generated in \gls{BPP} instances in comparison with \gls{BFCTP}. $18.09\%$ of the total cuts are generated in the initial single-objective Benders solve (Algorithm \ref{alg:BBSA}, line \ref{alg:single_solve}), where $99.04\%$ and $13.47\%$ of the feasibility and optimality cuts are generated initially.

In contrast with \gls{BFCTP}, \gls{BPP} instances have fewer degenerate pivots as shown in Figure \ref{fig:cuts_deg_BPP}.B for \gls{BM} where the number of degenerate pivots per \gls{BBSA} iteration is between $2.39$ and $4.79$, with a mean of $3.41$, which is significantly less than for \gls{BFCTP} instances. The \gls{BS} has one non-dominated point, and there is no pivoting required to confirm this, and as such, there are no degenerate pivots in \gls{BS}.

\subsection{\gls{BMKP} Results} \label{sec:BMKP_results}
\begin{table}[tb]
\centering
\caption{Summary of minimum, mean and maximum run statistics for \gls{BMKP} instances.}
\label{tab:BMKP_min_max_avg}
\begin{tabular}{lrrrrrrrrrr}
\hline
 & \multicolumn{1}{c}{$|\mathcal{Z}_n|$} & \multicolumn{1}{c}{FC} & \multicolumn{1}{c}{OC} & \multicolumn{1}{c}{AC} & \multicolumn{1}{c}{it} & \multicolumn{1}{c}{tBBSA} & \multicolumn{1}{c}{Ben\%} & \multicolumn{1}{c}{BS\%} & \multicolumn{1}{c}{BM\%} & DBA \\ \hline
min  & 475.00 & 0.00 & 54509.00 & 13827.00 & 733.00 & 61.57 & 0.22 & 43.92 & 45.49 & 1.96 \\
mean & 584.77 & 0.00 & 67505.53 & 20657.43 & 858.83 & 73.26 & 0.32 & 47.43 & 49.57 & 3.19 \\
max & 705.00 & 0.00 & 81766.00 & 27152.00 & 958.00 & 83.50 & 0.56 & 51.61 & 53.16 & 5.06
\\\hline
\end{tabular}
\end{table}

Table \ref{tab:BMKP_min_max_avg} summarises the numerical results of \gls{BMKP}. We report each column's minimum, maximum and mean value.
For \gls{BMKP} instances, the number of non-dominated extreme  points (column $|\Z_n|$ in Table \ref{tab:BMKP_min_max_avg}) is much higher than for the \gls{BFCTP} and the \gls{BPP} instances.
There is no clear trend in the relationship between runtime and $|\Z_n|$.

The percentage of the runtime spent in different components of the \gls{BBSA} is different from other problem instance types.
For \gls{BMKP} instances, solving \gls{BM} and \gls{BS} takes a similar amount of time, and the time to solve the initial single-objective Benders problems is very low.

\begin{figure}[tp]
    \centering
    \includegraphics[width=1\textwidth]{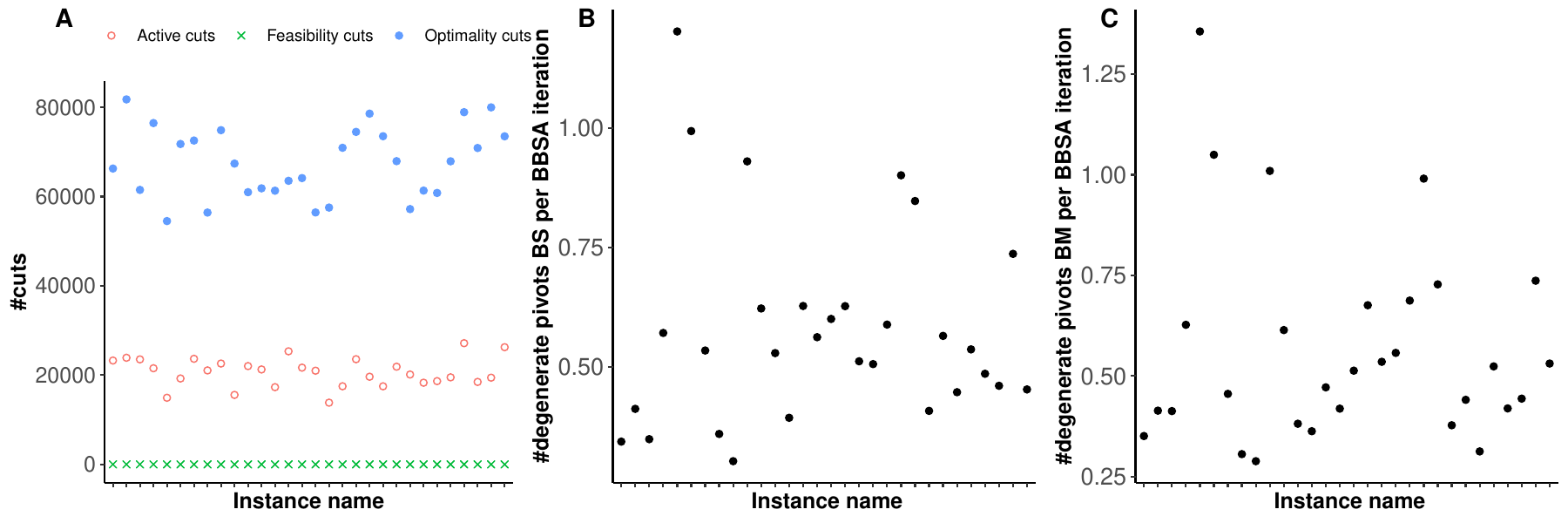}
    \caption{Analysis of \gls{BBSA} for \gls{BMKP} instances.}
    \label{fig:Statitical_Analysis_BMKP}
\end{figure}

Figure \ref{fig:Statitical_Analysis_BMKP}.A gives an overview of cut generation for \gls{BMKP} instances.  There are between $23.11\%$ and $39.89\%$ active cuts, with an average of $30.81\%$.  No feasibility cuts need to be generated for these instances. There are a large number of optimality cuts, only about $0.006\%$ of which are generated when solving the initial single-objective Benders problem.

The number of explored non-dominated points per non-dominated extreme  point is much lower compared to the other instance types. For \gls{BMKP} instances, only between $1.34$ and $1.63$ non-dominated points are required to be explored to identify an non-dominated extreme point, with an average of $1.47$.

Compared to the \gls{BPP} instances, there are even fewer degenerate pivots in these instances. From Figures \ref{fig:Statitical_Analysis_BMKP}.B and \ref{fig:Statitical_Analysis_BMKP}.C, it is clear that \gls{BS} has a similar number of degenerate pivots per \gls{BBSA} iteration  as the \gls{BM}. The number of degenerate pivots per \gls{BBSA} iteration in \gls{BS} is between $0.3$ and  $1.2$, with an average of $0.58$. For \gls{BM}, the number of degenerate pivots per iteration is between $0.29$ and $1.36$, with an average of $0.57$. These results show that, when the number of degenerate pivots in BS and BM is small enough, the \gls{BBSA} quickly identifies a set of non-dominated extreme points, even if there are many non-dominated extreme points.

Another interesting feature of the \gls{BMKP} instances is the number of weighted optimality cuts generated from \gls{BS}. In contrast to the \gls{BFCTP} and \gls{BPP} in stances, a significant number of weighted optimality cuts are generated in \gls{BMKP} instances. For \gls{BMKP} instances, between $65.86$ and $95.64$, with an average of $78.68$ weighted optimality cuts, are generated per iteration.

\subsection{Summary of experiments for \gls{BFCTP}, \gls{BPP} and \gls{BMKP} instances}

Several observations can be made from results obtained from the three different types of instances. For \gls{BFCTP} instances we observe that increasing problem size makes the problems more challenging to solve, and so does an increase of $\ratio$. Not surprisingly, the number of generated cuts generally increases for all problem types as the number of non-dominated points increases, although more feasibility cuts and active cuts are required for \gls{BPP} instances with lower $k$ that tend to have fewer non-dominated points.
It is interesting that \gls{BBSA} runtime spent in \gls{BM} and \gls{BS} varies by instance type. Runtime spent in \gls{BM} is higher for \gls{BFCTP} instances, and it is higher for \gls{BS} in \gls{BPP} instances, whereas both are similar for \gls{BMKP} instances. Runtime performance appears to depend on how many points need to be explored by the \gls{BBSA} compared to the number of final non-dominated extreme points. Runtime also depends on problem degeneracy where we observe that a high number of degenerate pivots are required for \gls{BFCTP} instances, whereas the other two instance types have much fewer degenerate pivots. In our experiments the highest runtimes are observed for the most challenging \gls{BFCTP} instances.
\gls{BS} has few efficient solutions for \gls{BFCTP} instances, often only one, meaning that few weighted optimality cuts are added in most iterations. In contrast, \gls{BMKP} instances have many efficient solutions in \gls{BS} and therefore many weighted optimality cuts.

\subsection{Comparison of the \gls{BBSA} and \gls{DBA}}\label{sec:Comparison}
We now compare the performance of the proposed \gls{BBSA} with \gls{DBA}.
The idea of \gls{DBA} is to iteratively  solve  weighted-sum problems \eqref{eq:BLP} with the standard (single-objective) Benders decomposition technique where the weights of $\lambda$ are identified from the dichotomic technique for \glspl{BLP} \citep[e.g.][]{Aneja,CCS79} to iteratively identify all non-dominated extreme  points. In the following, we compare the runtime (in seconds) of the proposed \gls{BBSA} with the \gls{DBA} on the \gls{BFCTP}, \gls{BPP}, and \gls{BMKP} instances.

Tables \ref{tab:RobertiDBA}-\ref{tab:BMKP_min_max_avg_DBA} give a summary of the DBA runs for our three instance types. Table \ref{tab:RobertiDBA} has the same first six columns as Table \ref{tab:Roberti}, and Tables \ref{tab:BPP_DBA} and \ref{tab:BMKP_min_max_avg_DBA} have the same first four columns as Tables \ref{tab:BPP} and \ref{tab:BMKP_min_max_avg}. The subsequent columns are the number of DBA iterations (itDBA), that is the number of times single-objective Benders decomposition was called to solve a weighted problem (note that $\text{itDBA} = 2|\Z_n|-1$); the DBA runtime (tDBA); the ratio of DBA and \gls{BBSA} runtime (tDBA/tBBSA) and the average DBA runtime per iteration (tDBA per it). If the runtime ratio tDBA/tBBSA exceeds 1, then \gls{BBSA} runs faster than DBA. Tables \ref{tab:BPP_DBA} and \ref{tab:BMKP_min_max_avg_DBA} show averages.

Solving our problem instances with the proposed \gls{BBSA} is more effective than iteratively solving single-objective weighted-sum standard Benders decomposition problems, as the columns tDBA/tBBSA in Tables \ref{tab:RobertiDBA}-\ref{tab:BMKP_min_max_avg_DBA} show very clearly, where the run-time ratio mostly exceeds $1$, often quite significantly. Some \gls{BFCTP} instances are challenging and a time limit of $10000$ seconds is applied to \gls{DBA}. The \gls{DBA} did not find all non-dominated extreme points in two instances of \gls{BFCTP} data set 3. There are only a few small instances where \gls{DBA} outperforms the \gls{BBSA}.

\begin{table}[t]
\centering
\caption{Summary of DBA run statistics for \gls{BFCTP} instance groups. Averages are shown for instance groups $15\_15$ and $20\_20$ over sets of 10 instances for each value of $\ratio$.} \label{tab:RobertiDBA}.
\begin{tabular}{clrrrrrrrr}
  \hline
\multicolumn{1}{c}{DS} & \multicolumn{1}{c}{$n\_m$} & \multicolumn{1}{c}{$\ratio$} & \multicolumn{1}{c}{$|\mathcal{Z}_n|$} & \multicolumn{1}{c}{FC} &\multicolumn{1}{c}{OC} & \multicolumn{1}{c}{itDBA} & \multicolumn{1}{c}{tDBA} & \multicolumn{1}{c}{tDBA/tBBSA} & \multicolumn{1}{c}{tDBA per it} \\
  \hline
 1 &\phantom{4}4\_3& -&4&60&27&7&1.20&0.76&0.17 \\
1 &10\_10& -&28&2324&317&55&27.83&5.38&0.51 \\
1&10\_12& -&22&2189&305&43&30.53&4.71&0.71 \\
1&12\_12& -&33&5756&631&65&90.69&4.54&1.40 \\
1&13\_13& -&33&6240&588&65&112.82&4.45&1.74 \\
1&16\_16& -&54&19258&1419&107&596.94&3.06&5.58 \\
1&\phantom{4}8\_32& -&60&14012&1919&119&445.38&2.17&3.74 \\
1&10\_26& -&43&11269&1234&85&358.34&1.47&4.22 \\
1&12\_21& -&60&16903&1753&119&495.90&1.53&4.17 \\
1&\phantom{4}6\_43& -&59&15148&1799&117&495.94&1.49&4.24 \\
1&14\_18& -&68&18157&2758&135&854.00&2.48&6.33 \\
1&17\_17& -&46&21277&2118&91&767.91&0.42&8.44 \\ \hline
2&15\_15&0.0&49.50&8050.10&96.60&98.00&190.55&16.02&1.94 \\
2&15\_15&0.2&43.50&11054.50&799.90&86.00&259.00&4.00&3.01 \\
2&15\_15&0.5&46.00&13377.50&1289.70&91.00&344.41&0.71&3.78 \\ \hline
3&20\_20&0.0&75.50&36387.60&146.80&150.00&3653.32&49.59&24.36 \\
3&20\_20&0.2&79.50&40932.30&1423.70&158.00&4204.77&5.99&26.61 \\
3&20\_20&0.5&72.10&83439.40&4226.40&143.20&5595.56&1.09&39.08 \\ \hline
\end{tabular}
\end{table}

\begin{table}[t]
\centering
\caption{Summary of DBA run statistics for \gls{BPP} instance groups, where the average values are shown.} \label{tab:BPP_DBA}.
\begin{tabular}{lrrrrrrr}
\hline
\multicolumn{1}{c}{$Instances$} & \multicolumn{1}{c}{$|\mathcal{Z}_n|$} & \multicolumn{1}{c}{FC} & \multicolumn{1}{c}{OC} & \multicolumn{1}{c}{itDBA} & \multicolumn{1}{c}{tDBA} & \multicolumn{1}{c}{tDBA/tBBSA} & \multicolumn{1}{c}{tDBA per it} \\ \hline
Bi-PP\_1&43.40&42588.00&16115.57&85.80&938.28&18.42&10.94\\
Bi-PP\_2&57.20&20062.20&23074.50&113.40&533.47&10.04&4.70\\
Bi-PP\_3&58.20&8848.60&24774.50&115.40&339.78&6.30&2.94\\
Bi-PP\_4&62.22&4925.40&26453.10&123.44&296.74&5.28&2.40\\
Bi-PP\_5&61.91&2678.00&27605.60&122.82&259.27&4.57&2.11\\
Bi-PP\_6&62.60&1532.40&28169.10&124.20&238.79&4.06&1.92\\
Bi-PP\_7&65.91&932.00&28790.80&130.82&239.91&3.78&1.83\\
Bi-PP\_8&60.33&611.50&28591.10&119.66&232.94&3.94&1.95\\
Bi-PP\_9&63.40&463.10&30529.50&125.80&238.28&3.94&1.89\\
Bi-PP\_10&64.10&341.50&30521.30&127.20&244.87&3.87&1.93\\
\hline
\end{tabular}
\end{table}

\begin{table}[tb]
\centering
\caption{Summary of minimum, mean and maximum \gls{DBA} run statistics and \gls{BMKP} instances.}
\label{tab:BMKP_min_max_avg_DBA}
\begin{tabular}{lrrrrrrrrrr}
\hline
 & \multicolumn{1}{c}{$|\mathcal{Z}_n|$} & \multicolumn{1}{c}{FC} & \multicolumn{1}{c}{OC} & \multicolumn{1}{c}{itDBA} & \multicolumn{1}{c}{tDBA} & \multicolumn{1}{c}{tDBA/tBBSA} & \multicolumn{1}{c}{tDBA per it} \\ \hline
min&475.00&0.00&5546.00&949.00&162.05&1.96&0.15\\
avg&584.77&0.00&8130.50&1168.53&232.86&3.19&0.20\\
max&705.00&0.00&13191.00&1409.00&387.56&5.06&0.31\\\hline
\end{tabular}
\end{table}

Depending on the type of instance, different factors affect \gls{DBA} performance. For \gls{BFCTP} instances, degeneracy makes solving problems with Benders decomposition slow in general, and this applies also to the single-objective weighted Benders problems solved by \gls{DBA}. We also observe that \gls{DBA} generates many more feasibility and optimality cuts than the \gls{BBSA}. In making this comparison it should be noted that different iterations of \gls{DBA} may generate the same cut and we might hence be counting some cuts multiple times, whereas cuts counted for the \gls{BBSA} are always unique cuts. For instance, in Data Set 2 of the \gls{BFCTP} instances, between 21.5 and 42.26 times as many feasibility cuts are generated by \gls{DBA} compared to the \gls{BBSA}, and this ratio is significantly higher for Data Set 3. \gls{DBA} also generates more optimality cuts although the ratio is lower (between 2.25 and 3.93 for instances with $\sigma > 0$). Both approaches would be similarly affected by degeneracy.

For \gls{BPP} to find a non-dominated extreme  point, lots of standard Benders decomposition iterations are required, particularly when $k$ is small as can be seen in Figure \ref{fig:runtime_com_exp_BPP}.A and Table \ref{tab:BPP} where the initial single-objective Benders problem solve makes up a large component of the run-time for small $k$, as well as in the time of \gls{DBA} per iteration (tDBA per it) in Table \ref{tab:BPP_DBA}. \gls{DBA} experiences large runtimes per iteration, especially in the more challenging instances with smaller values of $k$ where large numbers of cuts are generated across the iterations (Table \ref{tab:BPP_DBA}).

Finally, for the \gls{BMKP} instances, we observe that \gls{DBA} generates significantly fewer optimality cuts when compared with \gls{BBSA} -- only between 10.17\% and 16.13\%. \gls{DBA} also generates under 10 optimality cuts per iteration, hence converging to a solution for each single-objective problem quickly. Despite generating fewer cuts overall (and in each iteration), \gls{DBA} does show longer overall runtime. This is because the problem has many non-dominated extreme points. If there are $p$ non-dominated extreme points, \gls{DBA} needs to solve almost twice as many ($2p-1$) weighted single-objective problems with Benders decomposition. Even though individual problems solve fast, runtime accumulates here. Furthermore, the \gls{BBSA} is implemented to only retain currently active cuts to keep the size of the master problem \gls{BM} small (hence enabling us to solve it quickly in each iteration), and the subproblems \gls{BS} also solve very quickly due to low degeneracy.

\section{Conclusion and Future Work} \label{sec:conclusion}
In this paper, we have developed the necessary theory for applying Benders decomposition in a bi-objective setting. In particular, we present the Benders reformulation of a generic \gls{BLP} and argue that a finite subset of the constraints is sufficient to generate a complete set of efficient extreme solutions. Furthermore,  we discuss the need for so-called weighted optimality cuts and propose an iterative procedure that combines the bi-objective simplex algorithm with the cut generation procedure of the standard Benders decomposition algorithm.   We have shown that the proposed methodology can obtain a complete set of efficient extreme solutions and the corresponding set of non-dominated extreme  points to a \gls{BLP}.  Using sets of diverse problem instances, we have also identified factors that affect the performance of the proposed algorithm and showed that it performs well when compared to a dichotomic weighting approach combined with standard (single-objective) Benders decomposition.

There are several topics that can be explored in future research. Firstly, the results of \gls{BBSA} have shown that, as problem size increases, runtime of the \gls{BBSA} can increase quite quickly. It could therefore be beneficial to investigate ways to improve the performance of the Benders decomposition approach. Acceleration techniques for the single-objective case can be adapted to the bi-objective setting. Secondly, methods that exploit the bi-objective nature of the problem are also of interest. A bi-directional approach that simultaneously explores non-dominated points from left to right and right to left is one such possibility. So too are methods that can improve the accuracy of $\lambda$ weights considered. The bi-objective simplex algorithm determines the next $\lambda$ weight for the problem to which it is being applied. Investigating whether or not one can accurately determine the next $\lambda$ weight within a Benders decomposition is certainly interesting and would remove some of the redundant steps in the BBSA. Knowing the next $\lambda$ would enable one to simply iteratively solve
single-objective weighted-sum problems to obtain a complete set of efficient extreme solutions. From an application perspective, it could be interesting to apply the developed methodology to stochastic \glspl{BLP} with multiple subproblems. It would also be of interest to extend the proposed approach to multi-objective \glspl{LP} with three or more objectives. As long as an algorithm to solve such multi-objective \glspl{LP} is available, such as the multi-objective simplex algorithm \citep{ES73,RUV17} one should be able to develop a multi-objective version of Algorithm \ref{alg:BBSA} to iteratively identify efficient solutions in a multi-objective master problem, generate cuts in a multi-objective subproblem and, by doing so, iteratively identify a complete set of efficient solutions. These multi-objective simplex algorithms are more challenging to implement because efficient solutions of a multi-objective optimisation \gls{LP} have multiple neighbouring efficient solutions that all need to be explored.

\subsection*{Acknowledgement}
This research has been partially supported by the Marsden Fund, grant number 16-UOA-086. We thank the anonymous reviewers for constructive feedback that helped us improve this paper.

{ 

}

\newpage
\appendix
\section{Proof of Remark \ref{remark1}} \label{app:proof_remark1}

\begin{figure}{}
\begin{center}
\includegraphics[scale=0.8]{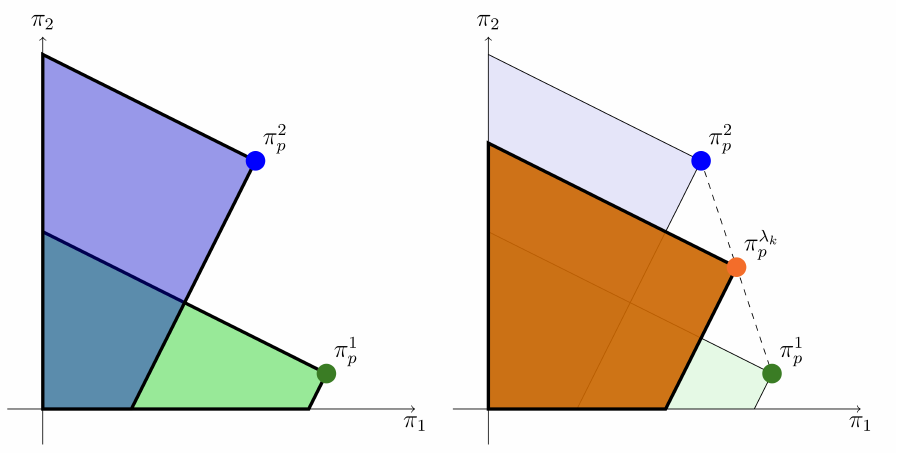}
\caption{A geometric illustration}
\label{fig:remark}
\end{center}
\end{figure}

\proof{}
Recall that for $\lambda_k\in(0,1)$ a weighted optimality cut \eqref{eq:weighted_optimality_cut} is of the form:
\begin{equation*}
    (b - By)^{\top} \pi_p^{\lambda_k} \leqq \lambda_k  \theta_1 + (1-\lambda_k)  \theta_2,
\end{equation*}
where $\pi_p^{\lambda_k}\in\mathcal{D}_p^{\lambda_k}$ is an extreme point in the feasible region of the dual formulation of $\lambda_k$-\gls{BS}, or $\lambda_k$-\eqref{eq:BSub}. If the two individual subproblem objectives are not contradictory, then for any master solution $\bar{y}$ that induces a feasible \gls{BS} \eqref{eq:BSub} there exists a solution $x^*$ that simultaneously minimises both objective functions. If $x^*$ respectively minimises each individual objective function, then it must also be an optimal solution to $\lambda_k$-\eqref{eq:BSub} for $\lambda_k\in(0,1)$. Let us assume that $\pi_p^1\in\D_p^1$ (resp. $\pi_p^2\in\D_p^0$) is an optimal dual solution for $x^*$ associated with the first (resp. second) objective function. Then $\pi_p^{\lambda_k}$ must be an optimal dual solution to $\lambda_k$-$\eqref{eq:BSub}$ for $\lambda_k\in(0,1)$, where
$$\pi_p^{\lambda_k} = \lambda_k \pi_p^1 + (1-\lambda_k) \pi_p^2.$$
Substituting this into~\eqref{eq:weighted_optimality_cut} gives the following 
$$\lambda_k(b - By)^{\top} \pi_p^1 +(1-\lambda_k) (b - By)^{\top}\pi_p^2\leqq \lambda_k \theta_1 + (1-\lambda_k) \theta_2.
$$
This implies that $\theta_1\geqq (b-By)^\top\pi_p^1$ and $\theta_2\geqq (b-By)^\top\pi_p^2$. Since the two single-objective subproblems agree for any choice of master solution $\bar{y}$, the weighted optimality cuts provide no additional information than the single-objective cuts defined by $\mathcal{D}_p^1$ and $\mathcal{D}_p^0$.  Consequently, they are not needed. Figure~\ref{fig:remark} provides a small example that illustrates this observation geometrically. The feasible region of the subproblem's dual formulation that is induced by a given master solution $\bar{y}$ is shown for $\lambda=1$ (in green) and for $\lambda=0$ (in blue). The orange feasible region corresponds to a $\lambda_k$-weighting of the two objective functions, where $\lambda_k\in(0,1)$. The optimal solution to the weighted problem is $\pi_p^{\lambda_k}$, and this is clearly a weighted combination of $\pi_p^1$ and $\pi_p^2$.
\endproof

\section{Solve Example \ref{ex:BLP-needs-weighted-cuts} with BBSA} \label{subsec:example}

To demonstrate how \gls{BBSA} can be employed to solve \glspl{BLP}, we apply it to Example \ref{ex:BLP-needs-weighted-cuts}.
We begin by decomposing~\eqref{eq:ex1-NP1} into a \gls{BM} and a \gls{BS} with decision variables $y=(y_{1},y_{2},y_{3})$ and $x=(x_1,x_2)$, respectively. The initial \gls{BM} is therefore:
\begin{equation}
\arraycolsep=1.8pt
\begin{array}{rrcrcrcrcrcrcrcrll}
\label{eq:ex1-BMASTER}
  \min & \left(\begin{array}{c} \theta_{1}\\ \phantom{.}\end{array}\right.&\left.\begin{array}{c} - \\ \phantom{.} \end{array}\right.&&&\left . \begin{array}{c} \\ \theta_{2}  \end{array} \right .&\left . \begin{array}{c}  \\ \end{array} \right .&&\left . \begin{array}{c} + \\ + \end{array} \right .&\left . \begin{array}{c} 2y_1 \\ 4y_1 \end{array} \right .&\left . \begin{array}{c}  \\ - \end{array} \right .&\left . \begin{array}{c}  \\ 6y_2 \end{array} \right .&\left . \begin{array}{c} - \\ - \end{array} \right .&\left . \begin{array}{c} 4 y_3 \\ 3y_3 \end{array} \right )&& \\
   \text{s.t.}& && && && & & y_{1} &-& 6y_{2} &-&4y_{3} &\geqq& -2,\\
  &&&&&&&&&y_{1},&&y_{2},&&y_{3},&\geqq& 0,
\end{array}
\end{equation}
while the \gls{BS} can be stated as:
\begin{equation}
\begin{array}{rrcrcrcrcrcrll}
\label{eq:ex1-BSUB}
  \min & \left ( \begin{array}{c} 4x_{1} \\ -2x_{1} \end{array} \right .&\left . \begin{array}{c} - \\ - \end{array} \right .&\left . \begin{array}{c} 2x_{2} \\ 2x_{2} \end{array} \right ) \\
   \text{s.t.}& - 2x_{1} & - & 6x_{2} &\geqq& -5 &+& 4\overline{y}_{1} &+& 3\overline{y}_{2}&+& 6\overline{y}_{3}  ,\\
  & -5x_{1} && &\geqq& -2 &+&  && 3\overline{y}_{2} &+& 5\overline{y}_{3},\\
  & x_{1}&,&x_{2}&\geqq& 0,
\end{array}
\end{equation}
where $\overline{y}=(\overline{y}_1,\overline{y}_2,\overline{y}_3)$ is the $y$-component of an efficient solution of~\eqref{eq:ex1-BMASTER}.
The \gls{BBSA} starts by finding the minimiser of the first objective function. To do this, we ignore the second objective function of~\eqref{eq:ex1-BMASTER} and~\eqref{eq:ex1-BSUB} and apply standard, single-objective Benders iterations (line \ref{alg:single_solve} Algorithm \ref{alg:BBSA}). Suppose that $\overline{y}^1=(\overline{y}^1_1,\overline{y}^1_2,\overline{y}^1_3)$ is the $y$-component of the solution that minimises the first objective function of~\eqref{eq:ex1-NP1}. \gls{BS} \eqref{eq:ex1-BSUB} is then solved with the bi-objective simplex algorithm given $\overline{y}^1$ (\texttt{Explore}$(\bar{y}^1)$ in line \ref{alg:Solve BS} Algorithm \ref{alg:BBSA}). In this example, two optimality cuts and one feasibility cut are generated
from these two steps and are added to the initial \gls{BM} \eqref{eq:ex1-BMASTER} (lines \ref{alg:single_solve} and \ref{alg:startcheck cuts}-\ref{alg:end check cuts} of Algorithm \ref{alg:BBSA}). The updated \gls{BM} then looks as follows:
\begin{equation}
\arraycolsep=1.8pt
\begin{array}{rrcrcrcrcrcrcrcrll}
\label{eq:ex1-updated_BMASTER_1}
      \min & \left ( \begin{array}{c} \theta_{1} \\ \phantom{.} \end{array} \right .&&&&\left . \begin{array}{c} \\ \theta_{2}  \end{array} \right .&&&\left . \begin{array}{c} + \\ + \end{array} \right .&\left . \begin{array}{c} 2y_1 \\ 4y_1 \end{array} \right .&\left . \begin{array}{c}  \\ - \end{array} \right .&\left . \begin{array}{c}  \\ 6y_2 \end{array} \right .&\left . \begin{array}{c} - \\ - \end{array} \right .&\left . \begin{array}{c} 4 y_3 \\ 3y_3 \end{array} \right )&& \\
   \text{s.t.}& && && && && y_{1} &-& 6y_{2} &-&4y_{3} &\geqq& -2,\\
  & && && && && && 3y_{2} &+& 5y_{3} &\leqq& 2,\\
  & -\theta_{1}&& && && &+& \frac{2}{3} y^{1} &+& \frac{1}{2} y^{2} &+& y_{3} &\leqq& \frac{5}{6},\\
  & && && -\theta_{2} &&&+& \frac{4}{3} y_{1} &+& \frac{9}{5} y_{1} &+& \frac{10}{3} y_{3} &\leqq& \frac{11}{5},\\
  &&&&&&&&&y_{1},&&y_{2},&&y_{3}&\geqq& 0.
\end{array}
\end{equation}
Next we identify the first efficient solution of \gls{BM} \eqref{eq:ex1-updated_BMASTER_1} in line \ref{alg:firstEff}, which becomes the last explored solution $(\theta^{last}_1, \theta^{last}_2, y^{last})$ with $z^{last} = (\frac{-61}{30},\frac{-31}{15})$. From here bi-objective simplex pivots are applied, and two consecutive non-dominated points $(\frac{-61}{30},\frac{-31}{15})$ and $(\frac{-16}{9},\frac{-23}{9})$ are found in this process (lines \ref{alg:startPivots}-\ref{alg:endPivots} of Algorithm \ref{alg:BBSA}). Figure~\ref{fig:ex1_iteration1} illustrates these points, connected by a solid blue line, as well as values of $\lambda$ for which the corresponding efficient solutions are optimal solutions of $\lambda$-\eqref{eq:ex1-updated_BMASTER_1}. While $\lambda$ is not explicitly derived in Algorithm \ref{alg:BBSA}, it can be computed as part of the simplex pivots in line \ref{alg:OnePivot}.
For reference, Figure~\ref{fig:ex1_iteration1} also indicates the set of non-dominated (extreme)  points of~\eqref{eq:ex1-NP1} connected by a dotted line. The figures also show the intervals of values of $\lambda$ for which the non-dominated extreme  points are optimal solutions of the corresponding weighted problem, where `?' indicates that the value is not known yet.

\begin{figure}[tb]
\centering
\begin{minipage}[t]{0.45\textwidth}
\scalebox{0.9}{
 \begin{tikzpicture}
    \begin{axis}[legend style= {nodes={scale=0.5, transform shape}},
        xlabel=$z_1$,
        ylabel=$z_2$,
        grid
      ]
      \addplot[color=black,mark=square,dashed,thick]
      coordinates{
        (-2.0333, -2.06667)
        (-1.7778, -2.5556)
        (-0.6667, -3.3333)
        (0.2, -3.6)
      };

      \addplot[->,color=blue, mark=x,line width=1.5pt]  coordinates {
        (-2.0333, -2.06667)
        (-1.7778, -2.5556)
        };
      \node (mark) [draw, green, circle, minimum size = 4pt, inner sep=2.5pt, thick]
      at (axis cs: -1.7778, -2.5556) {};
        \node (mark) [label=right:{$\left [\frac{44}{67},1 \right]$}]
      at (axis cs: -2.0333, -2.06667) {};
      \node (mark) [label=right:{$\left [?, \frac{44}{67} \right]$}]
      at (axis cs: -1.7778, -2.5556) {};

    \legend{non-dominated points of \eqref{eq:ex1-NP1}, non-dominated points of \eqref{eq:ex1-updated_BMASTER_1}}
    \end{axis}
  \end{tikzpicture}
  }
\caption{Iteration 1 and first two consecutive non-dominated points of \eqref{eq:ex1-updated_BMASTER_1}.}
\label{fig:ex1_iteration1}
\end{minipage}%
\begin{minipage}[t]{0.05\textwidth}
 $\quad$ %
 \end{minipage}%
\begin{minipage}[t]{0.45\textwidth}
\scalebox{0.9}{
 \begin{tikzpicture}
    \begin{axis}[legend style= {nodes={scale=0.5, transform shape}},
        xlabel=$z_1$,
        ylabel=$z_2$,
        grid
      ]
        \addplot[color=black,mark=square,dashed,thick]
      coordinates{
        (-2.0333, -2.06667)
        (-1.7778, -2.5556)
        (-0.6667, -3.3333)
        (0.2, -3.6)
      };
      \addplot[->,color=blue, mark=x,line width=1.5pt] coordinates {
        (-1.7778, -2.5556)
        (-0.67,-3.6)
        };

        \node (mark) [label=right:{$\left [\frac{47}{97}, \frac{44}{67} \right]$}]
      at (axis cs: -1.7778, -2.5556) {};
        \node (mark) [label=right:{$\left [?, \frac{47}{97} \right]$}]
      at (axis cs: -0.67,-3.6) {};
     \node (mark) [draw, green, circle, minimum size = 5pt, inner sep=2.5pt, thick]
      at (axis cs: -0.67, -3.6) {};
         \legend{non-dominated points of \eqref{eq:ex1-NP1}, non-dominated points of \eqref{eq:ex1-updated_BMASTER_1}}
    \end{axis}
  \end{tikzpicture}
  }
\caption{Two consecutive non-dominated points of problem \eqref{eq:ex1-updated_BMASTER_1} in iteration 2.}
\label{fig:ex1_iteration2}
\end{minipage}
\end{figure}

In the second iteration, the unexplored non-dominated point $\bar{z} = (\frac{-16}{9},\frac{-23}{9})$ (marked in green in Figure~\ref{fig:ex1_iteration1}) is explored.
The corresponding efficient solution is $(\overline{\theta}_{1},\overline{\theta}_{2},\overline{y})$.
We explore it by solving \gls{BS} \eqref{eq:ex1-BSUB}  (\texttt{Explore}($\bar{y}$) in line \ref{alg:Solve BS} of Algorithm \ref{alg:BBSA}). This time, however, only known cuts are generated, which confirms that $(\frac{-16}{9},\frac{-23}{9})$ is a non-dominated extreme  point of the original problem \eqref{eq:ex1-NP1}, and it becomes the last explored solution $z^{last}$.

To continue, bi-objective parametric simplex pivots are applied to $(\theta_1^{last},\theta_2^{last},y^{last})$ to identify the first unexplored non-dominated point with smaller second objective value (lines \ref{alg:startPivots}-\ref{alg:endPivots}) of Algorithm \ref{alg:BBSA}. We obtain the two non-dominated points $(\frac{-16}{9},\frac{-23}{9})$ and $(\frac{-2}{3},\frac{-18}{5})$, which are shown in blue in Figure~\ref{fig:ex1_iteration2}.
 The \gls{BBSA} sets the next solution to be explored to the one with objective vector $\bar{z}=(\frac{-2}{3},\frac{-18}{5})$ (line \ref{alg:updateNext}). This has a corresponding efficient solution that is optimal for $\lambda =  \frac{47}{97}$. When this solution is explored in the next iteration (line \ref{alg:Solve BS}), one cut is found and added to \gls{BM} (lines \ref{alg:startcheck cuts}-\ref{alg:end check cuts}), which gives:
\begin{equation}
\arraycolsep=1.8pt
\begin{array}{rrcrcrcrcrcrcrcrll}
\label{eq:ex1-updated_BMASTER_2}
        \min & \left ( \begin{array}{c} \theta_{1}^{+} \\ \phantom{.} \end{array} \right .&\left . \begin{array}{c} - \\ \phantom{.} \end{array} \right .&\left . \begin{array}{c} \theta_{1}^{-} \\ \phantom{.} \end{array} \right .&&\left . \begin{array}{c} \\ \theta_{2}^{+}  \end{array} \right .&\left . \begin{array}{c}  \\ - \end{array} \right .&\left . \begin{array}{c} \\ \theta_{2}^{-}  \end{array} \right .&\left . \begin{array}{c} + \\ + \end{array} \right .&\left . \begin{array}{c} 2y_1 \\ 4y_1 \end{array} \right .&\left . \begin{array}{c}  \\ - \end{array} \right .&\left . \begin{array}{c}  \\ 6y_2 \end{array} \right .&\left . \begin{array}{c} - \\ - \end{array} \right .&\left . \begin{array}{c} 4 y_3 \\ 3y_3 \end{array} \right )&& \\
   \text{s.t.}& && && && && y_{1} &-& 6y_{2} &-&4y_{3} &\geqq& -2,\\
  & && && && && && 3y_{2} &+& 5y_{3} &\leqq& 2,\\
  & -\theta_{1}^{+} &+& \theta_{1}^{-} && && &+& \frac{2}{3} y^{1} &+& \frac{1}{2} y^{2} &+& y_{3} &\leqq& \frac{5}{6},\\
  & && && -\theta^{+}_{2} &+& \theta^{-}_{2} &+& \frac{4}{3} y_{1} &+& \frac{9}{5} y_{1} &+& \frac{10}{3} y_{3} &\leqq& \frac{11}{5},\\
  & -\frac{4}{17}\theta_{1}^{+} &+&
  \frac{4}{17}\theta_{1}^{-} &-& \frac{13}{17}\theta_{2}^{+} &+& \frac{13}{17} \theta_{2}^{-} &+& \frac{20}{17} y_{1} &+& \frac{15}{17} y_{2} &+& \frac{30}{17} y_{3} &\leqq& \frac{25}{17},\\
  &\theta_{1}^{+}&,&\theta_{1}^{-}&,&\theta_{2}^{+}&,&\theta_{2}^{-}&,&y_{1}&,&y_{2}&,&y_{3} &\geqq& 0.
\end{array}
\end{equation}
The \gls{BBSA} then moves on to the unexplored point with objective vector $\bar{z}=(\frac{-2}{3},\frac{-10}{3})$ with associated weight $\lambda = \frac{7}{17}$ (line \ref{alg:updateNext}), which is explored in line \ref{alg:Solve BS} in the next iteration of the \gls{BBSA} Algorithm \ref{alg:BBSA}.
This iterative process continues until there are no remaining unexplored non-dominated points. Figure~\ref{fig:ex1-final points} shows the final set of non-dominated extreme  points that  is found by the BBSA,  and we can see that the reformulation correctly identifies all extreme
non-dominated points of the non-decomposed (original) problem.

\begin{figure}[tb]
\centering
\begin{minipage}[t]{0.45\linewidth}
\scalebox{0.9}{
 \begin{tikzpicture}
    \begin{axis}[legend style= {nodes={scale=0.5, transform shape}},
        xlabel=$z_1$,
        ylabel=$z_2$,
        grid
      ]
      \addplot[color=black,mark=square,dashed, thick]
      coordinates{
        (-2.0333, -2.06667)
        (-1.7778, -2.5556)
        (-0.6667, -3.3333)
        (0.2, -3.6)
      };
      \addplot[color=blue,mark=x] coordinates {
        (-1.7778, -2.5556)
        (-0.6667, -3.3333)
      };
      \addplot[color=pink,mark=x] coordinates {
        (-1.7778, -2.5556)
        (-0.6667, -3.6)
        };
        \draw [->,blue,thick,line width=2pt] (-0.6667, -3.6) -- (-1.7778, -2.5556) ;
      \draw [->,blue,thick,line width=2pt] (-1.7778, -2.5556) -- (-0.6667, -3.3333) ;
        \node (mark) [draw, green, circle, minimum size = 5pt, inner sep=2pt, thick]
      at (axis cs: -0.67, -3.333) {};

      \node (mark) [label=right:{$\left [\frac{7}{17}, \frac{44}{67} \right]$}]
      at (axis cs: -1.7778, -2.5556) {};
        \node (mark) [label=right:{$\left [?, \frac{7}{17} \right]$}]
      at (axis cs: -0.6667, -3.3333) {};

    \legend{non-dominated points of \eqref{eq:ex1-NP1},non-dominated points of \eqref{eq:ex1-updated_BMASTER_2}}
    \fill[black!40!white] (-1.7778, -2.5556) -- (-0.6667, -3.6) -- (-0.6667, -3.3333);

    \end{axis}
  \end{tikzpicture}
  }
\caption{Two consecutive non-dominated points of problem \eqref{eq:ex1-updated_BMASTER_2} in iteration 3.}
\label{fig:ex1-iteration3}
\end{minipage}
\begin{minipage}[t]{0.05\textwidth}
$\quad$
 \end{minipage}
\begin{minipage}[t]{0.45\linewidth}
\scalebox{0.9}{
 \begin{tikzpicture}
    \begin{axis}[legend style= {nodes={scale=0.5, transform shape}},
        xlabel=$z_1$,
        ylabel=$z_2$,
        grid
      ]
      \addplot[color=black,,mark=square,dashed, thick] coordinates {
        (-2.0333, -2.06667)
        (-1.7778, -2.5556)
        (-0.6667, -3.3333)
        (0.2, -3.6)
        };
      \addplot[color=blue,mark=x] coordinates {
        (-2.0333, -2.06667)
        (-1.7778, -2.5556)
        (-0.6667, -3.3333)
        (0.2, -3.6)
        };
        \node (mark) [label=right:{$\left [\frac{44}{67},1 \right]$}]
      at (axis cs: -2.0333, -2.06667) {};
      \node (mark) [label=right:{$\left [\frac{7}{17}, \frac{44}{67} \right]$}]
      at (axis cs: -1.7778, -2.5556) {};
       \node (mark) [label=left:{$\left [\frac{4}{17}, \frac{7}{17}\right]$}]
      at (axis cs: -0.6667, -3.3333) {};
      \node (mark) [label=left:{$\left [ 0, \frac{4}{17} \right]$}]
      at (axis cs: 0.2, -3.6) {};

         \legend{non-dominated points of \eqref{eq:ex1-NP1}, \gls{BBSA} non-dominated points \eqref{eq:ex1-updated_BMASTER_2}, font=\tiny}
    \end{axis}
  \end{tikzpicture}
  }
\caption{Final non-dominated extreme points found by the \gls{BBSA} for \eqref{eq:ex1-NP1}.}
\label{fig:ex1-final points}
\end{minipage}
\end{figure}

\section{Illustration of cases for proof of Proposition \ref{prop:allfound}} \label{sec:Appendix_Prop2}
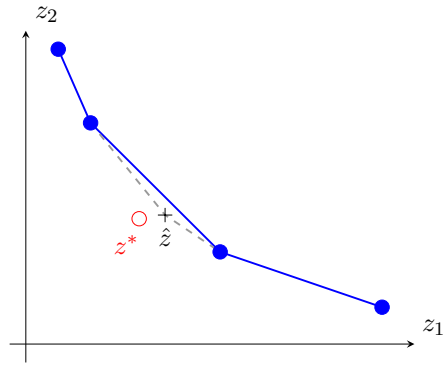
\begin{figure}[h!]
    \centering
\scalebox{0.9}{
  \begin{tikzpicture}
  \pgfplotsset{ticks=none}
      \begin{axis}[grid=none,ymin=-0.5,ymax=8.5,xmax=12,xmin=-0.5,xticklabel=\empty,yticklabel=\empty,minor tick num=0,axis lines = middle,xlabel=$z_1$,ylabel=$z_2$,label style =
               {at={(ticklabel cs:1.1)}}]
      \addplot[mark=*, only marks, color=blue,mark size=3pt] coordinates {
        (1, 8)
        (2, 6)
        (6, 2.5)
        (11, 1)
        (4.3, 10)
      };
      \addplot[mark=none, thick, color=blue, mark size=3pt] coordinates {
        (1, 8)
        (2, 6)
        (6, 2.5)
        (11, 1)
      };
      \addplot[mark=+, only marks, color=black,mark size=3pt] coordinates {
        (4.3, 3.5)
        (4.3,9)
      };
      \addplot[mark=none, thick, color=black!40, dashed, mark size=3pt] coordinates {
        (2, 6)
        (4.3, 3.5)
        (6, 2.5)
      };
      \addplot[mark=o, only marks, color=red,mark size=3pt] coordinates {
        (3.5, 3.4)
      };
    
    \node at (axis cs: 7,10) {non-dominated};
    \node at (axis cs: 6.15,9) {missing $\hat{z}$};
    \node at (axis cs: 4.3, 2.9) {$\hat{z}$};
    \node at (axis cs: 3.1, 2.7) {\textcolor{red}{$z^*$}};
    
    \end{axis}
  \end{tikzpicture}
  }
    \caption{Cases for ``missing'' a non-dominated extreme point.}
    \label{fig:Prop2_cases}
\end{figure}

\section{Problem instance formulations}
\subsection{\gls{FCTP}}
\label{sec:Appendix_FCTP}
Parameters for an \gls{FCTP} with $m$ sources and $n$ sinks are:
\begin{itemize}
    \item demand $d_{j}$ of each sink $j=1,2,\ldots,m$
    \item $s_{i}$ of each source  $i=1,2,\ldots n$
    \item per-unit cost $c_{ij}^k$ and fixed charge $f_{ij}^k$ for objectives $k=1,2$.
    \item $A$ is the set of all source-sink connections
\end{itemize}
A mathematical formulation of the \gls{BFCTP} is:
  \begin{equation}
\label{eq:bi-objective-fixed-charge}
\begin{array}{rrcrcrr}
   \min&\left ( \begin{array}{c} \sum\limits_{(i,j) \in A} c^{1}_{ij} x_{ij} \\ \sum\limits_{(i,j) \in A } c^{2}_{ij} x_{ij} \end{array} \right .&\left.\begin{array}{c} \underset{\phantom{A}}{+} \\  \underset{\phantom{A}}{+} \end{array} \right. &\left . \begin{array}{c} \sum\limits_{(i,j) \in A} f^{1}_{ij} y_{ij} \\  \sum\limits_{(i,j) \in A} f^{2}_{ij} y_{ij} \end{array} \right ) &&&\\
   \text{s.t.}& & & y_{ij} & \leqq & 1 & \; \forall (i,j) \in A,\\
   &\sum\limits_{j=1}^{m}x_{ij} & & & \leqq & s_{i} & \; \forall i=1\ldots,n,\\
    &\sum\limits_{i=1}^{n}x_{ij} & & & \geqq & d_{j} & \; \forall j=1\ldots,m,\\
    & x_{ij} & - & m_{ij}y_{ij} & \leqq & 0 & \; \forall (i,j) \in A,\\
    & x_{ij}, &  & y_{ij} & \geqq& 0 &
   \; \forall (i,j) \in A,
\end{array}
\end{equation}
 where $m_{ij} = \min \{ s_i, d_j\} \, \forall (i,j) \in A$. The problem is decomposed as in Section \ref{sec:full-Benders-reformulation-BLP} with variables $y$ in \gls{BM} and variables $x$ in \gls{BS}.

\subsection{\gls{BPP}}
\label{sec:Appendix_BPP}
The mathematical formulation of \gls{BPP} is as follows, adapted from \gls{PP} in \citet{qiu2014}:

\begin{equation}
\label{eq:bi-objective-portfolio}
\begin{array}{rrcrcrr}
   \min&\left ( \begin{array}{c} (c^{1})^{\top} x \\ (c^{2})^{\top} x \end{array} \right .&\left . \begin{array}{c} + \\ + \end{array} \right .&\left . \begin{array}{c} (f^{1})^{\top} y \\ (f^{1})^{\top} y \end{array} \right )&&&\\
   \text{s.t.}& & & e^{\top} y  & = & 1, & \;\\
   & r x_{i}& + & a_{i}^{\top} y & \geqq & r & \; \forall i=1\ldots,m,\\
    & \sum_{i=1}^{m} x_{i}& & & \leqq & k,& \;\\
    & x_{i}& & & \leqq & 1& \; \forall i=1\ldots,m,\\
    & x_{i}, &  & y_{j} & \geqq& 0 &

   \; \forall i = 1,\ldots m, j=1,\ldots,n.
\end{array}
\end{equation}

In \eqref{eq:bi-objective-portfolio}, variable $y_{j}$ corresponds to the amount of investment of the $j$-th asset ($j=1,2,\ldots,n$), and $x_{i}$ captures by how much the $i$-th return constraint is violated ($i=1,2,\ldots,m$). The first constraint is a budget constraint obtained by scaling the investment level to a unit budget. The second set of constraints states the overall return should be at least $r$, where $a_{i}$ corresponds to the return vector of the $n$ assets for the $i$-th scenario. The third set of constraints imposes that the total amount by which the return constraints are violated should not exceed a specific number $k$. We introduce a penalty when a return constraint is violated by adding a cost corresponding to the violated constraints in the objective functions. 
To apply the \gls{BBSA} to \gls{BPP}, we consider variables $y$ in the \gls{BM} and $x$ in the \gls{BS}. The original instances are available in \citet{Ahmed2015}.

\subsection{\gls{BMKP}}
\label{sec:Appendix_BMKP}

In the \gls{MKP}, there are three sets of items (called $I, J$, and $K$) can be packed into two sets of knapsacks (called $L$, and $M$). The problem is how to pack these items in two knapsacks to minimise the overall cost of the  items. The mathematical formulation of \gls{BMKP} is adapted from \citet{maher2021}, where a stochastic version of the problem is presented.
\begin{equation}
\label{eq:bi-objective-BMKP}
\begin{array}{rrcrcrcrcll}
   \min&\left ( \begin{array}{c} \sum\limits_{k \in K} c^{1}_{k} x_{k} \\ \sum\limits_{k \in K} c^{2}_{k} x_{k} \end{array} \right .&\left.\begin{array}{c} \underset{\phantom{K}}{+} \\  \underset{\phantom{K}}{+} \end{array} \right. &\left . \begin{array}{c} \sum\limits_{i \in I} f^{1}_{i} y_{i} \\  \sum\limits_{i \in I} f^{2}_{i} y_{i} \end{array} \right .\left.\begin{array}{c} \underset{\phantom{K}}{+} \\  \underset{\phantom{K}}{+} \end{array} \right. &\left . \begin{array}{c} \sum\limits_{j \in J} q^{1}_{j} z_{j} \\  \sum\limits_{j \in J} q^{2}_{j} z_{j} \end{array} \right ) &&&&\\
   \text{s.t.}& & & \sum\limits_{i \in I} a_{il} y_{i} &+& \sum\limits_{j \in J} b_{jl} z_{j}&\geqq &\gamma_{l} &\forall l \in L, \;\\
              & \sum\limits_{k \in K} f_{km}x_{k} &+& \sum\limits_{i \in I}e_{im}y_{i}&&&\geqq& \delta_{m}& \forall m \in M, \\
              & x_{k} &&&&&\leqq& 1& \forall k \in K,\\
              &&&y_{i}&&&\leqq& 1& \forall i \in I, \\
              &&&&&z_{j}&\leqq& 1& \forall j \in J,\\
              & x_{k},&&y_{i},&&z_{j}&\geqq& 0&\forall k \in K, \forall i \in I, \forall j \in J.
\end{array}
\end{equation}

The first set of constraints in \eqref{eq:bi-objective-BMKP} states that the total weight of items of sets $I$ and $J$ should be at least $\gamma_{l}$ in each knapsack $l \in L$. Similarly, the second set of constraints ensures that the total weight of items of sets $I$ and $K$ is at least $\delta_{m}$ in each knapsack $m \in M$. To apply the proposed \gls{BBSA}, we decompose \eqref{eq:bi-objective-BMKP} into \gls{BM} including $y$ and $z$ variables and \gls{BS} containing $x$ variables. In \citet{angulo2016}, $30$ instances of \eqref{eq:bi-objective-BMKP} with $|I| = |J| = |K| =120$ items of each type are generated, where \gls{BM} and \gls{BS} contain $|L| = 50$ and $|M| = 5$ constraints, respectively. 

 \subsection{Full numerical results of the \gls{BBSA} for \gls{BMKP} instances}
\label{sec:BMKP_Table}
\begin{table}[H]
\centering
\caption{Summary of run statistics for \gls{BMKP} instances.} \label{tab:MKP}.
\begin{tabular}{lrrrrrrrrrr}
\hline
$Instances$ & $|\mathcal{Z}_n|$ & FC & OC & AC & it & tBBSA & Ben\% & BS\% & BM\% & DBA \\ \hline
\gls{BMKP}\_1  & 653 & 0 & 66274 & 23272 & 924 & 75.08 & 0.24 & 45.90 & 51.24 & 2.80 \\
\gls{BMKP}\_2  & 613 & 0 & 61841 & 21267 & 939 & 75.02 & 0.23 & 43.92 & 53.16 & 3.18 \\
\gls{BMKP}\_3  & 575 & 0 & 57190 & 20137 & 808 & 65.79 & 0.26 & 46.37 & 50.78 & 2.98 \\
\gls{BMKP}\_4  & 573 & 0 & 60805 & 18654 & 733 & 61.57 & 0.41 & 50.15 & 46.87 & 3.64 \\
\gls{BMKP}\_5  & 573 & 0 & 67899 & 19500 & 756 & 67.55 & 0.40 & 51.61 & 45.49 & 3.17 \\
\gls{BMKP}\_6  & 595 & 0 & 78900 & 27152 & 825 & 78.05 & 0.22 & 49.85 & 47.48 & 3.10 \\
\gls{BMKP}\_7  & 571 & 0 & 70878 & 18474 & 870 & 73.25 & 0.41 & 46.29 & 47.78 & 2.87 \\
\gls{BMKP}\_8  & 538 & 0 & 79976 & 19433 & 931 & 82.33 & 0.47 & 49.59 & 47.36 & 2.93 \\
\gls{BMKP}\_9  & 605 & 0 & 73491 & 26256 & 889 & 77.82 & 0.46 & 49.00 & 48.01 & 3.89 \\
\gls{BMKP}\_10 & 559 & 0 & 81766 & 23889 & 899 & 83.50 & 0.26 & 48.99 & 48.31 & 2.30 \\
\gls{BMKP}\_11 & 586 & 0 & 61490 & 23530 & 836 & 67.72 & 0.27 & 47.67 & 49.45 & 3.30 \\
\gls{BMKP}\_12 & 496 & 0 & 76459 & 21538 & 931 & 82.78 & 0.31 & 47.67 & 49.50 & 1.96 \\
\gls{BMKP}\_13 & 475 & 0 & 54509 & 14928 & 777 & 62.60 & 0.56 & 46.12 & 50.65 & 2.72 \\
\gls{BMKP}\_14 & 569 & 0 & 71771 & 19249 & 941 & 82.25 & 0.41 & 46.08 & 50.98 & 2.18 \\
\gls{BMKP}\_15 & 574 & 0 & 72548 & 23660 & 821 & 74.16 & 0.23 & 48.77 & 48.49 & 3.15 \\
\gls{BMKP}\_16 & 586 & 0 & 56416 & 21052 & 844 & 64.99 & 0.28 & 46.62 & 50.39 & 3.85 \\
\gls{BMKP}\_17 & 544 & 0 & 74869 & 22598 & 864 & 76.40 & 0.22 & 48.53 & 48.76 & 4.06 \\
\gls{BMKP}\_18 & 585 & 0 & 67389 & 15572 & 818 & 71.65 & 0.40 & 48.08 & 48.96 & 2.85 \\
\gls{BMKP}\_19 & 589 & 0 & 60975 & 22024 & 824 & 69.03 & 0.25 & 46.02 & 51.09 & 2.63 \\
\gls{BMKP}\_20 & 524 & 0 & 61316 & 17303 & 830 & 68.17 & 0.25 & 46.44 & 50.70 & 3.78 \\
\gls{BMKP}\_21 & 624 & 0 & 63528 & 25339 & 954 & 76.64 & 0.26 & 44.95 & 52.14 & 5.06 \\
\gls{BMKP}\_22 & 571 & 0 & 64139 & 21679 & 807 & 69.05 & 0.25 & 47.39 & 49.82 & 3.07 \\
\gls{BMKP}\_23 & 629 & 0 & 56437 & 20994 & 834 & 66.80 & 0.27 & 44.39 & 52.65 & 3.07 \\
\gls{BMKP}\_24 & 641 & 0 & 57535 & 13827 & 867 & 66.85 & 0.40 & 45.85 & 51.02 & 3.20 \\
\gls{BMKP}\_25 & 509 & 0 & 70925 & 17502 & 865 & 73.09 & 0.34 & 49.54 & 47.50 & 3.05 \\
\gls{BMKP}\_26 & 600 & 0 & 74484 & 23568 & 958 & 81.67 & 0.36 & 48.88 & 48.18 & 3.04 \\
\gls{BMKP}\_27 & 671 & 0 & 78569 & 19640 & 883 & 83.39 & 0.41 & 47.22 & 49.95 & 3.01 \\
\gls{BMKP}\_28 & 587 & 0 & 73518 & 17496 & 850 & 76.40 & 0.42 & 48.34 & 48.69 & 2.67 \\
\gls{BMKP}\_29 & 705 & 0 & 67929 & 21888 & 852 & 74.89 & 0.24 & 46.62 & 50.59 & 3.66 \\
\gls{BMKP}\_30 & 624 & 0 & 61340 & 18302 & 835 & 69.22 & 0.25 & 46.01 & 51.11 & 3.49
\\\hline
\end{tabular}
\end{table}

\end{document}